\documentclass[10pt]{amsart}
\usepackage{amssymb}\usepackage{euscript}
\usepackage[all]{xy}

\newtheorem{prop}{Proposition}[section]
\newtheorem{prop:def}{Proposition-Definition}[section]

\newtheorem{lemma}{Lemma}[section]
\newtheorem{fact}{Fact}[section]
\newtheorem{thm}{Theorem}[section]
\newtheorem{cor}{Corollary}[section]
\theoremstyle{remark}

\begin{document}

\newcommand{\nc}{\newcommand} \nc{\on}{\operatorname}

\nc{\pa}{\partial}

\nc{\cA}{{\cal A}} \nc{\cB}{{\cal B}}\nc{\cC}{{\cal C}}
\nc{\cD}{{\cal D}}
\nc{\cE}{{\cal E}} \nc{\cG}{{\cal G}}\nc{\cH}{{\cal H}}
\nc{\cI}{{\cal I}} \nc{\cJ}{{\cal J}}\nc{\cK}{{\cal K}}
\nc{\cL}{{\cal L}} \nc{\cR}{{\cal R}} \nc{\cS}{\mathcal{S}}
\nc{\cV}{{\cal V}} \nc{\cX}{{\cal X}}

\nc{\¦}{{|}}

\nc{\sh}{\on{sh}}\nc{\Id}{\on{Id}}\nc{\Diff}{\on{Diff}}
\nc{\Perm}{\on{Perm}}\nc{\conc}{\on{conc}}\nc{\Alt}{\on{Alt}}
\nc{\ad}{\on{ad}}\nc{\Der}{\on{Der}}\nc{\End}{\on{End}}
\nc{\no}{\on{no\ }} \nc{\res}{\on{res}}\nc{\ddiv}{\on{div}}
\nc{\Sh}{\on{Sh}} \nc{\card}{\on{card}}\nc{\dimm}{\on{dim}}
\nc{\Sym}{\on{Sym}} \nc{\Jac}{\on{Jac}}\nc{\Ker}{\on{Ker}}
\nc{\Spec}{\on{Spec}}\nc{\Cl}{\on{Cl}}
\nc{\Imm}{\on{Im}}\nc{\limm}{\lim}\nc{\Ad}{\on{Ad}}
\nc{\ev}{\on{ev}} \nc{\Hol}{\on{Hol}}\nc{\Det}{\on{Det}}
\nc{\Bun}{\on{Bun}}\nc{\diag}{\on{diag}}\nc{\pr}{\on{pr}}
\nc{\Span}{\on{Span}}\nc{\Comp}{\on{Comp}}\nc{\Part}{\on{Part}}
\nc{\tensor}{\on{tensor}}\nc{\ind}{\on{ind}}\nc{\id}{\on{id}}
\nc{\Hom}{\on{Hom}}\nc{\Quant}{\on{Quant}} \nc{\Dequant}{\on{Dequant}}\nc{\Def}{\on{Def}}
\nc{\AutLBA}{\on{AutLBA}}\nc{\AutQUE}{\on{AutQUE}}
\nc{\LBA}{{\on{LBA}}}\nc{\Aut}{\on{Aut}}\nc{\QUE}{{\on{QUE}}}
\nc{\Lyn}{\on{Lyn}}\nc{\Cof}{\on{Cof}}\nc{\LCA}{\underline{\on{LCA}}}\nc{\FLBA}{\on{FLBA}}
\nc{\LA}{\on{LA}}\nc{\FLA}{\on{FLA}}\nc{\EK}{\on{EK}}
\nc{\class}{\on{class}}\nc{\br}{\on{br}}\nc{\co}{\on{co}}
\nc{\Prim}{\on{Prim}}\nc{\ren}{\on{ren}}\nc{\lbr}{\on{lbr}}
\nc{\SC}{\on{SC}}\nc{\ol}{\overline}\nc{\FL}{\on{FL}}
\nc{\FA}{\on{FA}}\nc{\alg}{{\on{alg}}}\nc{\KZ}{\on{KZ}}
\nc{\op}{{\on{op}}}\nc{\cop}{{\on{cop}}}
\nc{\Inn}{\on{Inn}}\nc{\OutDer}{\on{OutDer}}
\nc{\inv}{{\on{inv}}}\nc{\gr}{{\on{gr}}}\nc{\Lie}{{\on{Lie}}}
\nc{\Out}{{\on{Out}}}\nc{\univ}{{\on{univ}}}
\nc{\GT}{{\on{GT}}}\nc{\GRT}{{\on{GRT}}}
\nc{\GS}{{\on{GS}}} 
\nc{\restr}{{\on{restr}}}\nc{\lin}{{\on{lin}}}
\nc{\mult}{{\on{mult}}}\nc{\qt}{{\on{qt}}}
\nc{\compl}{{\on{compl}}} \nc{\Trees}{{\on{Trees}}}
\nc{\Ord}{{\on{Ord}}} \nc{\norm}{{\on{norm\ ord}}}
\nc{\fd}{{\on{fd}}}\nc{\Maps}{{\on{Maps}}}
\nc{\CYBE}{{\on{CYBE}}} \nc{\CY}{{\on{C}}} \nc{\can}{{\on{can}}}
\nc{\Alg}{{\underline{\on{Alg}}}} \nc{\Poisson}{{\underline{\on{Poisson}}}}
\nc{\Coalg}{{\underline{\on{Coalg}}}} \nc{\Vect}{{\underline{\on{Vect}}}}
\nc{\UE}{{\on{UE}}}\nc{\Bialg}{{\underline{\on{Bialg}}}}
\nc{\QTA}{{\underline{\on{QTA}}}}\nc{\QTQUE}{{\underline{\on{QTQUE}}}}
\nc{\QTLBA}{{\underline{\on{QTLBA}}}}\nc{\TA}{{\underline{\on{TA}}}}
\nc{\MT}{{\underline{\on{MT}}}}\nc{\QYBE}{{\underline{\on{QYBE}}}}
\nc{\QTBialg}{{\underline{\on{QTBialg}}}}
\nc{\qcocomm}{{\on{quasi-cocomm}}}
\nc{\qcomm}{{\on{quasi-comm}}}\nc{\cocomm}{{\on{cocomm}}}
\nc{\dual}{{\on{dual}}} \nc{\sym}{{\on{sym}}}
\nc{\cycl}{{\on{cycl}}} \nc{\Prop}{{\underline{\on{Prop}}}}
\nc{\uV}{{\underline{V}}} \nc{\uR}{{\underline{R}}}
\nc{\graphs}{{\on{graphs}}} \nc{\Isom}{{\on{Isom}}}
\nc{\cp}{{\on{co-Poisson}}}\nc{\nat}{{\on{nat}}}\nc{\loc}{{\on{loc}}}
\nc{\PolBMC}{{\on{PolBMC}}}\nc{\coPoisson}{{\on{co-Poisson}}}
\nc{\as}{{\on{antisymm}}}\nc{\Harrison}{{\on{Harrison}}}
\nc{\Mor}{{\on{Mor}}}\nc{\Ob}{{\on{Ob}}}
\nc{\Bil}{{\on{Bil}}} \nc{\invt}{{\on{invt}}}

\nc{\al}{\alpha}\nc{\de}{\delta}
\nc{\eps}{\epsilon}\nc{\la}{{\lambda}}
\nc{\si}{\sigma}\nc{\z}{\zeta}

\nc{\La}{\Lambda}

\nc{\ve}{\varepsilon} \nc{\vp}{\varphi}

\nc{\AAA}{{\mathbb A}}\nc{\BB}{{\mathbb B}}
\nc{\CC}{{\mathbb C}}\nc{\ZZ}{{\mathbb Z}}
\nc{\QQ}{{\mathbb Q}} \nc{\NN}{{\mathbb N}}\nc{\VV}{{\mathbb V}}
\nc{\KK}{{\mathbb K}}

\nc{\ff}{{\mathbf f}}\nc{\bg}{{\mathbf g}}
\nc{\ii}{{\mathbf i}}\nc{\kk}{{\mathbf k}}
\nc{\bl}{{\mathbf l}}\nc{\zz}{{\mathbf z}}
\nc{\pp}{{\mathbf p}}\nc{\qq}{{\mathbf q}}

\nc{\cF}{{\cal F}}\nc{\cM}{{\cal M}}\nc{\cO}{{\cal O}}
\nc{\cT}{{\cal T}}\nc{\cW}{{\cal W}}

\nc{\Assoc}{{\mathbf Assoc}}
\nc{\ul}{\underline}
\def\Sha{{\mathop{\scriptstyle\amalg\!\hspace{-1.8pt}\amalg}}}

\nc{\ub}{{\underline{b}}}
\nc{\uk}{{\underline{k}}}
\nc{\un}{{\underline{n}}} \nc{\um}{{\underline{m}}}
\nc{\up}{{\underline{p}}}\nc{\uq}{{\underline{q}}}
\nc{\ur}{{\underline{r}}}
\nc{\us}{{\underline{s}}}\nc{\ut}{{\underline{t}}}
\nc{\uw}{{\underline{w}}}
\nc{\uz}{{\underline{z}}}
\nc{\ual}{{\underline{\alpha}}}\nc{\ualpha}{{\underline{\alpha}}}
\nc{\ubeta}{{\underline{\beta}}}\nc{\ugamma}{{\underline{\gamma}}}
\nc{\ueps}{{\underline{\epsilon}}}\nc{\ueta}{{\underline{\eta}}}
\nc{\uzeta}{{\underline{\zeta}}}\nc{\ula}{{\underline{\lambda}}}
\nc{\umu}{{\underline{\mu}}}\nc{\unu}{{\underline{\nu}}}
\nc{\usigma}{{\underline{\sigma}}}\nc{\utau}{{\underline{\tau}}}
\nc{\uI}{{\underline{I}}}\nc{\uJ}{{\underline{J}}}
\nc{\uK}{{\underline{K}}}\nc{\uM}{{\underline{M}}}
\nc{\uN}{{\underline{N}}}

\nc{\A}{{\mathfrak a}}
\renewcommand{\t}{{\mathfrak t}}\nc{\G}{{\mathfrak g}}
\nc{\x}{{\mathfrak x}}
\nc{\B}{{\mathfrak b}} \nc{\C}{{\mathfrak c}}
\nc{\D}{{\mathfrak d}} \nc{\HH}{{\mathfrak h}}
\nc{\iii}{{\mathfrak i}}\nc{\mm}{{\mathfrak m}}\nc{\N}{{\mathfrak n}}
\nc{\ttt}{{\mathfrak{t}}}\nc{\U}{{\mathfrak u}}\nc{\V}{{\mathfrak v}}
\nc{\grt}{{\mathfrak grt}}\nc{\gt}{{\mathfrak gt}}
\nc{\SL}{{\mathfrak{sl}}}\nc{\out}{{\mathfrak{out}}}

\nc{\SG}{{\mathfrak S}}

\nc{\wt}{\widetilde} \nc{\wh}{\widehat}

\nc{\bn}{\begin{equation}}\nc{\en}{\end{equation}} \nc{\td}{\tilde}

\newcommand{\g}{\mathfrak{g}}
\newcommand{\ghat}{\tilde{\mathfrak{g}}}
\newcommand{\h}{\mathfrak{h}}
\newcommand{\nilp}{\mathfrak{n}}
\newcommand{\borel}{{\mathfrak{b}}}
\newcommand{\borelp}{{\mathfrak{b}_{+}}}
\newcommand{\borelm}{{\mathfrak{b}_{-}}}
\newcommand{\borelhat}{\tilde{ \mathfrak{b}}}
\newcommand{\p}[1]{\bar {#1}}
\newcommand{\Lba}{\text{\underline{LBA}}}
\newcommand{\EKh}{\mathcal{U}_{\hbar}^{EK}}
\newcommand{\EKm}{\mathcal{U}_{-\hbar}^{EK}}
\newcommand{\DJh}{\mathcal{U}_{\hbar}^{DJ}}
\newcommand{\Uhbhat}{{\EKh(\borelhat_{+}) } }
\newcommand{\Uhb}{{\EKh(\borelp) } }
\newcommand{\Uhatp}{\mathcal{\widetilde{U}_{+}}}
\newcommand{\Up}{\mathcal{U_{+}}}
\newcommand{\Upqd}{\mathcal{U_{+}^{*}}}
\newcommand{\E}{\mathbb{E}}
\newcommand{\F}{\mathbb{F}}
\newcommand{    \qd}{{}^{\star}}
\newcommand{    \qdo}{{}^{\star cop}}
\newcommand{\cato}{\mathcal{O}}

\title[Compatibility of quantization functors
with duality and doubling operations]{Compatibility of quantization
functors of Lie bialgebras with duality and doubling operations}

\begin{abstract}
We study the behavior of the Etingof-Kazhdan quantization functors
under the natural duality operations of Lie bialgebras and Hopf algebras.
In particular, we prove that these functors are "compatible with duality",
i.e., they commute with the operation of duality followed by replacing the
coproduct by its opposite. We then show that any quantization functor
with this property also commutes with the operation of taking doubles.
As an application, we show that the Etingof-Kazhdan quantization of
some affine Lie superalgebras coincide with their Drinfeld-Jimbo-type
quantizations.
\end{abstract}

\author{Benjamin Enriquez}
\address{B.E.: IRMA (CNRS), rue Ren\'e Descartes, F-67084 Strasbourg, France}
\email{enriquez@math.u-strasbg.fr}

\author{Nathan Geer}
\address{N.G.: School of Mathematics, Georgia Institute of Technology, Atlanta, GA
30332-0160, USA}
\email{geer@math.gatech.edu}

\maketitle

We fix a field $\kk$ of characteristic $0$. Unless specified otherwise,
``algebra'', ``vector space'', etc.,  means ``algebra over $\kk$", etc.

\section*{Introduction}

In \cite{EK1,EK2}, Etingof and Kazhdan solved the problem of
quantization of Lie bialgebras. For each Drinfeld associator
$\Phi$, they constructed a quantization functor $\tilde Q_{\Phi} :
\{$Lie bialgebras$\} \to \{$quantized universal enveloping (QUE)
algebras$\}$, right inverse to the semiclassical limit functor.
This functor can in fact be viewed as a morphism of props
$Q_{\Phi} : \on{Bialg} \to S({\bf LBA})$. It was proved in
\cite{EK1} that $\tilde
Q_{\Phi}$ is compatible with natural operations on finite
dimensional Lie bialgebras and QUE algebras. Namely, if $\A$ is a
finite dimensional Lie bialgebra, then $\tilde Q_{\Phi}
(\A^{*cop})\simeq \tilde Q_{\Phi}(\A)^{\star cop}$ and $\tilde
Q_{\Phi}({\mathfrak D}(\A)) \simeq D(\tilde Q_{\Phi}(\A))$. Here
for $U$ a QUE algebra with finite-dimensional associated Lie bialgebra, 
$U^{\star}$ denotes the QUE algebra associated
to the quantum formal series Hopf algebra
$\on{Hom}(U,\kk[[\hbar]])$ (see \cite{Dr,Gav}); $cop$ indicates
the opposite cobracket or coproduct; ${\mathfrak D}(\A)$ is the
double of $\A$ and $D(U)$ is the double QUE algebra of $U$.

In this paper, we prove that these isomorphisms also hold when $\A$ is 
infinite-dimensional. This will be proved as a consequence of 
statements about the prop morphism $Q_{\Phi} : \on{Bialg} \to S({\bf LBA})$.
One can explain the difference of our work with the relevant part of 
\cite{EK1} as follows: the statements of \cite{EK1} can be made into 
``propic'' statements, involving the morphism $\on{Bialg} \to S({\bf LBA}_{fin})$, 
where ${\bf LBA}_{fin}$ is a ``cyclic'' prop (it has the same generators 
and relations as ${\bf LBA}$, except that diagrams with cycles are allowed). 
Even though compatibility with duality and doubling operations can be 
formulated for the prop ${\bf LBA}$, the proof of \cite{EK1} uses diagrams
involving cycles.

Our scheme of proof is the following. We first prove that the EK quantization 
functors are compatible with duality. More precisely, we show that the 
behavior of $Q_{\Phi}$ w.r.t. a dihedral group of duality transformations 
is related to the duality transformation of associators. We then show that 
any quantization functor $Q$ which is compatible with duality
(i.e., satisfies the propic analogue of $\tilde Q(\A^{*cop}) \simeq 
\tilde Q(\A)^{\star cop}$) is automatically compatible with doubles. The proof 
is a propic version of the following argument. Let $\A$ be a finite 
dimensional Lie bialgebra and let ${\mathfrak D}'(\A)
:= \tilde Q^{-1}(D(\tilde Q(\A)))$; here $\tilde Q^{-1}$ is the 
dequantization functor $\{$QUE algebras$\} \to \{\kk[[\hbar]]$-Lie 
bialgebras$\}$ inverse to $\tilde Q$. Since we have morphisms 
$\tilde Q(\A)\to D(\tilde Q(\A))$, 
$\tilde Q(\A)^{\star cop}\to D(\tilde Q(\A))$ and $\tilde Q$ is compatible 
with duals, we have Lie bialgebra morphisms $\A\to {\mathfrak D}'(\A)$ and
$\A^{*cop}\to {\mathfrak D}'(\A)$. By flatness, we then have 
${\mathfrak D}'(\A) = \A\oplus\A^{*cop}$; the structure of 
${\mathfrak D}'(\A)$ is then uniquely determined by the brackets between
$\A$ and $\A^{*cop}$. While for a given $\A$, these brackets could be such 
that ${\mathfrak D}'(\A)\ncong {\mathfrak D}(\A)$, one can prove that the 
condition that these brackets are propic imply a propic statement, which 
implies that ${\mathfrak D}'(\A)\simeq {\mathfrak D}(\A)$.

More generally, one can prove using \cite{Enr:sha} that any quantization 
functor is compatible with duality; then Theorem \ref{thm:doubles} implies 
that it is also compatible with doubling operations.

The results of this paper have several corollaries. First, they are necessary
ingredients to prove that when $\A$
is a Kac-Moody Lie bialgebra, $\tilde Q_{\Phi}(\A) \simeq \DJh(\A)$,
where $\DJh(\A)$ is the Drinfeld-Jimbo quantization of $\A$ (\cite{EK3}).
In the final part of this paper, we show that these arguments can then be 
generalized to the case of some affine Lie superalgebras (these Lie superalgebras 
were introduced in \cite{K,vL} and their quantum versions are due to 
\cite{KhT,Yam}). 

Similarly to \cite{EK3}, this can be used for proving a Kohno-Drinfeld-type 
theorem for these affine Lie superalgebras, namely the braid group representations 
arising from the monodromy of the KZ equations for with values in $\g$-modules 
from category $\cato$ are equivalent to those arising from the 
quantum $R$-matrix of $\DJh(\G)$.

\section{Reminders on quantization functors}

In this section, we recall the basic definitions in the theory of props, and the 
construction of Etingof-Kazhdan quantization functors (see \cite{EH}).  
We then introduce some notions which we will use in this paper: 
prop bimodules, the dual of a prop,  and the biprops ${\bf LBA}_{2}$
and $\ul{\on{Bialg}}_{2}$.  

\subsection{Schur functors} Let $\on{Sch}$ be the category whose objects are
polynomial Schur functors, i.e., endofunctors $F$ of the category $\on{Vect}$
of finite dimensional vector spaces, of the form
$F(V) = \oplus_{n\geq 0}\oplus_{\pi\in\hat S_{n}} F_{n,\pi}\otimes
p_{\pi}(V^{\otimes n})$, where for each $n\geq 0$ and
irreducible $S_{n}$-module $\pi$, $p_{\pi}\in \QQ S_{n}$ is a chosen rank one
projector in $\on{End}(\pi)\subset \oplus_{\pi'\in\hat S_{n}}\on{End}(\pi')
= \QQ S_{n}$; in the above sum, the $F_{n,\pi}$ are finite dimensional vector
spaces,
which vanish for almost every $(n,\pi)$. Then $\on{Sch}$ is an abelian tensor
category, where morphisms are natural
transformations. The set $\on{Irr(Sch)}$ of irreducible Schur functors is in
bijection
with $\{(n,\pi) | n\geq 0,\pi\in\hat S_{n}\}$; the bijection takes $(n,\pi)$ to
$(V\mapsto p_{\pi}(V^{\otimes n}))$. The category $\on{Sch}$ is equipped
with an involution, $F^*(V) := F(V^*)^*$. As a tensor category, $\on{Sch}$ is
generated by the identity functor ${\bf id}$, such that ${\bf id}(V)=V$.
The neutral object is ${\bf 1}$ defined by ${\bf 1}(V) = \kk$.
We will also use the $n$th symmetric power functor $S^{n}$, the $n$th exterior
power functor $\wedge^{n}$, and the $n$th tensor power functor
$T_{n} = {\bf id}^{\otimes n}$. We say that
$F\in\on{Ob(Sch)}$
has degree $n$ if $F_{n',\pi'}=0$ for $n'\neq n$; we then denote by $|F|$ the
degree of $F$.

Any $R\in \on{Ob(Sch)}$ gives rise to an endofunctor of $\on{Sch}$ by
$\on{Ob(Sch)}\ni F \mapsto F \circ R \in \on{Ob(Sch)}$. For $f\in
\on{Sch}(F,G)$, the corresponding morphism is $f(R)=f\circ R\in
\on{Sch}(F\circ R,G\circ R)$. By functoriality, one also defines a
(in general nonlinear) operation $\on{Sch}(F,G) \to
\on{Sch}(R\circ F,R\circ G)$, $f\mapsto R(f)$.

We also define a completion ${\bf Sch}$ of $\on{Sch}$, where objects are
functors $\on{Vect} \to {\bf Vect}$ (where ${\bf Vect}$ is the category of
vector spaces), of the same form as above, except that the condition that
almost every $F_{n,\pi}$ vanishes is dropped. Then ${\bf Sch}$ is again a
tensor category equipped with an involution.
The symmetric algebra functor $S = \oplus_{n\geq 0} S^{n}$ is an example of
an object of ${\bf Sch}$. The composition $R\circ F$ is an object of 
${\bf Sch}$ if: (a) $R,F\in \on{Ob}({\bf Sch})$, and $F$ has vanishing 
zero degree component; or if (b) $R\in\on{Ob(Sch)}$, $F\in\on{Ob}({\bf Sch})$.

\subsection{Props} \label{sect:basics}

A prop $P$ is a symmetric tensor category generated by a single object; 
equivalently, it is equipped with a tensor functor $i_P : \on{Sch}\to P$
(also denoted $x\mapsto x^P$), inducing the identity on objects; so 
$\on{Ob(Sch)} = \on{Ob}(P)$. The prop $P$ is characterized by 
vector spaces $P(F,G)$ (for $F,G\in\on{Sch}$), composition maps 
$\circ : P(F,G) \otimes P(G,H) \to P(F,H)$, external product maps
$\boxtimes : P(F,G)\otimes P(F',G') \to P(F\otimes F',G\otimes G')$
and $i_P : \on{Sch}(F,G)\to P(F,G)$. 

A ${\bf Sch}$-prop ${\bf P}$ is a symmetric tensor category with the 
same properties, except that the tensor functor is now $i_{\bf P} : 
{\bf Sch} \to {\bf P}$. 

If $P$ is a prop and $H\in\on{Ob(Sch)}$, then we define a prop 
$H(P)$ by $H(P)(F,G) = P(F\circ H,G\circ H)$. If ${\bf P}$ is
a ${\bf Sch}$-prop and $H\in\on{Ob}({\bf Sch})$, then one defines similarly 
a prop $H({\bf P})$. 

If $P,Q$ are props (or ${\bf Sch}$-props), a morphism $f : P\to Q$
is a functor such that $f\circ i_P = i_Q$. If $P$ is a prop, then a 
prop ideal of $P$ is a collection $I(F,G)$, such that the projection 
$P/I$ defined by $(P/I)(F,G):= P(F,G)/I(F,G)$ is a prop and the canonical 
projection $P\to P/I$ is a prop morphism. If $I,J\subset P$ are prop 
ideals, then the product $IJ$ is the smallest prop ideal, such that 
$IJ(F,G)$ contains both $I(F,G)$ and $J(F,G)$. If $x\in P(F,G)$, we denote 
by $\langle x\rangle \subset P$ the smallest prop ideal of $P$ 
such that $x\in \langle x\rangle(F,G)$.

Similarly to algebras, props can be defined by generators and relations
(see \cite{EH}). 

If $\cS$ is a symmetric monoidal category, and $V\in\on{Ob}(\cS)$, 
then one defines a prop $\on{Prop}(V)$ by $\on{Prop}(V)(F,G):= 
\cS(F(V),G(V))$. A $P$-module (in the category $\cS$)
is an object $V$ equipped with a prop morphism $P\to\on{Prop}(V)$. 
The $P$-modules form a category; one reovers in this way the categories
of Lie bialgebras, QUE algebras, etc.  

Let $P$ be a prop (or a ${\bf Sch}$-prop). If $\xi\in P({\bf id},{\bf id})$, 
then there is a unique collection   
$(\xi_F)_{F\in\on{Ob(Sch)}}$, where $\xi_F\in P(F,F)$ is such that 
$\xi_{F\otimes G}=\xi_F\boxtimes \xi_G$, and for $x\in\on{Sch}(F,G)$, 
$\xi_G \circ x^P = x^P \circ \xi_F$ (this construction, called inflation, 
is carried out in a
the more general setup in Section \ref{prop:bimodules}). We also have 
$(\xi\circ\eta)_{F} = \xi_{F}\circ \eta_{F}$, so if $\xi$ is invertible
then so are the $\xi_{F}$ (recall that the $P(F,F)$ are algebras with unit). We have 
a unique automorphism $\on{Inn}(\xi)$ (we call it an inner automorphism), 
whose action on $P(F,G)$ is $a \mapsto \xi_G \circ a \circ \xi_F^{-1}$. 
The map $P({\bf id},{\bf id})^\times \to \on{Aut}(P)$ is a group morphism.

\subsection{The structure of some props}

The prop $\on{LA}$ of Lie algebras is generated by $\mu\in\on{LA}(\wedge^2,{\bf
id})$ and the Jacobi relation $\mu \circ (\mu\boxtimes 
\on{id}_{\bf id}) \circ ((123)+(231)+(312))=0$. The prop 
$\on{LCA}$ of Lie coalgebras is generated by $\delta\in\on{LCA}({\bf
id},\wedge^2)$ and the co-Jacobi relation $((123)+(231)+(312)) \circ
(\delta\boxtimes \on{id}_{\bf id}) = 0$. The prop $\on{LBA}$
of Lie bialgebras is generated by $\mu\in\on{LBA}(\wedge^2,{\bf id})$, 
$\delta\in\on{LBA}({\bf id},\wedge^2)$, the Jacobi and co-Jacobi relations, 
and the cocycle relation $\delta \circ \mu = ((12)-(21)) \circ (\mu\boxtimes
\on{id}_{\bf id}) \circ (\on{id}_{\bf id} \boxtimes \mu) 
\circ ((12)-(21))$. 

We then have prop morphisms $\on{LA}\to \on{LBA}$ and $\on{LCA}\to\on{LBA}$,
taking $\mu,\delta$ to their homonyms. Composing these inclusions with 
composition, we get a map
$$
\oplus_{Z\in\on{Irr(Sch)}} \on{LCA}(F,Z)\otimes \on{LA}(Z,G) \to \on{LBA}(F,G);  
$$ 
one shows that this map is a linear isomorphism (\cite{univ:der,Pos}). More 
generally, one shows that if $(F_i)_{i\in I}$, $(G_j)_{j\in J}$ are finite 
families of Schur functors, then the map 
\begin{equation} \label{isoms:LBA}
\oplus_{(Z_{ij})_{i,j}\in \on{Irr(Sch)}^{I\times J}}
(\otimes_{i\in I} \on{LCA}(F_i,\otimes_{j\in J}Z_{ij})) 
\otimes 
(\otimes_{j\in J} \on{LA}(\otimes_{i\in I}Z_{ij},G_j)) \to 
\on{LBA}(\otimes_{i\in I}F_i,\otimes_{j\in J}G_j), 
\end{equation}
taking $(\otimes_{i\in I} \kappa_i) \otimes (\otimes_{j\in J} \lambda_j)$
to $(\boxtimes_{j\in J}\lambda_j) \circ \sigma_{I,J}^{\on{LBA}} \circ 
(\boxtimes_{i\in I} \kappa_i)$ (where $\sigma_{I,J}\in \on{Sch}(
\otimes_{i\in I}(\otimes_{j\in J}Z_{ij}),
\otimes_{j\in J}(\otimes_{i\in I}Z_{ij}))$ is the map 
$\otimes_{i\in I}(\otimes_{j\in J}z_{ij}) \mapsto 
\otimes_{j\in J}(\otimes_{i\in I}z_{ij})$), is a linear isomorphism. 

The props $\on{LA}$ and $\on{LCA}$ give rise to ${\bf Sch}$-props
${\bf LA}$, ${\bf LCA}$, where for ${\bf F},{\bf G}\in \on{Ob}({\bf Sch})$, 
${\bf LA}({\bf F},{\bf G}) 
= \hat\oplus_{\alpha,\beta} \on{LA}(F_\alpha,G_\beta)$ 
($F_\alpha$ is the degree $\alpha$ component of ${\bf F}$)
and ${\bf LCA}({\bf F},{\bf G})$ is defined in the same way. 

The prop $\on{LBA}$ is graded by $\NN^2$ (with $|\mu|=(1,0)$ and 
$|\delta|=(0,1)$), which implies that it can be completed to a 
${\bf Sch}$-prop, given by  
$$
{\bf LBA}({\bf F},{\bf G}) = \hat\oplus_{Z\in\on{Irr(Sch)}}
{\bf LCA}({\bf F},Z) \otimes {\bf LA}(Z,{\bf G}) 
$$
(here $\hat\oplus$ is the direct product). We have also  
$$
{\bf LBA}(\otimes_{i\in I}{\bf F}_i,\otimes_{j\in J}{\bf G}_j) = 
\hat\oplus_{(Z_{ij})_{i,j}\in\on{Irr(Sch)}^{I\times J}}
(\otimes_{i\in I}{\bf LCA}({\bf F}_i,\otimes_{j\in J}Z_{ij})) 
\otimes (\otimes_{j\in J}{\bf LA}(\otimes_{i\in I}Z_{ij},{\bf G}_j)).  
$$

The prop $\on{Alg}$ of associative algebras is generated by 
$\eta\in\on{Alg}({\bf 1},{\bf id})$ and $m\in \on{Alg}(T_2,{\bf id})$, 
with relations $m \circ (\eta\boxtimes \on{id}_{\bf id}) = m \circ 
(\on{id}_{\bf id}\boxtimes \eta)=\on{id}_{\bf id}$, and 
$m \circ (m\boxtimes \on{id}_{\bf id}) = m\circ (\on{id}_{\bf id}\boxtimes m)$
(associativity). The prop $\on{Coalg}$ of associative coalgebras is similarly
generated by $\varepsilon\in\on{Coalg}({\bf id},{\bf 1})$ and $\Delta\in 
\on{Coalg}({\bf id},T_2)$, with relations $(\varepsilon\boxtimes
\on{id}_{\bf id}) \circ \Delta = (\on{id}_{\bf id} \boxtimes \varepsilon)
\circ \Delta = \on{id}_{\bf id}$ and $(\Delta\boxtimes\on{id}_{\bf id})\circ
\Delta = (\on{id}_{\bf id}\boxtimes\Delta)\circ \Delta$ (coassociativity). 
The prop $\on{Bialg}$
of bialgebras is generated by $m,\Delta,\eta,\varepsilon$ as above, 
with the same relations and the additional compatibility relation\footnote{We 
denote by $(\sigma(1)...\sigma(n))$ the permutation $\sigma$, which 
we view as an element of $\on{Sch}(F^{\otimes n},F^{\otimes n})$
is $F$ is a Schur functor (here $F={\bf id}$) and we denote in the same way 
its image by $i_{P}$ in $P(F^{\otimes n},F^{\otimes n})$ (here $P=\on{Bialg}$).} 
$\Delta \circ m = (m\boxtimes m) \circ (1324) \circ (\Delta\boxtimes\Delta)$. 

Applying to $\on{id_{\bf id}^{\on{Bialg}}} - \eta\circ\varepsilon\in 
\on{Bialg}({\bf id},{\bf id})$ the inflation procedure (see Section
\ref{sect:basics}), we construct for any $Z\in\on{Irr(Sch)}$ an element 
$(\on{id}_{\bf id}^{\on{Bialg}} - \eta\circ\varepsilon)_Z\in 
\on{Bialg}(Z,Z)$ (shortly denoted $(\on{id} - \eta\circ\varepsilon)_Z$). 
We have prop morphisms $\on{Alg}\to \on{Bialg}$ and $\on{Coalg} \to 
\on{Bialg}$, taking $m,\eta$ and $\Delta,\varepsilon$ to their homonyms, 
which give rise to a map 
\begin{equation}\label{1.5}
\oplus_{Z\in\on{Irr(Sch)}} \on{Coalg}(F,Z)\otimes \on{Alg}(Z,G) \to 
\on{Bialg}(F,G)
\end{equation}
by $\oplus_{Z} c_Z \otimes a_Z \mapsto \sum_Z a_Z \circ
(\on{id}-\eta\circ\varepsilon)_Z \circ c_Z$. One can show that this map 
is a linear isomorphism. More generally, if 
$(F_i)_{i\in I}$ and $(G_j)_{j\in J}$ are finite families of 
objects of $\on{Sch}$, the map 
\begin{equation} \label{isom:bialg}
\oplus_{(Z_{ij})_{i,j}\in\on{Irr(Sch)}^{I\times J}} 
(\otimes_{i\in I}\on{Coalg}(F_i,\otimes_{j\in J}Z_{ij}))
\otimes (\otimes_{j\in J}\on{Alg}(\otimes_{i\in I}Z_{ij},G_j)) \to 
\on{Bialg}(\otimes_{i\in I}F_i,\otimes_{j\in J}G_j)
\end{equation}
direct sum of $(\otimes_{i\in I}c_i)\otimes(\otimes_{j\in J}a_j)\mapsto 
(\boxtimes_{j\in J}a_j) \circ \sigma_{I,J}^{\on{Bialg}}
\circ \boxtimes_{i\in I}(\boxtimes_{j\in J}
(\on{id}-\eta\circ\varepsilon)_{Z_{ij}})
\circ (\boxtimes_{i\in I}c_i)$ is a linear isomorphism. 

We then define the completed prop $\ul{\on{Bialg}} := \lim_{\leftarrow}
\on{Bialg}/\langle\on{id}_{\bf id}^{\on{Bialg}}
- \eta\circ\varepsilon\rangle^n$ (the prop ideal $\langle x\rangle^n$ has been 
defined in Section \ref{sect:basics}). Then one can prove that the
isomorphism (\ref{1.5})
$$
\ul{\on{Bialg}}(F,G) \simeq \hat\oplus_{Z\in\on{Irr(Sch)}} \on{Coalg}(F,Z)
\otimes \on{Alg}(Z,G).
$$
More generally, the isomorphism (\ref{isom:bialg}) extends to an 
isomorphism 
$$
\ul{\on{Bialg}}(\otimes_{i\in I}F_i,\otimes_{j\in J}G_j)
\simeq \hat\oplus_{(Z_{ij})_{i,j}\in\on{Irr(Sch)}^{I\times J}} 
(\otimes_{i\in I}\on{Coalg}(F_i,\otimes_{j\in J}Z_{ij}))
\otimes (\otimes_{j\in J}\on{Alg}(\otimes_{i\in I}Z_{ij},G_j)). 
$$
If $\cS = \{$topologically free $\kk[[\hbar]]$-modules$\}$, then 
the category of $\cS$-modules over $\ul{\on{Bialg}}$ is $\{$QUE algebras$\}$; 
if $\cS = \{\kk[[\hbar]]$-modules of the form $V[[\hbar]]$, where $V$ is a 
pro-vector space$\}$, then this category is $\{$QFSH algebras$\}$. 

We define Hopf as the prop with the same generators and relations
as Bialg, with the additional generator
$a\in \on{Hopf}({\bf id},{\bf id})$ (the antipode), and additional relations:
$a$ is invertible,   
$$
m \circ (\on{id}_{{\bf id}} \boxtimes a) \circ \Delta = 
m \circ (a\boxtimes \on{id}_{{\bf id}}) \circ \Delta = 
\on{id}_{{\bf id}}, 
$$
These relations imply $m \circ a^{\boxtimes 2} = a \circ m \circ (21)$, 
$a^{\boxtimes 2} \circ \Delta = (21) \circ \Delta \circ a$, 
$\varepsilon\circ a = \varepsilon$, $a\circ \eta = \eta$
(see \cite{K}).

The natural morphism $\on{Bialg}\to \ul{\on{Bialg}}$ factors as 
$\on{Bialg} \to \on{Hopf} \to \underline{\on{Bialg}}$, where the morphism
$\on{Hopf} \to \underline{\on{Bialg}}$ is such that $a$ maps to
$\eta \circ \varepsilon - (\on{id}_{\bf id}^{\ul{\on{Bialg}}} - \eta \circ 
\varepsilon)
+ m^{(2)} \circ (\on{id}_{\bf id}^{\ul{\on{Bialg}}} 
- \eta \circ \varepsilon)^{\boxtimes 2} \circ
\Delta^{(2)} - m^{(3)} \circ (\on{id}_{\bf id}^{\ul{\on{Bialg}}} -
\eta \circ \varepsilon)^{\boxtimes 3} \circ
\Delta^{(3)} + ...$, where $m^{(3)} = m\circ (m\boxtimes \on{id}_{{\bf id}})$,
$\Delta^{(3)} = (\Delta\boxtimes \on{id}_{{\bf id}}^{\ul{\on{Bialg}}}) 
\circ \Delta$, etc.

\subsection{Definition of quantization functors}

A quantization functor of Lie bialgebras is a prop morphism $Q : \on{Bialg}
\to S({\bf LBA})$, such that: (a) the composed morphism $\on{Bialg} \to 
S({\bf LBA})\to S({\bf Sch})$ (the second morphism is $\mu\to 0$, 
$\delta\to 0$) is $m\mapsto m^{S}$, $\Delta\mapsto \Delta^{S}$, 
$\eta\mapsto pr_{0}$, $\varepsilon\mapsto inj_{0}$, where
$m^{S}\in S({\bf Sch})(T_{2},{\bf id}) = {\bf Sch}(S^{\otimes 2},S)$
and $\Delta^{S}\in S({\bf Sch})({\bf id},T_{2}) = {\bf Sch}(S,S^{\otimes 2})$
are the universal versions of the product and coproduct map of the symmetric
algebra $S(V)$, where $V$ is a vector space; and $pr_{0}\in {\bf Sch}(S,{\bf 1})$, 
$inj_{0}\in {\bf Sch}({\bf 1},S)$ are the canonical projection and injection; 
(b) the element\footnote{We denote by $pr_{k}\in{\bf Sch}(S,S^{k})$ and 
$inj_{k}\in {\bf Sch}(S^{k},S)$ the canonical projection and injection.} 
$pr_{1}^{\bf LBA} \circ Q(m) \circ (inj_{1}^{{\bf LBA}})^{\boxtimes 2}
\circ inj_{\wedge^{2}} \in {\bf LBA}(\wedge^{2},{\bf id})$
equals $\mu$, the element $pr_{\wedge^{2}}\circ 
(pr_{1}^{\bf LBA})^{\boxtimes 2}\circ Q(\Delta)\circ inj_{1}^{{\bf LBA}}
\in {\bf LBA}({\bf id},\wedge^{2})$ equals $\delta$ (here $pr_{\wedge^{2}}
\in \on{Sch}(T_{2},\wedge^{2})$ and $inj_{\wedge^{2}}\in
\on{Sch}(\wedge^{2},T_{2})$ are the canonical morphisms 
corresponding to the decomposition $T^{2}=S^{2}\oplus\wedge^{2}$). 

One proves that any quantization functor $Q:\on{Bialg} \to S({\bf LBA})$ factors
uniquely as $\on{Bialg}\to \ul{\on{Bialg}}\to S({\bf LBA})$, where 
the induced morphism $\ul{\on{Bialg}}\to S({\bf LBA})$ (also denoted $Q$)
is an isomorphism (see \cite{EK2,EE}).

\subsection{Schur multifunctors}

If $I$ is a finite set, then $\on{Sch}_{I}$ is the category of
polynomial functors $\on{Vect}^{I}\to \on{Vect}$, i.e., its objects are 
the functors of the
form $(V_{i})_{i\in I} \mapsto \oplus_{(n_{i},\pi_{i})_{i\in I}}
F_{(n_{i},\pi_{i})_{i\in I}} \otimes
(\otimes_{i\in I} p_{i}(V_{i}^{\otimes n_{i}}))$, where the
$F_{(n_{i},\pi_{i})_{i\in I}}$ are almost all zero.
We let ${\bf Sch}_{I}$ be the category whose objects are functors
$\on{Vect}^{I}\to {\bf Vect}$ of the same form as above, except that the
vanishing condition is dropped. 
If $I = \emptyset$, then $\on{Ob}(\on{Sch}_\emptyset) \subset 
\on{Ob}(\on{Vect})^{(\on{Ob}(\on{Vect})^\emptyset)}
= \on{Ob}(\on{Vect})$, and by convention $\on{Ob}(\on{Sch}_\emptyset) 
= \on{Ob}({\bf Sch}_\emptyset) = \on{Ob}(\on{Vect})$; we denote by 
${\mathfrak 1} \in \on{Ob}(\on{Sch}_\emptyset)$ the object corresponding 
to the base field $\kk$. 

The category $\on{Sch}_{I}$ (as well as ${\bf Sch}_{I}$) is equipped 
with a duality autofunctor
$F^{*}((V_{i})_{i\in I}) := F((V_{i}^{*})_{i\in I})^{*}$.
We set $\on{Sch}_{n} := \on{Sch}_{\{1,...,n\}}$,
${\bf Sch}_{n} := {\bf Sch}_{\{1,...,n\}}$.

Define a natural transformation $\boxtimes : \on{Sch}_I \times
\on{Sch}_J \to \on{Sch}_{I\sqcup J}$, $(F\boxtimes G)((V_k)_{k\in I\sqcup J})
:= F((V_{i})_{i\in I}) \otimes G((V_j)_{j\in J})$. The irreducible
objects of $\on{Sch}_{I}$ are then the $\boxtimes_{i\in I}Z_{i}$, where
$(Z_{i})_{i\in I}\in \on{Irr(Sch)}^{I}$.

We define a functor $\otimes : \on{Sch}_{I}\to \on{Sch}$,
taking $A : \on{Vect}_{I}\to \on{Vect}$ to $A\circ diag$, where
$diag : \on{Vect}\to \on{Vect}^{I}$ is the diagonal embedding.
Then $\otimes(\boxtimes_{i\in I} F_{i}) =  \otimes_{i\in I} F_{i}$.

Define the tensor category of Schur multifunctors
$\on{Sch}_{{(1)}}$ by 
$\on{Ob}(\on{Sch}_{{(1)}}):= \sqcup_{n\geq 0}
\on{Ob}(\on{Sch}_{n})$, $\on{Sch}_{(1)}(F,G)=\on{Sch}_{n}(F,G)$
if $F,G\in\on{Ob}(\on{Sch}_{n})$ and $\on{Sch}_{(1)}(F,G)=0$
if $F\in\on{Ob}(\on{Sch}_{n})$, $G\in\on{Ob}(\on{Sch}_{n'})$ with $n'\neq n$, 
and the tensor product is $(F,G)\mapsto F\boxtimes G$, where the 
identification $\{1,...,n\}\sqcup\{1,...,m\} \simeq\{1,...,n+m\}$
is given by $i\mapsto i$ for $i\in\{1,...,n\}$ and $i\mapsto n+i$
for $i\in \{1,...,m\}$.

We define the tensor category of Schur multibifunctors $\on{Sch}_{(1+1)}$
by $\on{Ob}(\on{Sch}_{(1+1)}) = \sqcup_{p,q\geq 0}\on{Sch}_{p+q}$. 
We denote by $\ul\boxtimes : \on{Ob}(\on{Sch}_{(1)})^2 \to
\on{Ob}(\on{Sch}_{(1+1)})$ the map obtained from $\boxtimes : 
\on{Sch}_p \times \on{Sch}_q \to \on{Sch}_{p+q}$. For $F,G\in 
\on{Ob}(\on{Sch}_{(1)})$, 
$$
\on{Sch}_{(1+1)}(F\ul\boxtimes G,F'\ul\boxtimes G') := 
\on{Sch}_{(1)}(F,F') \ul\boxtimes \on{Sch}_{(1)}(G,G'). 
$$
The tensor product is defined by $(F\ul\boxtimes G)\otimes (F'\ul\boxtimes G') 
:= (F\boxtimes G)\ul\boxtimes (F'\ul\boxtimes G')$. The tensor category 
${\bf Sch}_{(1+1)}$ is defined in the same way. 

\subsection{(Quasi)multibiprops}

A multibiprop $P_{(1+1)}$ is a tensor category equipped with a
tensor functor $i_{P_{(1+1)}}\on{Sch}_{(1+1)} \to P_{(1+1)}$, inducing the identity 
on objects. More explicitly, we have composition maps $P_{(1+1)}(F,G)
\otimes P_{(1+1)}(G,H) \to P_{(1+1)}(F,H)$, external product maps 
$P_{(1+1)}(F,G) \otimes P_{(1+1)}(F',G') \to P_{(1+1)}(F\otimes F',
G\otimes G')$ and maps $i_{P_{(1+1)}} : \on{Sch}_{(1+1)}(F,G) \to 
P_{(1+1)}(F,G)$, satisfying some axioms. 

In the case of a quasi-multibiprop, the composition maps are 
$P_{(1+1)}(F,G) \otimes P_{(1+1)}(G,H) \supset D(F,G,H)
\to P_{(1+1)}(F,H)$, where $D(F,G,H)$ is a vector subspace of   
$P_{(1+1)}(F,G) \otimes P_{(1+1)}(G,H)$; the above axioms 
(e.g., associativity) should be satisfied when the involved expressions 
make sense.

\subsection{Traces in some props}

If $I,J$ are finite sets and $\Sigma\subset I \times J$, and for $A\in 
\on{Ob(Sch}_{I})$, $B\in \on{Ob(Sch}_{J})$, we define
$$\on{LBA}^{\Sigma}(A,B)\subset \on{LBA}(\otimes(A),\otimes(B))
$$ as follows.
If $A = \boxtimes_{i} F_{i}$ and $B  = \boxtimes_{j} G_{j}$, where
$F_{i},G_{j}\in\on{Irr(Sch)}$, then
$$
\on{LBA}^{\Sigma}(A,B) \simeq \oplus_{Z_{ij} | Z_{ij}={\bf 1}\on{\ if\ }
(i,j)\notin \Sigma} \on{LCA}(F_{i},\otimes_{j\in J}Z_{ij}) \otimes
\on{LA}(\otimes_{i\in I}Z_{ij},G_{j}).
$$
under decomposition (\ref{isoms:LBA}). 
We extend this definition to any $A,B$ by linearity.
We define in the same way ${\bf LBA}^{\Sigma}({\bf A},{\bf B})$
for ${\bf A}\in \on{Ob}({\bf Sch}_{I})$, ${\bf B}\in \on{Ob}({\bf Sch}_{J})$.

If $I,J$ are finite sets, $a\in I\cap J$, and $\Sigma\subset I\times J$
is such that $(a,a)\notin \Sigma$, define $I_{a}:= I - \{a\}$, $J_{a}:= 
J - \{a\}$ and $\Sigma(a)\subset I_{a}\times J_{a}$ 
as $\Sigma(a):= \{(i,j)|(i,j)\in\Sigma$, or $(i,a)$ and $(a,j)\in
\Sigma(a)\}$. 

Then for $A\in \on{Ob(Sch}_{I_{a}})$, $B\in \on{Ob(Sch}_{J_{a}})$, 
$F\in \on{Ob(Sch)}$, we view $I$ as $I_{a}\sqcup \{a\}$, $J$ as 
$J_{a}\sqcup \{a\}$, $A\boxtimes F$ (resp., $B\boxtimes F$) as an object in
$\on{Sch}_{I}$ (resp., $\on{Sch}_{J}$) and we define a
trace map
$$
\on{tr}_{F} : \on{LBA}^{\Sigma}(A\boxtimes F,B\boxtimes F) \to
\on{LBA}^{\Sigma(a)}(A,B)
$$
by the condition that the diagram 
$$
\xymatrix{
\on{LBA}^{\Sigma}(A\boxtimes F,B\boxtimes F)  
\ar[r]^{\ \ \ \ \  \ \ \on{tr}_{F}}& 
\on{LBA}^{\Sigma(a)}(A,B) \\
*\txt{
$\oplus_{Z_{ij}|Z_{ij}={\bf 1}\on{if}(i,j)\in\bar\Sigma}$\\
$(\otimes_{i\in I}\on{LCA}(A_{i},\otimes_{j\in\bar J} Z_{ij}))
\otimes \on{LCA}(F,\otimes_{j\in\bar J}Z_{aj})$\\
$\otimes (\otimes_{j\in J}\on{LA}(\otimes_{i\in \bar I}Z_{ij},B_{j}))
\otimes \on{LA}(\otimes_{i\in\bar I}Z_{ia},F)$ 
}\ar[u]^{\sim}
\ar[ur]
}
$$
commutes, where the diagonal map takes $(\otimes_{i\in I}\kappa_{i})\otimes 
\kappa_{a} \otimes (\otimes_{j}\lambda_{j}) \otimes \lambda_{a}$
to $(\boxtimes_{j\in J}\lambda_{j}) \circ (\tau_{J,I,a}^{-1})^{\on{LBA}} \circ 
(\sigma_{I,J}^{\on{LBA}} \boxtimes (\kappa_{a}\circ\lambda_{a}))
\circ \tau_{I,J,a}^{\on{LBA}}\circ (\boxtimes_{i\in I}\kappa_{i})
\in \on{LBA}(\otimes(A),\otimes(B))$ (see \cite{EH}). Here 
$$\tau_{I,J,a}\in 
\on{Sch}(\otimes_{i\in I_{a}}(\otimes_{j\in J}Z_{ij}),
\otimes_{i\in I_{a}}(\otimes_{j\in J_{a}}Z_{ij})\otimes 
(\otimes_{i\in I_{a}}Z_{ia}))
$$
is the map $\otimes_{i\in I_{a}}(\otimes_{j\in \bar J_{a}}z_{ij}) \mapsto 
\otimes_{i\in I_{a}}(\otimes_{j\in J_{a}}z_{ij}) \otimes (\otimes_{i\in I_{a}}
z_{ia})$.  

This map extends to a trace map 
$\on{tr}_{F} : {\bf LBA}^{\Sigma}({\bf A}\boxtimes {\bf F},{\bf B}
\boxtimes {\bf F}) \to {\bf LBA}^{\Sigma(a)}({\bf A},{\bf B})$ for 
${\bf A}\in\on{Ob}({\bf Sch}_{I})$, ${\bf F}\in\on{Ob}({\bf Sch}_{J})$, 
${\bf F}\in\on{Ob}({\bf Sch})$. 

Trace maps can also be constructed in the case of the prop 
$\on{Bialg}$. We define similarly $\on{Bialg}^{\Sigma}(A,B) 
\subset \on{Bialg}(\otimes(A),\otimes(B))$ by selecting the 
analogous components in the decomposition (\ref{isom:bialg}). 
We then define a map 
$$
\on{tr}_{F} : \on{Bialg}^{\Sigma}(A\boxtimes F,B\boxtimes F)
\to \on{Bialg}^{\Sigma(a)}(A,B); 
$$
the diagonal map is then 
\begin{align*}
& (\otimes_{i\in I_{a}}c_{i})\otimes c_{a} \otimes (\otimes_{j\in J_{a}}a_{j}) 
\otimes a_{a} \mapsto 
(\boxtimes_{j\in J_{a}} a_{j}) \circ (\tau_{J,I,a}^{-1})^{\on{Bialg}} 
\circ 
(\on{id}_{\otimes_{j\in J_{a}}(\otimes_{i\in I_{a}}Z_{ij})}\boxtimes 
(\boxtimes_{j\in J_{a}}(\on{id}-\eta\circ\varepsilon)_{Z_{aj}}))
\\
& \circ (\sigma_{I_{a},J_{a}}^{\on{Bialg}} \boxtimes (c_{a}\circ a_{a}))
\circ \tau_{I,J,a}^{\on{Bialg}}
\circ (\boxtimes_{i\in I_{a}}(\boxtimes_{j\in J}
(\on{id}-\eta\circ\varepsilon)_{Z_{ij}}))
\circ (\boxtimes_{i\in I_{a}} c_{i})
\in \on{Bialg}(\otimes(A),\otimes(B))
\end{align*}
This map extends to a trace map 
$$
\on{tr}_{F} : \ul{\on{Bialg}}^{\Sigma}(A\boxtimes F,B\boxtimes F)
\to \ul{\on{Bialg}}^{\Sigma(a)}(A,B), 
$$
where $\ul{\on{Bialg}}^{\Sigma}(A,B)$ is the direct product of the
summands of ${\on{Bialg}}^{\Sigma}(A,B)$. 

If now $I,J$ are finite sets, $a,b\in I\cap J$ are distinct, and $\Sigma\subset
I\times J$ is such $(a,a)\notin \Sigma$ and $(b,b)\notin\Sigma(a)$, 
then $(b,b)\notin\Sigma$ and $(a,a)\notin\Sigma(b)$, so that 
$\Sigma(a)(b)$ and $\Sigma(b)(a)$ are both defined; moreover, 
$\Sigma(a)(b)=\Sigma(b)(a)$; we denote by $\Sigma(a,b)$ this set. 
If now $A\in\on{Ob}(\on{Sch}_{I_{a,b}})$, 
$B\in \on{Ob}(\on{Sch}_{J_{a,b}})$, $F,G\in\on{Ob}(\on{Sch})$, 
then $\on{tr}_{F} \circ \on{tr}_{G} = \on{tr}_{G}\circ \on{tr}_{F}$ 
(equality of maps $\on{LBA}^{\Sigma}(A\boxtimes F \boxtimes G, 
B\boxtimes F\boxtimes G)\to \on{LBA}^{\Sigma(a,b)}(A,B)$). 
The maps $\on{tr}_{F}$, $\on{tr}_{G}$ defined in the case of 
$\on{Bialg}$ (and of ${\bf LBA}$, $\ul{\on{Bialg}}$)
have the same property. 

This property allows to generalize the trace map to the 
following setup. Let $I,J$ be finite sets and $K\subset I\cap J$. 
Let $\Sigma\subset I\times J$ be such that there is no finite sequence
$(k_1,...,k_p)$ of elements of $K$ such that $(k_1,k_2)\in \Sigma$, ...,
$(k_{p-1},k_p)\in\Sigma$, $(k_p,k_1)\in\Sigma$. Define $\Sigma(K) \subset 
(I-K)\times (J-K)$ be the set of $(i,j)$ such that there exists a 
(possibly empty) sequence $(k_1,...,k_p)$ of elements of $K$, 
such that $(i,k_1)\in\Sigma$, $(k_1,k_2)\in\Sigma$, ..., $(k_p,j)\in\Sigma$. 

If $A\in \on{Ob}(\on{Sch}_{I-K})$, 
$B\in \on{Ob}(\on{Sch}_{J-K})$, and $F_k\in \on{Irr(Sch)}$ for each 
$k\in K$, we define a trace map 
$$
\on{tr}_{\boxtimes_{k\in K}F_k} : 
\on{LBA}^{\Sigma}(A\boxtimes(\boxtimes_{k\in K}F_k),B\boxtimes
(\boxtimes_{k\in K}F_k)) \to \on{LBA}^{\Sigma(K)}(A,B), 
$$
by $\on{tr}_{\boxtimes_{k\in K}F_k} = \circ_{k\in K}\on{tr}_{F_k}$
(the order is irrelevant). This generalizes by linearity to a trace map 
$\on{tr}_F : \on{LBA}^{\Sigma}(A\boxtimes F,B\boxtimes F) 
\to \on{LBA}^{\Sigma(K)}(A,B)$ for any $F\in \on{Sch}_K$. 
One defines similar trace maps in the cases of ${\bf LBA}$, 
$\on{Bialg}$ and $\ul{\on{Bialg}}$.

The above property then generalized as follows: if $I,J$ are finite
sets, if $K,L\subset I\cap J$ are disjoint, and if $\Sigma\subset I\times J$ 
is such that $\Sigma(J\cup K)$ exists, then for 
$A\in \on{Ob}(\on{Sch}_{I-(K\cup L)})$, 
$B\in \on{Ob}(\on{Sch}_{J-(K\cup L)})$, $F\in\on{Ob}(\on{Sch}_K)$, 
$G\in\on{Ob}(\on{Sch}_L)$, we have the equality  
\begin{equation} \label{trtr}
\on{tr}_F \circ \on{tr}_G = \on{tr}_G \circ \on{tr}_F 
= \on{tr}_{F\boxtimes G}
\end{equation}
of maps $\on{LBA}^\Sigma(A\boxtimes F \boxtimes G,
B\boxtimes F \boxtimes G) \to \on{LBA}^{\Sigma(K\cup L)}(A,B)$.

\subsection{The quasi-multibiprop \boldmath$\Pi$\unboldmath} 

We define a quasi-multibiprop $\Pi$ over $\on{Sch}_{(1+1)}$
as follows. For $F,...,G'\in\on{Irr}(\on{Sch}_{(1)})$, we set 
$$
\Pi(F\ul\boxtimes G,F'\ul\boxtimes G'):= 
\on{LBA}(\otimes(F)\otimes \otimes(G')^{*}, 
\otimes(F')\otimes \otimes(G)^{*}). 
$$
The external product map and the functor $i_{\Pi}:
\on{Sch}_{(1+1)}\to\Pi$ are obvious. Let us define the 
composition. If $I,J,I',J'$ are disjoint sets, if 
$F\in\on{Ob}(\on{Sch}_{I})$, ...,
$G'\in\on{Ob}(\on{Sch}_{J'})$, and if 
$\Sigma\subset (I\sqcup J')\times (I'\sqcup J)$, 
define 
$$
\Pi^{\Sigma}(F\ul\boxtimes G,F'\ul\boxtimes G') := 
\on{LBA}^{\Sigma}(F\boxtimes G^{\prime*},F'\boxtimes G^{*}). 
$$
We attach an oriented graph $\Gamma(\Sigma)$ to $\Sigma$: vertices are 
the elements of $I\sqcup ...\sqcup J'$; there is an edge from $x$ 
to $y$ iff $(x,y)\in\Sigma$. 

Then if $I'',J''$ are disjoint and disjoint of $I,...,J'$, and if 
$\Sigma'\subset (I'\sqcup J'')\times (I''\sqcup J')$, then 
$\Sigma$ and $\Sigma'$ are called composable if the graph 
obtained by glueing $\Gamma(\Sigma)$ and $\Gamma(\Sigma')$
has no loops. If this is the case, define $\Sigma''\subset (I\sqcup J'')
\times (I''\sqcup J)$ as the set of all pairs of vertice $(x,y)$ 
related by a sequence of edges of the big graph. We then 
define a linear map 
$$
\Pi^{\Sigma}(F\ul\boxtimes G,F'\ul\boxtimes G')
\otimes \Pi^{\Sigma'}(F'\ul\boxtimes G',F''\ul\boxtimes G'')
\to \Pi^{\Sigma''}(F\ul\boxtimes G,F''\ul\boxtimes G'')
$$
as the map $\on{LBA}^{\Sigma}(F\boxtimes G^{\prime*},F'\boxtimes
G^{*}) \otimes \on{LBA}^{\Sigma'}(F'\boxtimes G^{\prime\prime*},
F''\boxtimes G^{\prime*}) \to 
\on{LBA}^{\Sigma''}(F\boxtimes G^{\prime\prime*},F''\boxtimes
G^{*})$, $x\otimes y \mapsto 
\on{tr}_{G^{\prime*}}\big( (y\boxtimes \on{id}^{\on{LBA}}_{\otimes(G^{*})})
\circ (\on{id}^{\on{LBA}}_{\otimes(F')}\boxtimes 
\sigma^{\on{LBA}}_{\otimes(G^{*}),\otimes(G^{\prime\prime*})})
\circ (x\boxtimes\on{id}^{\on{LBA}}_{\otimes(G^{\prime\prime*})}) 
\big)$, where $\sigma_{F,G}\in\on{Sch}(F\otimes G,G\otimes F)$
is $x\otimes y \mapsto y\otimes x$.  This partially defined composition 
equips $\Pi$ with the structure of a quasi-multibiprop; in particular, the 
associativity of the composition follows from (\ref{trtr}). 

In the same way, one defines a quasi-multibiprop 
\boldmath$\Pi$\unboldmath\  over 
\boldmath$\on{Sch}$\unboldmath$_{(1+1)}$ by 
\boldmath$\Pi$\unboldmath$({\bf F}\underline\boxtimes {\bf G},
{\bf F'}\underline\boxtimes {\bf G'}) =
{\bf LBA}(\otimes ({\bf F}) \otimes \otimes ({\bf G'})^*, 
\otimes ({\bf F'}) \otimes \otimes ({\bf G})^*)$
and quasi-multibiprops $\tilde\Pi$, $\ul{\tilde\Pi}$ by 
$$
\tilde\Pi(F\ul\boxtimes G,F'\ul\boxtimes G'):= \on{Bialg}
(\otimes(F)\otimes \otimes(G^{\prime*}),\otimes(F')\otimes \otimes(G^{*})), 
$$
$$
\ul{\tilde\Pi}(F\ul\boxtimes G,F'\ul\boxtimes G'):= \ul{\on{Bialg}}
(\otimes(F)\otimes \otimes(G^{\prime*}),\otimes(F')\otimes \otimes(G^{*})). 
$$ 

We will need the linear isomorphism 
$$
\Pi(F\underline\boxtimes G,F'\underline\boxtimes G')
\to \Pi((G')^*\underline\boxtimes (F')^*, G^*\underline\boxtimes F^*),
\quad x\mapsto x^*, 
$$
where $F,...,G'\in \on{Ob}(\on{Sch}_{(1)})$, given by the sequence of 
maps 
$\Pi(F\underline\boxtimes G,F'\underline\boxtimes G') \simeq
\on{LBA}(\otimes(F) \otimes \otimes (G^{\prime*}),\otimes (F')
\otimes \otimes (G^*)) \to 
\on{LBA}(\otimes (G^{\prime *}) \otimes \otimes(F),\otimes(G^*)
\otimes \otimes(F')) \simeq \Pi((G')^*\underline\boxtimes (F')^*,
G^*\underline\boxtimes F^*)$, where the middle map is 
$x\mapsto \sigma_{\otimes(F'),\otimes(G^{*})}^{\on{LBA}}\circ x \circ 
\sigma_{\otimes(G^{\prime*}),\otimes(F)}^{\on{LBA}}$, 
where $\sigma_{A,B}\in\on{Sch}(A\otimes B,B\otimes A)$ is the 
permutation operator. 

We have $(x\circ y)^* = y^* \circ x^*$. The isomorphism $x\mapsto x^{*}$
extends when $\Pi$ is replaced by \boldmath$\Pi$\unboldmath. 

We then define an element
$m_\Pi\in $\boldmath$\Pi$\unboldmath$((S\ul\boxtimes S)^{\otimes 2},
S\ul\boxtimes S)$; for $\A$ a finite dimensinal Lie bialgebra, it specializes
to the map $(S(\A)\otimes S(\A^{*}))^{\otimes 2} \simeq 
U({\mathfrak D}(\A))^{\otimes 2}
{\to} U({\mathfrak D}(\A))\simeq S(\A)\otimes S(\A^{*})$, where 
$S(\A)\otimes S(\A^{*})\simeq U({\mathfrak D}(\A))$ is the tensor product 
of inclusions followed by product and the middle map is the product  
of ${\mathfrak D}(\A)$.  A graph for $m_{\Pi}$ is given by the following set 
of edges, if the indices of the functors in the initial object 
$(S\ul\boxtimes S)^{\otimes 2} \simeq (S\boxtimes S)\ul\boxtimes
(S\boxtimes S)$ are $(1+,2+,1-,2-)$ and those of the final object 
$S\ul\boxtimes S$ are $(+,-)$:  $\{$edges of $\Gamma\} = \{(1+,+),
(2+,+),(-,1-),(-,2-),(2+,1-)\}$. 
It follows that both sides of the associativity identity 
$m_{\Pi}\circ (m_{\Pi}\boxtimes 
\on{id}^{\Pi}_{S\ul\boxtimes S})
= m_{\Pi}\circ ( \on{id}^{\Pi}_{S\ul\boxtimes S} \boxtimes m_{\Pi})$ 
make sense; the identity holds. One shows that the possible definitions of the $n$fold
iteration of $m_{\Pi}$ all make sense and coincide; we denote by $m_{\Pi}^{(n)}\in
$\boldmath$\Pi$\unboldmath$((S\ul\boxtimes S)^{\otimes n},S\ul\boxtimes S)$
this $n$fold iteration.

We define $\Delta_0\in $\boldmath$\Pi$\unboldmath$(S\ul\boxtimes S,
(S\ul\boxtimes S)^{\otimes 2})
\simeq {\bf LBA}(S\otimes S^{\otimes 2},S^{\otimes 2}\otimes S)$
as the image of $(\Delta^{S}\boxtimes m^{S})^{{\bf LBA}}$, where 
$\Delta^{S}\in
{\bf Sch}(S,S^{\otimes 2})$ and $m^{S}\in{\bf Sch}(S^{\otimes 2},S)$
are the universal versions of the coproduct and product map of $S(V)$, $V$ a vector
space. Then $\Delta_{0}$ specializes to the standard coproduct of 
$U({\mathfrak D}(\A))$, if $\A$ is a finite dimensional Lie bialgebra. If the 
indices of the functors are as before, a graph for $\Delta_{0}$ is given by the 
set of edges $\{(+,1+),(+,2+),(1-,1),(2-,-)\}$.  Hence both sides of the 
associativity identity 
$(\Delta_{0}\boxtimes\on{id}^{\Pi}_{S\ul\boxtimes S}) \circ\Delta_{0}= 
(\on{id}^{\Pi}_{S\ul\boxtimes S}\boxtimes \Delta_{0})\circ\Delta_{0}$
make sense, and the identity holds. As in the case of $m_{\Pi}$, one denotes
by $\Delta_{0}^{(n)}$ the $n$fold iteration of $\Delta_{0}$. 

We now define a ``universal algebra'' ${\bf U}_n := 
$\boldmath$\Pi$\unboldmath$({\mathfrak 1}\ul\boxtimes {\mathfrak 1}, 
(S\ul\boxtimes S)^{\otimes n})$. For $x,y\in{\bf U}_{n}$, the composition 
$m_{\Pi}^{\boxtimes n}\circ \sigma_{n,2}(S\ul\boxtimes S)^{\Pi} 
\circ (x\boxtimes y)$ makes sense\footnote{Here and below, 
$\sigma_{i,j}\in\on{Sch}(({\bf id}^{\otimes i})^{\otimes j},
({\bf id}^{\otimes j})^{\otimes i})$ corresponds to 
$\otimes_{j'=1}^{j}(\otimes_{i'=1}^{i}x_{i'j'})
\mapsto \otimes_{i'=1}^{i}(\otimes_{j'=1}^{j}x_{i'j'}).$} 
 and we set 
$x\star y:= m_{\Pi}^{\boxtimes n}\circ \sigma_{n}^{\Pi} 
\circ (x\boxtimes y)$. 
This defines an associative product on ${\bf U}_{n}$. 

To a partially defined map $\{1,...,m\}\supset 
D_{\phi}\stackrel{\phi}{\to} \{1,...,n\}$, one associates
$\Delta_{0}^{\phi}\in\Pi((S\ul\boxtimes S)^{\otimes n},
(S\ul\boxtimes S)^{\otimes m})$; then for any $x\in{\bf U}_{n}$, 
the composition $\Delta_{0}^{\phi}\circ x\in {\bf U}_{m}$ is well-defined; 
the resulting (insertion-coproduct) map $x\mapsto x^{\phi}$ is an algebra morphism. 
We also denote $x^{\phi}$ by $x^{\phi^{-1}(1),...,\phi^{-1}(n)}$. 

If $\A$ is a finite dimensional Lie bialgebra, then the collection of 
algebra morphisms ${\bf U}_{n}\to U({\mathfrak D}(\A))^{\otimes n}$ is 
compatible with the insertion-coproduct maps on both sides.  

Let $r\in{\bf U}_{2}$ be the image of the element corresponding
to the identity in $\Pi({\mathfrak 1}\ul\boxtimes{\mathfrak 1},
({\bf id}\ul\boxtimes{\bf 1})\otimes ({\bf 1}\ul\boxtimes
{\bf id}))\subset $\boldmath$\Pi$\unboldmath$
({\mathfrak 1}\ul\boxtimes{\mathfrak 1},
(S\ul\boxtimes S)^{\otimes 2})$ and $r_{ij}:= r^{i,j}\in 
{\bf U}_{n}$. Then the $r_{ij}$ satisfy the classical Yang-Baxter 
relation $[r_{ij},r_{ik}]+[r_{ij},r_{jk}]+[r_{ik},r_{jk}]=0$ for 
$i,j,k$ distinct. 

Let $\t_{n}$ be the infinitesimal braid Lie algebra (\cite{Dr:QH}),
generated by $t_{ij}$, $1\leq i\neq j\neq n$ with relations
$t_{ij}=t_{ji}$, $[t_{ij},t_{ik}+t_{jk}]=0$, $[t_{ij},t_{kl}]=0$
for $i,j,k,l$ distinct. We have Lie algebra morphisms
$\t_{n}\to\t_{m}$ for each partially defined map 
$\{1,...,m\}\supset D_{\phi}\stackrel{\phi}{\to}\{1,...,n\}$, 
$x\mapsto x^{\phi}$, defined by 
$t_{ij}\mapsto \sum_{i'\in\phi^{-1}(i),j'\in\phi^{-1}(j)}t_{i'j'}$. 
Then we have an algebra morphism $\widehat{U(\t_{n})} \to {\bf U}_{n}$, 
defined by $t_{ij}\mapsto r_{ij}+r_{ji}$ (the hat denotes the degree completion 
of $U(\t_{n})$, obtained by assigning degree $1$ to each $t_{ij}$).

\subsection{The Etingof-Kazhdan construction of quantization functors}

Let $\Phi(A,B)$ be a Drinfeld associator; this is a series in noncommuting 
variables $A,B$, which may be identified with an element $\Phi(t_{12},t_{23})
\in \widehat{U(\t_{3})}$. We recall the construction of a quantization 
functor $Q_\Phi = Q : \on{Bialg}\to S({\bf LBA})$ attached to $\Phi$.
We have $Q(\varepsilon) = pr_{0}\in {\bf LBA}(S,{\bf 1})$ 
and $Q(\eta)=inj_{0}\in {\bf LBA}({\bf 1},S)$. It remains to 
construct $Q(m)$ and $Q(\Delta)$.

One first constructs $\on{J}\in{\bf U}_{2}$ of the form $\on{J}=1-r/2+...$,
solution of 
$$
\on{J}^{1,2}\star\on{J}^{12,3}=\on{J}^{2,3}\star\on{J}^{1,23}\star\Phi; 
$$
we denote by $\on{J}^{-1}$ the inverse of $\on{J}$ in the group of 
invertible elements of $({\bf U}_{2},\star)$. One then defines  $\on{Ad(J)} \in
\hbox{\boldmath$\Pi$\unboldmath}((S\underline\boxtimes S)^{\otimes 2},
(S\underline\boxtimes S)^{\otimes 2})$ by
$$\on{Ad}(\on{J}) := (m_{\Pi}^{(3)})^{\boxtimes 2} \circ 
\sigma_{2,3}(S\ul\boxtimes S)^{\Pi} \circ (\on{J}\boxtimes \Delta_{0}
\boxtimes\on{J}^{-1})\in 
\hbox{\boldmath$\Pi$\unboldmath}
(S\ul\boxtimes S,(S\ul\boxtimes S)^{\otimes 2})
$$
A graph for $\on{Ad}(\on{J})$ is then given by the following 
set of edges (if the indices for the factors of the 
source and target $(S\ul\boxtimes S)^{\otimes 2}\simeq 
S^{\boxtimes 2}\ul\boxtimes S^{\boxtimes 2}$ are $1+,2+,1-,2-$
and $1'+,2'+,1'-,2'-$): $(i+,j'+)$, $(i'-,j-)$, $(i'-,j'+)$
for $i,j\in\{1,2\}$.  

It follows that the composition 
$$
\Delta_{\Pi}:= \on{Ad}(\on{J})\circ\Delta_{0}\in
\Pi(S\ul\boxtimes S,(S\ul\boxtimes S)^{\otimes 2})
$$
is well defined and admits the graph (if the indices of the factors of the 
source $S\ul\boxtimes S$ and target $(S\ul\boxtimes S)^{\otimes 2}
\simeq S^{\boxtimes 2}\ul\boxtimes S^{\boxtimes 2}$ are $+,-$ 
and $1+,2+,1-,2-$): $(+,i+)$, $(i-,-)$, $(i-,j+)$ for $i,j\in\{1,2\}$. 

We also construct $\on{R} := \on{J}^{2,1}\star\on{exp}_{\star}(t/2)
\star\on{J}^{-1} \in{\bf U}_{2}$ ($\on{exp}_{\star}$ is the exponential 
in the algebra 
$({\bf U}_{2},\star)$). Then we have $\on{R} = (\on{R}_+ \boxtimes \on{R}_-) 
\circ \on{can}_{S\underline \boxtimes
{\mathfrak 1}}$, where $\on{can}_{S\underline\boxtimes {\mathfrak 1}}
\in \hbox{\boldmath$\Pi$\unboldmath}
({\mathfrak 1} \underline\boxtimes {\mathfrak 1},
(S\underline\boxtimes {\mathfrak 1}) \otimes ({\mathfrak 1}
\underline\boxtimes S))\simeq {\bf LBA}(S,S)$ corresponds to 
$\on{id}_{S}$, and $\on{R}_+ \in \hbox{\boldmath$\Pi$\unboldmath}
(S\underline\boxtimes {\mathfrak 1},
S\underline\boxtimes S)$, $\on{R}_- \in
\hbox{\boldmath$\Pi$\unboldmath}
({\mathfrak 1}\underline\boxtimes S,S\underline\boxtimes S)$.

There exist $\on{R}_+^{(-1)}\in \hbox{\boldmath$\Pi$\unboldmath}
(S\underline\boxtimes {\mathfrak 1},
S\underline\boxtimes S)$ and 
$\on{R}_-^{(-1)}\in \hbox{\boldmath$\Pi$\unboldmath}
({\mathfrak 1} \underline
\boxtimes S,S\underline\boxtimes S)$, such that $\on{R}_+^{(-1)} \circ
\on{R}_+ = \on{id}_{S\underline\boxtimes {\mathfrak 1}}^{\Pi}$, and
$\on{R}_-^{(-1)} \circ \on{R}_-
= \on{id}_{{\mathfrak 1}\underline\boxtimes S}^{\Pi}$. It follows from the
quasitriangular identities satisfied by $\on{R}$ that
$$
\on{R}_+^{(-1)} \circ m_\Pi \circ \on{R}_+^{\boxtimes 2}
= \on{R}_-^* \circ \Delta_\Pi^* \circ (21) \circ
(\on{R}_-^{(-1)*})^{\boxtimes 2},
\quad
(\on{R}_+^{(-1)})^{\boxtimes 2} \circ \Delta_\Pi \circ \on{R}_+
= (\on{R}_-^*)^{\boxtimes 2} \circ m_\Pi^* \circ
(\on{R}_-^{(-1)})^*.
$$
Then we define $m_a\in \hbox{\boldmath$\Pi$\unboldmath}
(S\underline\boxtimes {\mathfrak 1},
(S\underline\boxtimes {\mathfrak 1})^{\otimes 2})$ as the common
value of both sides of the first identity, and
$\Delta_a\in \hbox{\boldmath$\Pi$\unboldmath}
((S\underline\boxtimes {\mathfrak 1})^{\otimes 2},
S\underline\boxtimes {\mathfrak 1})$ as the common
value of both sides of the second identity.

We then define $Q(m)\in {\bf LBA}(S^{\otimes 2},S)$ as the element
corresponding
to $m_a$, and $Q(\Delta)\in {\bf LBA}(S,S^{\otimes 2})$ as the element
corresponding to $\Delta_a$.

\subsection{Prop bimodules} \label{prop:bimodules}
Let $P,Q$ be props, then a $(Q,P)$-prop bimodule $M$ is a symmetric tensor 
category bimodule over the symmetric tensor categories $Q$ and $P$,
equipped with a tensor functor $i_M : \on{Sch} \to M$, such that
$(i_Q,i_M,i_P)$ define a tensor category bimodule morphism from
$\on{Sch}$, equipped with its obvious $(\on{Sch},\on{Sch})$-bimodule
structure, to $M$, equipped with its $(Q,P)$-bimodule structure.

More explicitly, $M$ is a biadditive assignment
$(F,G) \mapsto M(F,G)$ for any $F,G\in \on{Ob}(\on{Sch})$,
equipped with composition maps $P(F,G) \otimes M(G,H) \to M(F,H)$
and $M(F,G) \otimes Q(G,H) \to M(F,H)$ denoted $x\otimes y \mapsto y\circ x$,
external product maps $M(F,G) \otimes M(F',G') \to M(F\otimes F',G\otimes G')$
denoted $x\otimes y \mapsto x\boxtimes y$ and functorial biadditive maps
$i_M(F,G) : \on{Sch}(F,G)\to M(F,G)$, satisfying natural conditions;
for example, the natural diagrams associated to
$P(F,G) \otimes P(G,H) \otimes M(H,K) \to M(F,K)$,
$P(F,G) \otimes M(G,H) \otimes Q(H,K) \to M(F,K)$ and
$M(F,G) \otimes Q(G,H) \otimes Q(H,K) \to M(F,K)$ all commute.

One defines in the same way a prop left (or right) module 
over a  prop $P$. 

If $M$ is a prop module over a prop $P$, and $H$ is a Schur functor, then
$H(M)$ is the prop module over $H(P)$ given by $H(P)(F,G) =
P(F\circ H,G\circ H)$ and $\psi : M_1 \to M_2$ is a prop module
morphism, then $H(\psi) : H(M_1) \to H(M_2)$ is defined by
$H(\psi)(F,G) = \psi(F\circ H,G\circ H)$.

Let $M$ be a $(Q,P)$-prop bimodule.

\begin{prop} \label{inner}
For any $\xi\in M({\bf id},{\bf id})$, there is a unique
system $(\xi_{F})_{F\in \on{Sch}}$, such that
$\xi_{F}\in M(F,F)$, and: if $x\in \on{Sch}(F,G)$,
then $i_{Q}(x)\circ \xi_{F} = \xi_{G}\circ i_{P}(x)$,
$\xi_{F\otimes G} = \xi_{F} \boxtimes \xi_{G}$, and
$\xi_{{\bf id}} = \xi$.
\end{prop}

{\em Proof.} The construction of $\xi_F$ is a generalization of that
of \cite{EH}, Section 1.11. If $Z\in \on{Irr(Sch)}$ has degree $n$,
we denote by $pr_Z \in \on{Sch}(T_n,Z)$ and $inj_Z\in \on{Sch}(Z,T_n)$
morphisms such that $pr_Z \circ inj_Z = \on{id}_Z$. Then we set
$\xi_Z := i_Q(pr_Z)\circ \xi^{\boxtimes n}\circ i_P(inj_Z)$.
If $F\in \on{Sch}$ decomposes as $F = \oplus_{Z\in \on{Irr(Sch)}}
F_Z \otimes Z$, where $F_Z$ is a multiplicity space, then
$\xi_F := \sum_{Z\in \on{Irr(Sch)}} \on{id}_{F_Z} \otimes \xi_Z
\in \oplus_{Z\in \on{Irr(Sch)}} \on{End}(F_Z) \otimes M(Z,Z)
\subset M(F,F)$. \hfill \qed \medskip

\subsection{The dual prop, (anti)automorphisms of a prop}

If $P$ is a prop, define its dual prop $P^*$ by $P^*(F,G) : = P(G^*,F^*)$.
For $x\in P(F,G)$, we define $x^*\in P^*(G^*,F^*)$ as the element
$x\in P(F,G)$. The prop operations of $P^*$ are defined  by
$y^* \circ x^* := (x\circ y)^*$, $y^* \boxtimes x^* := (y\boxtimes x)^*$,
$i_{P^*} = i_P \circ t$, where $t : \on{Sch} \to \on{Sch}^*$ is
induced by transposition; more precisely, for $F,G\in \on{Sch}$, 
$t_{F,G} : \on{Sch}(F,G) \to \on{Sch}^*(F,G) = \on{Sch}(G^*,F^*)$
is the map $\oplus_{Z\in\on{Irr(Sch)}}\on{Hom}(F_Z,G_Z) \to 
\oplus_{Z\in\on{Irr(Sch)}}\on{Hom}(G_Z^*,F_Z^*)$ induced by transposition.

A prop morphism $f : P \to Q$ induces a dual morphism
$f^* : P^* \to Q^*$, defined by $f^*(x^*) = f(x)^*$, with
$(g\circ f)^* = g^* \circ f^*$. We denote by $\on{Aut}(P) = \on{Aut}_+(P)$
the group of prop automorphisms of $P$, and by $\on{Aut}_-(P):=
\on{Iso}(P,P^*)$ the set of prop isomorphisms from $P$ to $P^*$.
Set $\on{Aut}_\pm(P):= \on{Aut}_+(P) \sqcup \on{Aut}_-(P)$. Then elements of
$\on{Aut}_\pm(P)$ can be composed by $g \bullet f := g\circ f$
if $f\in \on{Aut}_+(P)$, $g \bullet f := g^*\circ f$
if $f\in \on{Aut}_-(P)$. Then $\on{Aut}_\pm(P)$ is a group, and we have
an exact sequence $1\to \on{Aut}_+(P) \to \on{Aut}(P) \to \{\pm1\}\to 1$.

\subsection{Biprops}

We define $\on{Sch}_{1+1}$ as the symmetric tensor category whose objects are
Schur bifunctors, i.e., objects of $\on{Sch}_2$, and
$$
\on{Sch}_{1+1}(F\ul\boxtimes G,F'\ul\boxtimes G') :=
\on{Sch}(F,F') \otimes \on{Sch}^{*}(G,G') \simeq
\on{Sch}(F,F') \otimes \on{Sch}(G,G').
$$
as transposition gives rise to a prop isomorphism $\on{Sch}\simeq \on{Sch}^{*}$.  
We denote $\ul\boxtimes : (\on{Sch})^2 \to \on{Sch}_{1+1}$
the natural functor, so the object $F\boxtimes G$ of $\on{Sch}_2$
is denoted $F\ul\boxtimes G$ when viewed as an object of $\on{Sch}_{1+1}$
(and $\ul\boxtimes: \on{Sch}(F,F') \otimes \on{Sch}(G,G') \to
\on{Sch}_{1+1}(F\ul\boxtimes G,F'\ul\boxtimes G')$ is the natural map).
The tensor product is given by $(F\ul\boxtimes G)\otimes (F'\ul\boxtimes G')
:= (F\ul\boxtimes F')\ul\boxtimes (G\ul\boxtimes G')$. One defines 
similarly a tensor category ${\bf Sch}_{1+1}$, where the objects are those of
${\bf Sch}_2$. We define functors $\ul\Delta,\ul\Delta_{\otimes} : \on{Sch} \to
\on{Sch}_{1+1}$, by $\ul\Delta(F) := ((V,W) \mapsto F(V\oplus W))$ and
$\ul\Delta_{\otimes}(F) := ((V,W) \mapsto F(V\otimes W))$.

A $\on{Sch}_{1+1}$-biprop is a symmetric tensor category $\Pi_{1+1}$,
equipped with a morphism $\on{Sch}_{1+1} \to \Pi_{1+1}$, inducing the identity
on objects. One defines similarly ${\bf Sch}_{1+1}$-biprops.

If $B\in \on{Ob}({\bf Sch}_{1+1})$ and 
$\Pi_{1+1}$ is a ${\bf Sch}_{1+1}$-biprop,
then $B(\Pi_{1+1})$ is a $\on{Sch}$-prop defined by
$B(\Pi_{1+1})(F,G):= \Pi_{1+1}(F \circ B,G\circ B)$.

\subsection{The biprop $P_2$}

Let $P$ be one of the props (or ${\bf Sch}$-prop) $\on{LBA}$, ${\bf LBA}$, 
$\on{Bialg}$, $\ul{\on{Bialg}}$. For $F,G,F',G'\in\on{Irr(Sch)}$, 
we set 
$$
P_{2}(F\ul\boxtimes G,F'\ul\boxtimes G'):= 
P^{\Sigma}(F\boxtimes G^{\prime*},F'\boxtimes G^{*}),
$$
where $\Sigma = \{(F,F'),(F,G^{*}),(G^{\prime*},G^{*})\}$. 
We extend this definition by linearity  to define $P_{2}(B,B')$
for $B,B'\in\on{Ob}(\on{Sch}_{1+1})$ ($\on{Ob}({\bf Sch}_{1+1})$
in the case of ${\bf LBA}$). 

More explicitly, let $C\to P$, $O\to P$ be the prop morphisms and 
$e\in P({\bf id},{\bf id})$ the element such that 
for any $(F_i)_{i\in I}$,
$(G_j)_{j\in J}$ in $\on{Irr}(\on{Sch})$, the map
$$
\oplus_{Z_{ij}\in \on{Irr(Sch)}} (\otimes_{i\in I} C(F_i,\otimes_{j\in J}
Z_{ij})) \otimes (\otimes_{j\in J} O(\otimes_{i\in I} Z_{ij},G_j))
\to P(\otimes_{i\in I} F_i, \otimes_{j\in J} G_j),
$$
$(\otimes_{i\in I} c_i)\otimes (\otimes_{j\in J} o_j) \mapsto
(\boxtimes_{j\in J} o_j) \circ \sigma_{I,J} \circ 
\boxtimes_{i\in I}(\boxtimes_{j\in J} e_{Z_{ij}}) 
\circ (\boxtimes_{i\in I} c_i)$ is a linear isomorphism
($\oplus$ is replaced by $\hat\oplus$ in the case of ${\bf LBA}$, 
$\ul{\on{Bialg}}$). Then for $F,G,F',G'\in\on{Irr(Sch)}$,
$P_2(F\ul\boxtimes G,F'\ul\boxtimes G') \subset
P(F \otimes G^{\prime*}, F'\otimes G^*)$
is given by 
\begin{align*}
& P_2(F\ul\boxtimes G,F'\ul\boxtimes G')
\simeq \oplus_{Z_{XY}\in\on{Irr(Sch)} | Z_{G'F'} = {\bf 1}}
\\ & C(F, Z_{FF'}\otimes Z_{FG}) \otimes
C(G^{\prime *},Z_{G'F'}\otimes Z_{G'G})
\otimes O(Z_{FF'} \otimes Z_{G'F'},F') \otimes
O(Z_{FG} \otimes Z_{G'G},G^*)
\end{align*}
as the sum of all summands where $Z_{G'F'} = {\bf 1}$;
this definition is extended by linearity.

In particular, $\ul{\on{Bialg}}_2(F\ul\boxtimes G,F'\ul\boxtimes G')$
is the closure of
$\on{Bialg}_2(F\ul\boxtimes G,F'\ul\boxtimes G')
\subset \on{Bialg}(F\boxtimes G^{\prime*},F'\boxtimes G^*)$
for the $(\on{id}_{\bf id}^{\on{Bialg}} - \eta\circ \varepsilon)$-adic
topology of the latter space.

A structure of biprop is defined on $P_{2}$ as follows. 
We have $\on{Ob}(\on{Sch}_{1+1}) \subset \on{Ob}(\on{Sch}_{(1+1)})$
(and a similar inclusion replacing $\on{Sch}$ by ${\bf Sch}$), and for 
any $B,B'\in \on{Ob}(\on{Sch}_{1+1})$ (resp., ${\bf B},{\bf B'}
\in\on{Ob}({\bf Sch}_{1+1})$) we have inclusions 
$$
\on{LBA}_{2}(B,B')\subset \Pi(B,B'), \; 
{\bf LBA}_{2}(B,B')\subset \hbox{\boldmath$\Pi$\unboldmath}
({\bf B},{\bf B'}),
$$
$$ 
\on{Bialg}_{2}(B,B') \subset \tilde\Pi(B,B'), \; 
\ul{\on{Bialg}}_{2}(B,B') \subset \ul{\tilde\Pi}(B,B'). 
$$
The composition of the diagrams $\{(F,F'),(F,G^{*}),(G^{\prime*},G^{*})\}$
and $\{(F',F''),(F',G^{\prime*}),(G^{\prime\prime*},G^{\prime*})\}$
does not give rise to any loop, so the composition is well defined on 
$P_{2}(B,B')\otimes P_{2}(B',B'')$, and one checks that it takes its values in 
$P_{2}(B,B'')$. The external product also restricts to $P_{2}$. This defines a 
biprop structure on $P_{2}$. 

If $P,Q$ are any props, define the biprop $P\ul\boxtimes Q$ by
$(P\ul\boxtimes Q)(F\ul\boxtimes G,F'\ul\boxtimes G') :=
P(F,F') \otimes Q(G,G')$, for $F,...,G'\in \on{Irr(Sch)}$,
which we extend by linearity.

Then for $P\in\{\on{LBA},{\bf LBA},\on{Bialg},\ul{\on{Bialg}}\}$, 
we have a biprop morphism $P\ul\boxtimes P^{*} \to P_{2}$, given by
$P(F,F') \otimes P^{*}(G,G') \subset P_{2}(F\ul\boxtimes F',G\ul\boxtimes G')$
for any $F,...,G'\in \on{Sch}$.

When $P = S({\bf LBA})$, we set $P_2 := S^{\ul\boxtimes 2}({\bf LBA}_2)$.

\section{Compatibility of quantization functors with duality}

We first prove that the Etingof-Kazhdan functors are ``almost compatible with duality'', 
i.e., the duality diagram commutes up to conjugation by an inner automorphism. 
We then prove that any such functor is equivalent to a functor (i.e., can be transformed
into such a functor by composition with an inner automorphism) which is 
compatible with duality (i.e., for which the duality diagram commutes). 

\subsection{Almost compatibility with duality}

The dihedral group $D_4$ can be presented as the group with generators
$op$, $cop$ and $*$ and relations $op^2 = cop^2 = (op,cop)=*^2 = 1$,
$* \cdot cop \cdot * = op$, $* \cdot op \cdot * = cop$ (the product in
$D_4$ is denoted by $\cdot$ and $(a,b) = a\cdot b \cdot a^{-1} \cdot b^{-1}$).
We set $a_{1}...a_{n}:= a_{n}\cdot ... \cdot a_{1}$. So we have e.g.
$*cop = cop \cdot *$.

We have two morphisms $\eps,\eps' : D_{4}\to \ZZ/2\ZZ$, where
$\eps$ is defined by  $op, cop,*\mapsto -1$ and $\eps'$ by 
$op, cop\mapsto 1$, $*\mapsto -1$. Then $\on{Ker}(\eps) \simeq \ZZ/4\ZZ$, 
while $\on{Ker}(\eps')\simeq (\ZZ/2\ZZ)^{2}$. 
 
We have a unique morphism $D_{4}\to \on{Aut}_{\pm}(\on{LBA})$, 
$\theta\mapsto\theta_{\on{LBA}}$
such that $op_{\on{LBA}} : (\mu,\delta) \mapsto (-\mu,\delta)$, 
$cop_{\on{LBA}} : (\mu,\delta) \mapsto (\mu,-\delta)$ and 
$*_{\on{LBA}} : (\mu,\delta) \mapsto (\delta^{*},\mu^{*})$. 
It is such that the diagram 
$$
\xymatrix 
{
D_4 \ar[rd]^{\eps'} \ar[rr]_{} & &
\on{Aut}_\pm(\on{LBA})\ar[dl]_{}
\\ & \ZZ/2\ZZ & 
}
$$
commutes. We also have a unique morphism 
$D_{4}\to \on{Aut}_{\pm}(\on{Bialg})$, $\theta
\mapsto\theta_{\on{Bialg}}$
such that $op_{\on{Bialg}} : (m,\Delta) \mapsto (m\circ (21),\Delta)$, 
$cop_{\on{Bialg}} : (m,\Delta) \mapsto (m,(21)\circ \Delta)$ and 
$*_{\on{Bialg}} : (m,\Delta) \mapsto (\Delta^{*},m^{*})$. 
It is such that the diagram 
$$
\xymatrix 
{
D_4 \ar[rd]^{\eps'} \ar[rr]_{} & &
\on{Aut}_\pm(\on{Bialg})\ar[dl]_{}
\\ & \ZZ/2\ZZ & 
}
$$
commutes. 

Denote by $\on{Assoc}(\kk)$ the set of associators defined over $\kk$. This
set consists of series $\Phi(A,B)$ in noncommutative variables $A,B$.
It is equipped with an action of $\{\pm 1\}$, where $(-1) \cdot \Phi = \Phi'$
given by $\Phi'(A,B) = \Phi(-A,-B)$. (The fixed points of this action are
called the even associators.)

\begin{thm} \label{thm:duality}
For each $\theta\in D_{4}$, there exists $\xi_{\theta}$, where
$\xi_{\theta}\in S({\bf LBA})({\bf id},{\bf id})^{\times}$ if
$\eps'(\theta)=1$ and $\xi_{\theta}\in S({\bf LBA}^{*})
({\bf id},{\bf id})^{\times}$ if $\eps'(\theta)=-1$, such that: 
$$
\theta_{S({\bf LBA})} \circ Q_{\Phi} = \on{Inn}(\xi_{\theta})
\circ Q_{\eps(\theta)\cdot \Phi} \circ \theta_{\on{Bialg}}\; \on{if}\;
\eps'(\theta)=1, 
$$
$$
\theta_{S({\bf LBA})} \circ Q_{\Phi} = \on{Inn}(\xi_{\theta})
\circ Q^{*}_{\eps(\theta)\cdot \Phi} \circ \theta_{\on{Bialg}}\; \on{if}\;
\eps'(\theta)=-1.  
$$
\end{thm}

We say that a quantization functor $Q$ is almost compatible with duality if 
for some $\xi\in S({\bf LBA}^{*})({\bf id},{\bf id})^{\times}$, 
\begin{equation} \label{id:quasicomp}
*cop_{S({\bf LBA})} \circ Q = \on{Inn}(\xi)\circ 
Q^{*}\circ *cop_{\on{Bialg}}.
\end{equation} 
So each $Q_{\Phi}$ is almost compatible with duality. 

{\em Proof.} The subset of $D_4$ of elements $\theta$ for which the 
result holds is a subgroup. Since $D_{4}$ is generated by $cop$ and $*op$, 
it suffices to prove it for these elements. When $\theta=cop$, 
the result was proved in \cite{EH}. So we now prove it for $\theta = *op 
= (*cop)^{-1}$. 
We set in this section, $Q:= Q_{\Phi}$, $\xi:= \xi_{*op}$, so we should find
$\xi$ such that 
$$
*op_{S({\bf LBA})} \circ Q = \on{Inn}(\xi)\circ 
Q^{*}\circ *op_{\on{Bialg}}.
$$
This is written as follows
$$
*op_{S({\bf LBA})}(Q(m)) = \xi \circ 
(Q(\Delta))^{*}\circ (\xi^{-1})^{\boxtimes 2}, \quad 
*op_{S({\bf LBA})}(Q(\Delta)) = \xi^{\boxtimes 2} \circ (21)\circ
Q(m)^{*}\circ \xi^{-1},  
$$
i.e., applying $x\mapsto x^{*}$, as follows
\begin{equation} \label{duality:result}
Q(m)^{\tau} = \tilde\xi^{\boxtimes 2} \circ Q(\Delta) \circ
\tilde\xi^{-1}, \quad 
Q(\Delta)^{\tau} = \tilde\xi\circ Q(m)\circ (21) \circ 
(\tilde\xi^{-1})^{\boxtimes 2}, 
\end{equation}
where $x^{\tau} = (*op(x))^{*}$ and $\tilde\xi = (\xi^{*})^{-1}$. 

The map $x\mapsto x^{\tau}:= (*op(x))^{*}$ is a linear map ${\bf LBA}(A,B)\to 
{\bf LBA}(B^{*},A^{*})$, such that $(y\circ x)^{\tau} = x^{\tau}
\circ y^{\tau}$. It is uniquely defined by this condition, the assignments
$(\mu,\delta)\mapsto (\delta,-\mu)$, and 
$(x\boxtimes y)^{\tau} = x^{\tau}\boxtimes y^{\tau}$, 
$i_{{\bf LBA}}(z)^{\tau}=i_{{\bf LBA}}(z^{\tau})$ for 
$z\in {\bf Sch}(A,B)$, where $z\mapsto z^{\tau}$ denotes the 
map ${\bf Sch}(A,B)\to {\bf Sch}(B^{*},A^{*})$ induced by the 
transposition. 

We then construct a linear map 
$$
\hbox{\boldmath$\Pi$\unboldmath}
(F\underline\boxtimes G,F'\underline\boxtimes G') 
\to\hbox{\boldmath$\Pi$\unboldmath}
(G\underline\boxtimes F,G'\underline\boxtimes F'), \quad  
x\mapsto x^{\tau},
$$ 
induced by the isomorphisms 
$\hbox{\boldmath$\Pi$\unboldmath}
(F\underline\boxtimes G,F'\underline\boxtimes G') \simeq
{\bf LBA}((\otimes F) \otimes (\otimes G')^*,(\otimes F')
\otimes (\otimes G)^*)$, 
$\hbox{\boldmath$\Pi$\unboldmath}(G\underline\boxtimes F,$ 
$G'\underline\boxtimes F') \simeq
{\bf LBA}((\otimes G) \otimes (\otimes F^{\prime *}),(\otimes G')
\otimes (\otimes F^{*}))$, and the map $x\mapsto x^{\tau}$ on 
${\bf LBA}$-spaces. We then have (for $x,y$ in 
\boldmath$\Pi$\unboldmath-spaces) $(y\circ x)^{\tau}
=y^{\tau}\circ x^{\tau}$, $(x\boxtimes y)^{\tau}=x^{\tau}
\boxtimes y^{\tau}$. 

We finally define a linear map 
$$
\hbox{\boldmath$\Pi$\unboldmath}
(F\underline\boxtimes G,F'\underline\boxtimes G') 
\to\hbox{\boldmath$\Pi$\unboldmath}
(F^{\prime*}\underline\boxtimes G^{\prime*},F^*\underline\boxtimes G^*), \quad  
x\mapsto x^t := (x^\tau)^* = (x^*)^\tau; 
$$
then $(x\circ y)^t = y^t \circ x^t$. 

Let us denote by $x\mapsto \underline x$ the canonical map
${\bf LBA}(S^{\otimes p},S^{\otimes q}) \to
\hbox{\boldmath$\Pi$\unboldmath}
((S\underline\boxtimes {\mathfrak 1})^{\otimes p},
(S\underline\boxtimes {\mathfrak 1})^{\otimes q})$. Then one
checks that $\underline{x^{\tau}} = (\underline x)^t$.

Applying the operation $x\mapsto \ul{x}$ to (\ref{duality:result}), this
condition translates in terms of \boldmath$\Pi$\unboldmath-spaces as 
$$
m_{a}^t = \ul{\tilde\xi}^{\boxtimes 2} \circ \Delta_{a} \circ
\ul{\tilde\xi}^{-1}, \quad 
\Delta_{a}^t = \ul{\tilde\xi}\circ m_{a}\circ (21) \circ 
(\ul{\tilde\xi}^{-1})^{\boxtimes 2}. 
$$

Let us now compute $m_a^\tau$, $\Delta_a^\tau$.

\begin{lemma}
Define $\omega_S\in {\bf Sch}(S,S)^\times$ as $\hat\oplus_{n\geq 0}
(-1)^n \on{id}_{S^n}$.
Then
$$
m_\Pi^\tau = (\on{id}_S \underline\boxtimes \omega_S) \circ 
m_\Pi \circ (21)
\circ (\on{id}_S \underline\boxtimes \omega_S)^{\boxtimes 2}.
$$
\end{lemma}

{\em Proof of Lemma.}
Let $\A$ be a finite dimensional Lie bialgebra, and let
$\G$ be its double. Then $\G = \A\oplus \A^{*cop}$ (as Lie
coalgebras; $\A$ and $\A^{*cop}$ are also sub-Lie bialgebras of
$\G$). Then
$m_\Pi$ is the propic version of the map $(S(\A) \otimes S(\A^*))^{\otimes 2}
\to S(\A) \otimes S(\A^*)$ induced by the product of $U(\G)$
and the isomorphism $S(\A) \otimes S(\A^*) \to U(\G)$ given by
$x \otimes \xi \mapsto \on{sym}_\A(x)\on{sym}_{\A^*}(\xi)$, where
$\on{sym}_\x :
S(\x) \to U(\x)$ is the symmetrization map for the Lie algebra $\x$.

Let now $\G'$ be the double of $\A^{*cop}$. Then
$\G' = (\A^*)^{cop} \oplus (\A^{*cop})^{*cop}
= \A^{*cop} \oplus \A^{op,cop}$. Then
$m_\Pi^\tau$ is the propic version of the map
$(S(\A) \otimes S(\A^*))^{\otimes 2}
\stackrel{(m_\Pi)_\A}{\to} S(\A) \otimes S(\A^*)$ induced
by the product of $U(\G')$
and the isomorphism $S(\A) \otimes S(\A^*) \to U(\G')$,
$x\otimes \xi\mapsto \on{sym}_{\A^*}(\xi) \on{sym}_{\A^{{op}}}(x)$.

An isomorphism $\iota : \G\to \G'$ is given by $(a,\ell) \mapsto (\ell,-a)$.
Then the map $S_{U(\G')} \circ U(\iota) : U(\G) \to U(\G')$ is a linear
isomorphism, where $S_{U(\G')}$ is the antipode of $U(\G')$.
Since the diagram
$$
\xymatrix{
S(\A) \ar[r]^{\on{sym}_{\A}} \ar[d]_{S(-\on{id}_{\A})}& U(\A)
\ar[d]^{\scriptstyle{U(\iota)}}\\
S(\A) \ar[r]^{\on{sym}_{\A^{{op}}}}& U(\A^{{op}})}
$$
commutes, the composed map $S(\A) \otimes S(\A^*) \to U(\G) 
\simeq U(\G') \leftarrow
S(\A) \otimes S(\A^*)$ is $x\otimes \xi \mapsto x\otimes \bar\xi$
(where $v\mapsto \bar v$ denotes the automorphism $S(-\on{id}_V)$
of $S(V)$).

In the commutative diagram
$$
\xymatrix@!0 @R=4ex @C=21ex{
x_1\otimes \xi_1 \otimes x_2 \otimes \xi_2 
& & & 
x'_1\otimes \xi'_1 \otimes x'_2 \otimes \xi'_2 
\\ 
\bigcap \mkern-13mu | &&& \bigcap \mkern-13mu | \\
(S(\A) \otimes S(\A^*))^{\otimes 2} \ar[r] \ar[ddd]_{\scriptstyle{(m_\Pi)_\A}}
& U(\G)^{\otimes 2} \ar@{<->}[r]^\simeq \ar[ddd]_{\scriptstyle{\on{prod}}} &
U(\G')^{\otimes 2} \ar[ddd]_{\scriptstyle{\on{prod}}} & \ar[l] (S(\A)\otimes
S(\A^*))^{\otimes 2} \ar[ddd]^{\scriptstyle{(m_\Pi^\tau)_\A}} 
\\
\\
\\
S(\A) \otimes S(\A^*) \ar[r] & U(\G)\ar@{<->}[r]^\simeq &
U(\G') &\ar[l]  S(\A)\otimes S(\A^*)  }
$$
the image of $x_1\otimes \xi_1 \otimes x_2 \otimes \xi_2$ in $U(\G')$ is
$\on{sym}_{\A^*}(\bar\xi_2) \on{sym}_{\A^{{op}}}(x_2)
\on{sym}_{\A^*}(\bar\xi_1) \on{sym}_{\A^{{op}}}(x_1)$, while the image of
$x'_1\otimes \xi'_1 \otimes x'_2 \otimes \xi'_2$ in $U(\G')$ is
$\on{sym}_{\A^*}(\xi'_1) \on{sym}_{\A^{{op}}}(x'_1)
\on{sym}_{\A^*}(\xi'_2)
\on{sym}_{\A^{{op}}}(x'_2)$.

It follows that the diagram
$$
\xymatrix{
(S(\A) \otimes S(\A^*))^{\otimes 2} \ar[r] \ar[d]_{\scriptstyle{(m_\Pi)_\A}}
&
(S(\A) \otimes S(\A^*))^{\otimes 2}
\ar[d]^{\scriptstyle{(m_\Pi^\tau)_{\A}}} \\
S(\A) \otimes S(\A^*) \ar[r] & S(\A) \otimes S(\A^*)
}
$$
commutes, where the horizontal maps are $x_1\otimes \xi_1 \otimes x_2
\otimes
\xi_2 \mapsto x_2 \otimes \bar\xi_2 \otimes x_1 \otimes \bar\xi_1$
and $x\otimes \xi \mapsto x\otimes \bar\xi$.

The statement of the lemma is a propic version of the commutativity of this
diagram; one shows that the proof can be carried to the propic setting.
\hfill \qed \medskip

\begin{lemma}
$\Delta_0^\tau = \Delta_0$.
\end{lemma}

{\em Proof.} $\Delta_0\in\hbox{\boldmath$\Pi$\unboldmath}
(S\underline\boxtimes S,(S\underline\boxtimes
S)^{\otimes 2}) \subset {\bf LBA}(S \otimes S^{\otimes 2},S^{\otimes 2}
\otimes S)$ is $\Delta^S \boxtimes m^S$, where $\Delta^S\in
{\bf Sch}(S,S^{\otimes 2})$ and $m^S \in {\bf Sch}(S^{\otimes 2},S)$
are the propic versions of the coproduct and product of symmetric algebras.
The statement then follows from the fact that 
$m^S$ and $\Delta^S$ are interchanged by the maps 
${\bf Sch}(S,S^{\otimes 2}) \to
{\bf Sch}(S^{\otimes 2},S)$
and ${\bf Sch}(S^{\otimes 2},S) \to
{\bf Sch}(S,S^{\otimes 2})$ defined by $x\mapsto x^\tau$.
\hfill \qed \medskip

\begin{lemma}
For $Y$ in the image of $\t_n \subset U(\t_n) \to {\bf U}_n
= \hbox{\boldmath$\Pi$\unboldmath}
({\mathfrak 1}\underline\boxtimes {\mathfrak 1},
(S\underline\boxtimes S)^{\otimes n})$,
we have $Y^\tau = -(\on{id}_S \underline\boxtimes \omega_S)^{\boxtimes n}
\circ Y$.
\end{lemma}

{\em Proof.} We prove this by induction on the degree of $Y$
w.r.t. the generators $t_{ij}$ of $\t_n$. If $Y = t_{ij}$,
then $Y^\tau = Y$, while
$(\on{id}_S \underline\boxtimes \omega_S)^{\boxtimes n}
\circ Y = -Y$. Let us assume that $Y = [Y',Y''] =
m_\Pi^{\boxtimes n} \circ \sigma_{n,2}(S\ul\boxtimes S)^{\Pi} 
\circ (Y' \boxtimes Y''- Y'' \boxtimes Y')$. Then
\begin{align*}
& Y^\tau = (m_\Pi^\tau)^{\boxtimes n} \circ  
\sigma_{n,2}(S\ul\boxtimes S)^{\Pi} \circ
(Y^{\prime\tau} \boxtimes Y^{\prime\prime\tau} - Y^{\prime\prime\tau} 
 \boxtimes Y^{\prime\tau})
\\ &
=
\big( (\on{id}_S \underline\boxtimes \omega_S) \circ m_\Pi \circ (21)
\circ (\on{id}_S \underline\boxtimes \omega_S)^{\boxtimes 2} 
\big)^{\boxtimes n}
\circ  \sigma_{n,2}(S\ul\boxtimes S)^{\Pi}
\\ & \circ
\Big( \big( (\on{id}_S \underline\boxtimes \omega_S)^{\boxtimes n}
\circ Y' \big)
\boxtimes
\big( (\on{id}_S \underline\boxtimes \omega_S)^{\boxtimes n}
\circ Y'' \big)
- \big( (\on{id}_S \underline\boxtimes \omega_S)^{\boxtimes n}
\circ Y'' \big)
\boxtimes
\big( (\on{id}_S \underline\boxtimes \omega_S)^{\boxtimes n}
\circ Y' \big) \Big)
\\ &
= \big( (\on{id}_S \underline\boxtimes \omega_S) \circ m_\Pi
\big)^{\boxtimes n} \circ  \sigma_{n,2}(S\ul\boxtimes S)^{\Pi}
 \circ (Y''\boxtimes Y' - Y'\boxtimes Y'' )
= - (\on{id}_S \underline\boxtimes \omega_S)^{\boxtimes n} \circ Y.
\end{align*}
\hfill \qed \medskip

\begin{lemma}
For $X$ in the image of $\on{exp}(\t_n) \subset \wh{U(\t_n)} \to {\bf U}_n
= \hbox{\boldmath$\Pi$\unboldmath}
({\mathfrak 1}\underline\boxtimes {\mathfrak 1},
(S\underline\boxtimes S)^{\otimes n})$,
we have $X^\tau = (\on{id}_S \underline\boxtimes \omega_S)^{\boxtimes n}
\circ X^{-1}$.
\end{lemma}

{\em Proof.} Let $Y\in \t_n$ and $X := \on{exp}(Y)$. Then
$X = \sum_{k\geq 0} (k!)^{-1} (m_\Pi^{(n)})^{\boxtimes k} \circ 
 \sigma_{n,k}(S\ul\boxtimes S)^{\Pi} \circ
Y^{\boxtimes k}$. Then
\begin{align*}
& X^\tau =
\sum_{k\geq 0} (k!)^{-1} (m_\Pi^{(k)\tau})^{\boxtimes n}
\circ  \sigma_{n,k}(S\ul\boxtimes S)^{\Pi}
\circ (Y^\tau)^{\boxtimes k}
\\ &
= \sum_{k\geq 0} (k!)^{-1} \Big(
(\on{id}_S \underline\boxtimes \omega_S) \circ m_\Pi^{(k)} \circ (k,...,2,1)
\circ (\on{id}_S \underline\boxtimes \omega_S)^{\boxtimes k}
\Big)^{\boxtimes n} \circ \sigma_{n,k}(S\ul\boxtimes S)^{\Pi} \circ
\big(- (\on{id}_S \underline\boxtimes \omega_S)^{\boxtimes n} \circ Y
\big)^{\boxtimes k}
\\ &
=
\sum_{k\geq 0} (k!)^{-1} \Big(
(\on{id}_S \underline\boxtimes \omega_S) \circ m_\Pi^{(k)}
\Big)^{\boxtimes n} \circ \sigma_{n,k}(S\ul\boxtimes S)^{\Pi}\circ
\big(- Y
\big)^{\boxtimes k} =
(\on{id}_S \underline\boxtimes \omega_S)^{\boxtimes n}
\circ X^{-1}.
\end{align*}
\hfill \qed \medskip

\begin{lemma}
There exists $u\in {\bf U}_1^\times$, such that
$(\on{id}_S \underline\boxtimes
\omega_S)^{\boxtimes 2} \circ \on{J}^\tau = u^{12} * 
(\on{J}^{2,1})^{-1} * (u^1 *u^2)^{-1}$.
\end{lemma}

{\em Proof.} $\on{J}$ satisfies $\on{J}^{1,2} \star \on{J}^{12,3} =
\on{J}^{2,3} \star \on{J}^{1,23} \star \Phi$, which is rewritten as
$$
m_\Pi^{\boxtimes 3} \circ \sigma_{3,2}(S\ul\boxtimes S)^{\Pi} 
\circ (\on{J}^{1,2} \boxtimes \on{J}^{1,23})
=
(m_\Pi^{(3)})^{\boxtimes 3} \circ \sigma_{3,3}(S\ul\boxtimes S)^{\Pi}
\circ (\on{J}^{2,3}\boxtimes \on{J}^{12,3} \boxtimes \Phi).
$$
Applying $x\mapsto x^\tau$, we get
$$
(m_\Pi^\tau)^{\boxtimes 3} \circ \sigma_{3,2}(S\ul\boxtimes S)^{\Pi}
\circ ((\on{J}^{1,2})^\tau\boxtimes (\on{J}^{1,23})^\tau)
=
(m_\Pi^{(3)\tau})^{\boxtimes 3} \circ 
 \sigma_{3,3}(S\ul\boxtimes S)^{\Pi} \circ
\Big((\on{J}^{2,3})^{\tau}
\boxtimes (\on{J}^{12,3})^\tau \boxtimes 
\big( (\on{id}_S \underline\boxtimes
\omega_S)^{\boxtimes 3} \circ \Phi^{-1} \big) \Big).
$$
So
\begin{align*}
& \Big( (\on{id}_S\underline\boxtimes\omega_S) \circ m_\Pi\Big)^{\boxtimes 3}
\circ  \sigma_{3,2}(S\ul\boxtimes S)^{\Pi}\circ \Big(
\big( (\on{id}_S \underline\boxtimes \omega_S)^{\boxtimes 3} \circ
(\on{J}^{1,23})^\tau \big) \boxtimes
\big( (\on{id}_S \underline\boxtimes \omega_S)^{\boxtimes 3} \circ
(\on{J}^{1,2})^\tau \big) \Big)
\\ & =
( (\on{id}_S \underline\boxtimes \omega_S) \circ m_\Pi^{(3)})^{\boxtimes 3} 
\circ  \sigma_{3,3}(S\ul\boxtimes S)^{\Pi}\circ
\Big( \Phi^{-1} \boxtimes
\big( (\on{id}_S \underline\boxtimes
\omega_S)^{\boxtimes 3} \circ
(\on{J}^{12,3})^\tau \big) \boxtimes
\big( (\on{id}_S \underline\boxtimes
\omega_S)^{\boxtimes 3} \circ
(\on{J}^{2,3})^\tau \big) \Big).
\end{align*}
Simplifying $(\on{id}_S \underline\boxtimes \omega_S)^{\boxtimes 3}$, we get
$$
\big( (\on{id}_S \underline\boxtimes \omega_S)^{\boxtimes 2}
\circ \on{J}^\tau \big)^{1,23}
\star\big( (\on{id}_S \underline\boxtimes \omega_S)^{\boxtimes 2}
\circ \on{J}^\tau \big)^{2,3}
= \Phi^{-1}\star
\big( (\on{id}_S \underline\boxtimes \omega_S)^{\boxtimes 2}
\circ \on{J}^\tau \big)^{12,3}
\star\big( (\on{id}_S \underline\boxtimes \omega_S)^{\boxtimes 2}
\circ \on{J}^\tau \big)^{1,2}.
$$
So $\on{J}' := \big( (\on{id}_S \underline\boxtimes \omega_S)^{\boxtimes 2}
\circ \on{J}^\tau \big)^{-1}$ satisfies the same equation as $\on{J}$,
namely $\on{J}^{1,2} \star \on{J}^{12,3} =
\on{J}^{2,3} \star\on{J}^{1,23} \star\Phi$. We have $r^\tau = r^{2,1}$, therefore
$\on{J}' = 1 + r^{2,1}/2 + ...$, so according to \cite{E}, Remark 6.7
(see also \cite{E}, Thm. 2.1), $\on{J}'$ is gauge-equivalent to
$\on{J}^{2,1}$. \hfill \qed \medskip

\begin{lemma}
The elements $(\on{id}_S \underline\boxtimes \omega_S)^{\boxtimes 2}
\circ \on{J}^\tau$ and $(\on{id}_S \underline\boxtimes \omega_S)^{\boxtimes 2}
\circ (\on{J}^{-1})^\tau$ of ${\bf U}_2$ are inverse of each other.
\end{lemma}

{\em Proof.} We have $m_\Pi^{\boxtimes 2} \circ (1324) \circ (\on{J}
\boxtimes \on{J}^{-1}) = 1$. Applying $x\mapsto x^\tau$, we get
$(m_\Pi^\tau)^{\boxtimes 2} \circ (\on{J}^\tau \boxtimes
(\on{J}^{-1})^\tau) = 1$. Therefore
$((\on{id}_S \underline\boxtimes \omega_S) \circ m_\Pi)^{\boxtimes 2}
\circ (1324) \circ
\big( (\on{id}_S \underline\boxtimes \omega_S)^{\boxtimes 2}
\circ (\on{J}^{-1})^\tau
\boxtimes (\on{id}_S \underline\boxtimes \omega_S)^{\boxtimes 2}
\circ \on{J}^\tau \big) = 1$. Simplifying
$(\on{id}_S \underline\boxtimes \omega_S)^{\boxtimes 2}$, we get the result.
\hfill \qed \medskip

\begin{lemma}
$\on{R}^\tau = (\on{id}_S \underline\boxtimes \omega_S)^{\boxtimes 2}
\circ \big( u^1\star u^2 \star (\on{R}^{2,1})^{-1} \star(u^1\star
u^2)^{-1} \big)$.
\end{lemma}

{\em Proof.} We have $\on{R} = \on{J}^{2,1}\star e^{t/2}\star\on{J}^{-1}$, i.e.,
$\on{R} = \sum_{k\geq 0} (k!)^{-1} (m_\Pi^{(k+2)} \boxtimes m_\Pi^{(k+2)})
\circ \sigma_{2,k+2}(S\ul\boxtimes S)^{\Pi}
\circ (\on{J}^{2,1} \boxtimes (t/2)^{\boxtimes k} \boxtimes
\on{J}^{-1})$. Then
\begin{align*}
& \on{R}^\tau = \sum_{k\geq 0}(k!)^{-1}
(m_\Pi^{(k+2)\tau} \boxtimes m_\Pi^{(k+2)\tau})
\circ  \sigma_{2,k+2}(S\ul\boxtimes S)^{\Pi} 
\circ ( (\on{J}^{2,1})^\tau \boxtimes 
(t^\tau/2)^{\boxtimes k} \boxtimes (\on{J}^{-1})^\tau)
\\ & =
\sum_{k\geq 0}(k!)^{-1}
\big( (\on{id}_S \underline\boxtimes \omega_S) \circ
m_\Pi^{(k+2)} \circ (k+2,...,1) \circ
(\on{id}_S \underline\boxtimes \omega_S)^{\boxtimes k+2} \big)^{\boxtimes 2}
\\ & \circ  \sigma_{2,k+2}(S\ul\boxtimes S)^{\Pi}
\circ ((\on{J}^{2,1})^\tau \boxtimes (t/2)^{\boxtimes k}
\boxtimes (\on{J}^{-1})^\tau)
\\ &
= \sum_{k\geq 0}(k!)^{-1}
\big( (\on{id}_S \underline\boxtimes \omega_S) \circ
m_\Pi^{(k+2)} \big)^{\boxtimes 2}
\circ  \sigma_{2,k+2}(S\ul\boxtimes S)^{\Pi}
\\ & \circ \Big(
\big( (\on{id}_S \underline\boxtimes \omega_S)^{\boxtimes 2} \circ
(\on{J}^{-1})^\tau \big)
\boxtimes
(-t/2)^{\boxtimes k}
\boxtimes
\big( (\on{id}_S \underline\boxtimes \omega_S)^{\boxtimes 2} \circ
(\on{J}^{2,1})^\tau \big)
\Big)
\\ &
= \sum_{k\geq 0}(k!)^{-1}
\big( (\on{id}_S \underline\boxtimes \omega_S) \circ
m_\Pi^{(k+2)} \big)^{\boxtimes 2}
\circ  \sigma_{2,k+2}(S\ul\boxtimes S)^{\Pi}
\\ & \circ \Big(
(u^1\star u^2 *\star\on{J}^{2,1} \star(u^{12})^{-1})
\boxtimes
(-t/2)^{\boxtimes k}
\boxtimes
(u^{12} \star\on{J}^{-1}\star(u^1\star u^2)^{-1})
\Big)
\\ &
= (\on{id}_S \underline\boxtimes \omega_S)^{\boxtimes 2} \circ
(u^1\star u^2 \star\on{J}^{2,1} \star(u^{12})^{-1} \star e^{-t/2}\star 
u^{12} \star \on{J}^{-1} 
\star (u^1\star u^2)^{-1})
\\ & = (\on{id}_S \underline\boxtimes \omega_S)^{\boxtimes 2} \circ
(u^1\star u^2 \star (\on{R}^{2,1})^{-1} \star (u^1\star u^2)^{-1}).
\end{align*}
\hfill \qed \medskip

Recall that for some $\sigma\in \hbox{\boldmath$\Pi$\unboldmath}
({\mathfrak 1}\underline\boxtimes S,
{\mathfrak 1} \underline\boxtimes S)^\times$, we have
$\on{R}^{-1} = (\on{R}_+ \boxtimes (\on{R}_- \circ \sigma))
\circ \on{can}_{S\underline\boxtimes {\mathfrak 1}}$. Then
$(\on{R}^{-1})^{2,1} = ((\on{R}_- \circ \sigma) \boxtimes \on{R}_+)
\circ \on{can}_{{\mathfrak 1}\underline\boxtimes S}$. Therefore
$$
\on{R}^\tau = (\on{id}_S \underline\boxtimes \omega_S)^{\boxtimes 2}
\circ \on{Ad}(u)^{\boxtimes 2} \circ
((\on{R}_- \circ \sigma) \boxtimes \on{R}_+)
\circ \on{can}_{{\mathfrak 1}\underline\boxtimes S},
$$
where $\on{Ad}(u) = m_\Pi^{(3)} \circ
(u \boxtimes \on{id}_{S\underline\boxtimes
{\mathfrak 1}} \boxtimes u^{-1})\in \hbox{\boldmath$\Pi$\unboldmath}
(S\underline\boxtimes {\mathfrak 1},
S\underline\boxtimes {\mathfrak 1})^\times$.
On the other hand, $\on{R} = (\on{R}_+ \boxtimes \on{R}_-) \circ
\on{can}_{S\underline\boxtimes {\mathfrak 1}}$ implies that
$$
\on{R}^\tau = (\on{R}_+^\tau \boxtimes \on{R}_-^\tau) \circ
\on{can}^\tau_{S\underline\boxtimes{\mathfrak 1}}
=(\on{R}_+^\tau \boxtimes \on{R}_-^\tau) \circ
\on{can}_{{\mathfrak 1}\underline\boxtimes S}.
$$
Comparing these formulas, we get the existence of $\psi\in
\hbox{\boldmath$\Pi$\unboldmath}({\mathfrak
1}\underline\boxtimes S,{\mathfrak 1} \underline\boxtimes S)^\times$, such that
$$
\on{R}_+^\tau =
(\on{id}_S \underline\boxtimes \omega_S) \circ \on{Ad}(u) \circ \on{R}_-
\circ \sigma \circ \psi^{-1}, \quad
\on{R}_-^\tau = (\on{id}_S \underline\boxtimes \omega_S)
\circ \on{Ad}(u) \circ \on{R}_+ \circ\psi^*.
$$

It then follows that
$$
(\on{R}_+^{(-1)})^\tau = \psi\circ\sigma^{-1} \circ \on{R}_-^{(-1)}
\circ \on{Ad}(u)^{-1} \circ (\on{id}_S \underline\boxtimes \omega_S),
\quad
(\on{R}_-^{(-1)})^\tau = (\psi^*)^{-1} \circ \on{R}_+^{(-1)}
\circ \on{Ad}(u)^{-1} \circ (\on{id}_S \underline\boxtimes \omega_S).
$$

\begin{lemma}
$(\on{Ad(J)})^\tau = (\on{id}_S \underline\boxtimes \omega_S)^{\boxtimes 2}
\circ \on{Ad}(u^1\star u^2 \star\on{J}^{2,1}\star(u^{12})^{-1})
\circ (\on{id}_S \underline\boxtimes \omega_S)^{\boxtimes 2}$.
\end{lemma}

{\em Proof.} We have
\begin{align*}
& (\on{Ad(J)})^\tau = \big( (m_\Pi^{(3)} \boxtimes m_\Pi^{(3)})
\circ  \sigma_{3,3}(S\ul\boxtimes S)^{\Pi} \circ (\on{J}
\boxtimes \on{id}_{S\underline\boxtimes S} \boxtimes \on{J}^{-1})\big)^\tau
\\ &
= \big( (\on{id}_S \underline\boxtimes \omega_S) \circ m_\Pi^{(3)} \circ (321)
\circ (\on{id}_S \underline\boxtimes \omega_S)^{\boxtimes 3}\big)^{\boxtimes 2}
\circ   \sigma_{3,3}(S\ul\boxtimes S)^{\Pi} \circ (\on{J}^\tau
\boxtimes \on{id}_{S\underline\boxtimes S} \boxtimes (\on{J}^{-1})^\tau)
\\ &
= \big( (\on{id}_S \underline\boxtimes \omega_S) \circ m_\Pi^{(3)}
\big)^{\boxtimes 2}
\circ   \sigma_{3,3}(S\ul\boxtimes S)^{\Pi}\circ \Big(
\big( (\on{id}_S \underline\boxtimes \omega_S)^{\otimes 2}
\circ \on{J}^\tau \big)
\boxtimes (\on{id}_S \underline\boxtimes \omega_S)^{\otimes 2}
\boxtimes \big( (\on{id}_S \underline\boxtimes \omega_S)^{\otimes 2}
\circ (\on{J}^{-1})^\tau \big)
\Big)
\\ &
= \big( (\on{id}_S \underline\boxtimes \omega_S) \circ m_\Pi^{(3)}
\big)^{\boxtimes 2}
\circ   \sigma_{3,3}(S\ul\boxtimes S)^{\Pi} \\ & \circ \big(
\{u^1\star u^2\star \on{J}^{2,1} \star(u^{12})^{-1}\}
\boxtimes (\on{id}_S \underline\boxtimes \omega_S)^{\otimes 2}
\boxtimes
\{u^{12} \star (\on{J}^{2,1})^{-1} \star(u^1\star u^2)^{-1}\} \big)
\\ &
=
(\on{id}_S \underline\boxtimes \omega_S)^{\boxtimes 2}
\circ \on{Ad}(u^1\star u^2\star\on{J}^{2,1}\star(u^{12})^{-1})
\circ (\on{id}_S \underline\boxtimes \omega_S)^{\boxtimes 2}.
\end{align*}
\hfill \qed \medskip

\begin{lemma}
$\Delta_\Pi^\tau = (\on{id}_S \underline\boxtimes \omega_S)^{\boxtimes 2}
\circ \on{Ad}(u)^{\boxtimes 2} \circ (21) \circ \Delta_\Pi \circ
\on{Ad}(u)^{-1}
\circ (\on{id}_S \underline\boxtimes \omega_S)$.
\end{lemma}

{\em Proof.}
We have
\begin{align*}
& \Delta_\Pi^\tau = \on{Ad(J)}^\tau \circ \Delta_0^\tau
= (\on{id}_S \underline\boxtimes \omega_S)^{\boxtimes 2} \circ
\on{Ad}(u^1\star u^2\star \on{J}^{2,1} \star(u^{12})^{-1}) \circ
(\on{id}_S \underline\boxtimes
\omega_S)^{\boxtimes 2} \circ \Delta_0
\\ & =
 (\on{id}_S \underline\boxtimes \omega_S)^{\boxtimes 2} \circ
\on{Ad}(u^1\star u^2 \star\on{J}^{2,1} \star(u^{12})^{-1}) \circ
 \Delta_0 \circ (\on{id}_S \underline\boxtimes
\omega_S)
\\ & =
(\on{id}_S \underline\boxtimes \omega_S)^{\boxtimes 2}
\circ \on{Ad}(u)^{\boxtimes 2} \circ (21) \circ \Delta_\Pi \circ
\on{Ad}(u)^{-1}
\circ (\on{id}_S \underline\boxtimes \omega_S).
\end{align*}
\hfill \qed \medskip

Now
\begin{align*}
& m_a^\tau = 
(\on{R}_+^{(-1)} \circ m_\Pi \circ \on{R}_+^{\boxtimes 2})^\tau
=
(\on{R}_+^{(-1)})^\tau  \circ m_\Pi^\tau \circ 
(\on{R}_+^\tau)^{\boxtimes 2}
\\ &
= \psi\circ \sigma^{-1} \circ \on{R}_-^{(-1)} \circ \on{Ad}(u)^{-1}
\circ m_\Pi \circ (21) \circ  (\on{Ad}(u) \circ \on{R}_-
\circ \sigma \circ \psi^{-1})^{\boxtimes 2}
\\ &
= \psi\circ \sigma^{-1} \circ \on{R}_-^{(-1)}
\circ m_\Pi \circ (21) \circ  (\on{R}_-
\circ \sigma \circ \psi^{-1})^{\boxtimes 2}
\\ &
= (\psi\circ \sigma^{-1})  \circ
(\on{R}_-^{(-1)}   \circ m_\Pi \circ \on{R}_-^{\boxtimes 2}) \circ (21)
\circ (\sigma \circ \psi^{-1})^{\boxtimes 2}
\end{align*}
and
\begin{align*}
& \Delta_a^\tau = \big( (\on{R}_+^{(-1)})^{\boxtimes 2} \circ
\Delta_\Pi \circ \on{R}_+\big)^\tau
= ((\on{R}_+^{(-1)})^\tau)^{\boxtimes 2} \circ
\Delta_\Pi^\tau \circ \on{R}_+^\tau
\\ & =
(\psi\circ\sigma^{-1} \circ \on{R}_-^{(-1)} )^{\boxtimes 2}
\circ (21) \circ \Delta_\Pi
\circ \on{R}_- \circ \sigma\circ \psi^{-1}
\\ & =
(\psi\circ \sigma^{-1})^{\boxtimes 2} \circ
(21) \circ
(\on{R}_-^{(-1)} )^{\boxtimes 2}
\circ \Delta_\Pi
\circ \on{R}_- \circ (\sigma\circ \psi^{-1}).
\end{align*}

Now
\begin{align*}
& m_a^t = (m_a^\tau)^*  =
((\sigma \circ \psi^{-1})^*)^{\boxtimes 2} \circ (21) \circ
\big((\on{R}_-^*)^{\boxtimes 2} \circ m_\Pi^* \circ \on{R}_-^{(-1)*} \big)
\circ (\psi\circ \sigma^{-1})^*
\\ &
= ((\sigma \circ \psi^{-1})^*)^{\boxtimes 2} \circ (21) \circ \Delta_a
\circ (\psi\circ \sigma^{-1})^*,
\end{align*}
and
\begin{align*}
& \Delta_a^t = (\Delta_a^\tau)^* =
(\sigma\circ \psi^{-1})^* \circ \on{R}_-^* \circ \Delta_\Pi^* \circ
(\on{R}_-^{(-1)*})^{\boxtimes 2} \circ (21) \circ
((\psi\circ \sigma^{-1})^*)^{\boxtimes 2}
\\ &
= (\sigma\circ \psi^{-1})^* \circ m_a  \circ
((\psi\circ \sigma^{-1})^*)^{\boxtimes 2}.
\end{align*}

Let $S_a\in\hbox{\boldmath$\Pi$\unboldmath}
(S\underline\boxtimes {\mathfrak 1},
S\underline\boxtimes {\mathfrak 1})$ be the antipode for
$(m_a,\Delta_a)$. Then $m_a \circ (21) = S_a \circ m_a \circ
(S_a^{-1})^{\boxtimes 2}$,
$(21) \circ \Delta_a = S_a^{\boxtimes 2} \circ \Delta_a
\circ S_a^{-1}$. Therefore
$$
m_a^t = ((\sigma\circ\psi^{-1})^* \circ S_a)^{\boxtimes 2} \circ \Delta_a
\circ ((\sigma\circ\psi^{-1})^* \circ S_a)^{-1},
$$
$$
\Delta_a^t = ((\sigma\circ\psi^{-1})^* \circ S_a) 
\circ m_a \circ (21)
\circ (((\sigma\circ\psi^{-1})^* \circ S_a)^{-1})^{\boxtimes 2},
$$
as wanted.
\hfill \qed \medskip

\subsection{From almost compatible to compatible (with duality)
quantization functors}

Let $Q$ be any quantization functor. Recall that is gives rise to a prop
isomorphism $Q^* : \ul{\on{Bialg}}^* \to S({\bf LBA})^*$. We say that
$Q$ is almost compatible with duality iff there exists
$\xi\in S({\bf LBA})({\bf id},{\bf id})^\times$, of the form 
$\xi = \on{id}_{\bf id}^{S({\bf LBA})}$ + higher order terms, such that
$Q^* \circ *cop_{\on{Bialg}} = \on{Inn}(\xi^*) \circ S(*cop_{\on{LBA}})\circ Q$
(recall that $\xi\mapsto \xi^*$ is the tautological map
$P(F,G) \to P^*(G^*,F^*)$), and $Q$ is compatible with duality
if this identity holds with $\xi = \on{id}_{{\bf id}}^{S(\on{LBA})}$,
i.e., if the diagram of props
$$
\xymatrix{
\on{Bialg} \ar[r]^{Q}\ar[d]_{*cop_{\on{Bialg}}}
& S(\on{LBA})\ar[d]^{S(*cop_{\on{LBA}})} \\
\on{Bialg}^* \ar[r]^{Q^*} &
S({\bf LBA})^* & \simeq S({\bf LBA}^*)
}
$$
commutes.

\begin{prop}
If $Q$ is almost compatible with duality, then there exists
$\xi_0\in S({\bf LBA})({\bf id},{\bf id})^\times$ such that
$\on{Inn}(\xi_0) \circ Q$ is compatible with duality.
\end{prop}

{\em Proof.} Let $\xi$ be such that $Q^* \circ *cop_{\on{Bialg}} =
\on{Inn}(\xi^*) \circ *cop_{S({\bf LBA})} \circ Q$. Since $\on{Inn}(\xi^*)^*
= \on{Inn}(\xi^{-1})$, we get $Q \circ (*cop_{\on{Bialg}})^* = 
\on{Inn}(\xi^{-1}) \circ (*cop_{S({\bf LBA})})^* \circ Q^*$. 
Since $(*cop_{S({\bf LBA})})^* \circ *cop_{S({\bf LBA})} 
= op\ cop_{S({\bf LBA})}$
and $(*cop_{\on{Bialg}})^* \circ *cop_{\on{Bialg}} = op\ cop_{\on{Bialg}}$, 
we get $Q \circ op\ cop_{\on{Bialg}} = \on{Inn}(\xi^{-1}) \circ 
(*cop_{S({\bf LBA})})^* \circ \on{Inn}(\xi^*) \circ *cop_{S({\bf LBA})} 
\circ Q = \on{Inn}(\xi^{-1} \circ (*cop_{S({\bf LBA})})^*(\xi^*)) 
\circ op\ cop_{S({\bf LBA})} \circ Q$, and since $(*cop_{S({\bf LBA})})^*(\xi^*)
= (*cop_{S({\bf LBA})}(\xi))^*$, we get 
$$
Q \circ op\ cop_{\on{Bialg}} =\on{Inn}(\xi_1) \circ 
op\ cop_{S({\bf LBA})} \circ Q,  
$$
where $\xi_1 := \xi^{-1} \circ(*cop_{S({\bf LBA})}(\xi))^*$. Since 
$(op\ cop)^2 = 1$, we get 
$$
Q = \on{Inn}(\xi_1 \circ op\ cop_{S({\bf LBA})}(\xi_1)) \circ Q, 
$$ 
so $\xi_1 \circ op\ cop_{S({\bf LBA})}(\xi_1)=\on{exp}(\alpha(\mu\circ\delta))$
for some scalar $\alpha$. Since $op\ cop_{{\bf LBA}}(\mu\circ\delta) =
\mu\circ\delta$, if we set 
$\xi_2:= \on{exp}(-\alpha(\mu\circ\delta)/2)
\circ \xi_1$, we have $\xi_2 \circ op\ cop_{S({\bf LBA})}(\xi_2)
=\on{id}_{{\bf id}}^{S({\bf LBA})}$, and 
$$
Q \circ op\ cop_{\on{Bialg}} =\on{Inn}(\xi_2) \circ 
op\ cop_{S({\bf LBA})} \circ Q.   
$$ 
Since $\xi_2 = \on{id}_{\bf id}^{S({\bf LBA})}$ + higher order terms, 
there exists a unique $\xi_2^{1/2}$ of the same form, with 
$\xi_2 = (\xi_2^{1/2})^2$. Since $\xi_2 = 
op\ cop_{S({\bf LBA})}(\xi_2^{-1})$, we have 
$\xi_2 = 
(op\ cop_{S({\bf LBA})}((\xi_2^{1/2})^{-1}))^2$ so 
$\xi_2^{1/2} = 
op\ cop_{S({\bf LBA})}(\xi_2^{1/2})^{-1}$, hence 
$$
\xi_2 = \xi_2^{1/2} \circ op\ cop_{S({\bf LBA})}(\xi_2^{1/2})^{-1}, 
$$
so if we set $Q_1 := \on{Inn}(\xi_2^{1/2})^{-1}\circ Q$, we get 
$$
Q_1 \circ op\ cop_{\on{Bialg}}= op\ cop_{S({\bf LBA})} \circ Q_1.  
$$
We also have 
$$
Q_1^* \circ *cop_{\on{Bialg}} =
\on{Inn}(\xi_3^*) \circ *cop_{S({\bf LBA})} \circ Q_1,  
$$
where $\xi_3 = (*cop_{S({\bf LBA})}(\xi_2^{1/2}))^* \circ \xi 
\circ \xi_2^{1/2}$.  Then 
$(*cop_{S({\bf LBA})}(\xi_3))^* = (*cop_{S({\bf LBA})}(\xi_2^{1/2}))^* \circ 
*cop_{S({\bf LBA})}(\xi)^* \circ 
*cop_{{\bf LBA}}(*cop_{S({\bf LBA})}(\xi_2^{1/2})^*)^*$. 
Now $\xi \mapsto *cop_{S({\bf LBA})}(\xi)^*$ is a prop 
antiautomorphism of $S({\bf LBA})$, so its square is a 
prop automorphism of $S({\bf LBA})$; one check on generators 
this it is equal to $op\ cop_{S({\bf LBA})}$. Therefore 
\begin{align*}
& (*cop_{S({\bf LBA})}(\xi_3))^* 
= (*cop_{S({\bf LBA})}(\xi_2^{1/2}))^* \circ 
(*cop_{S({\bf LBA})}(\xi)^* \circ op\ cop_{S({\bf LBA})}(\xi_2^{1/2})
\\
 & = (*cop_{S({\bf LBA})}(\xi_2^{1/2}))^* \circ 
*cop_{S({\bf LBA})}(\xi)^* \circ \xi_2^{-1/2} 
= 
(*cop_{S({\bf LBA})}(\xi_2^{1/2}))^* \circ 
\xi \circ \xi_1 \circ \xi_2^{-1/2}  \\
& =
(*cop_{S({\bf LBA})}(\xi_2^{1/2}))^* \circ 
\xi \circ \xi_2^{1/2} \circ \on{exp}(\alpha(\mu\circ\delta)/2) = 
\xi_3 \circ \on{exp}(\alpha(\mu\circ\delta)/2).
\end{align*} 

Set $\xi_4 := \xi_3 \circ \on{exp}(\alpha(\mu\circ\delta)/4)$, then 
$(*cop_{S({\bf LBA})}(\xi_4))^* = \xi_4$, as $(*cop_{{\bf LBA}}
(\mu\circ\delta))^* = -\mu\circ\delta$, and we have 
$$
Q_1^* \circ *cop_{\on{Bialg}} =
\on{Inn}(\xi_4^*) \circ *cop_{S({\bf LBA})} \circ Q_1. 
$$ 
Let $\xi_{5}:= \xi_{4}^{1/2}$. Then $\xi_{5}=*cop_{S({\bf LBA})}(\xi_{5})^{*}$,
so $\xi_{4}^{*}=\xi_{5}^{*}\circ *cop_{S({\bf LBA})}(\xi_{5})$, so 
$\on{Inn}(\xi_{5}^{*})^{-1}\circ Q_{1}^{*}\circ *cop_{\on{Bialg}}
=\on{Inn}(*cop_{S({\bf LBA})}(\xi_{5})) \circ *cop_{S({\bf LBA})}
\circ Q_{1} = *cop_{S({\bf LBA})}\circ \on{Inn}(\xi_{5})\circ Q_{1}$.
So if we set $Q_{2}:= \on{Inn}(\xi_{5})\circ Q_{1}$, we have 
$Q_{2}^{*}\circ *cop_{\on{Bialg}}=*cop_{S({\bf LBA})}\circ Q_{2}$, 
while $Q_{2}=\on{Inn}(\xi_{0})\circ Q$, and $\xi_{0}=\xi_{5}
\circ \xi_{2}^{-1/2}$.  \hfill \qed \medskip

\section{Compatibility of quantization functors with doubling operations}

In this section, we formulate propic versions of the (Drinfeld) double constructions of
Lie bialgebras and Hopf algebras. We then express a condition on a quantization
functor $Q : \ul{\on{Bialg}} \to S({\bf LBA})$, which we call ``compatibility
with doubling operations''. One then proves that that if $Q$ is compatible with duality,
then it is compatible with doubling operations.

\subsection{Doubles of Lie bialgebras} \label{sec:double:LBA}

In this subsection, we set $P:= \on{LBA}$ or $\ul{\on{LBA}}$, ${\bf LBA}$.

We define a prop $D_{add}(P)$ by $D_{add}(P)(F,G) :=
P_2(\ul\Delta(F),\ul\Delta(G))$.
We also define prop $(D_{add}(P),P)$-bimodules
$M^\pm_{add}(P)$
by $M_{add}^+(P)(F,G) := P_2(F \ul\boxtimes {\bf 1},\ul\Delta(G))$
and $M_{add}^-(P)(F,G) := P_2({\bf 1}\ul\boxtimes F,\ul\Delta(G))$.

The structures of left prop $D_{add}(P)$-bimodules on
$M_{add}^\pm(P)$ are obvious. Let us define the structures of
right prop $P$-modules on $M_{add}^\pm(P)$.

In the case of $M_{add}^+(P)$, the composition
$P(F,G) \otimes P_2(G \ul\boxtimes {\bf 1},
\ul\Delta(H)) \to P_2(F\ul\boxtimes {\bf 1},\ul\Delta(H))$ is $x \otimes X
\mapsto X \circ (x\ul\boxtimes {\rm id}_{{\bf id}}^{\on{LBA}})$.

In the case of $M_{add}^-(P)$, the analogous composition is the map
$P(F,G) \otimes P_2({\bf 1} \ul\boxtimes G,\ul\Delta(H)) \to
P_2({\bf 1}\ul\boxtimes F,\ul\Delta(H))$ such that $x \otimes X
\mapsto X \circ ({\rm id}_{{\bf 1}}^{P} \ul\boxtimes *cop(x))$, where
$*cop = cop \circ * : P(F,G)\to P^*(F,G)$.

$M_{add}^-(P)$ also has a $(D_{add}(P),P^*)$-prop bimodule
structure, induced by
$$
P^*(F,G) \otimes P_2({\bf 1} \ul\boxtimes G,
\ul\Delta(H)) \to P_2({\bf 1} \ul\boxtimes F,\ul\Delta(H)), \quad
x^* \otimes X \mapsto X \circ ({\rm id}_{{\bf 1}}^{P} \ul\boxtimes x^*).
$$

Note that if $V$ is a finite dimensional module over $P$, i.e., we have a prop morphism
$P \to \on{Prop}(V)$, then $V\oplus V^*$ is a module over $D_{add}(P)$.

We then have
\begin{align*}
& D_{add}(\on{LBA})(\wedge^2,{\bf id})
\\ &
= \on{LBA}_2(\wedge^2\ul\boxtimes {\bf 1},{\bf id} \ul\boxtimes {\bf 1})
\oplus \on{LBA}_2({\bf id} \ul\boxtimes {\bf id}, {\bf id} \ul\boxtimes {\bf 1})
\oplus \on{LBA}_2({\bf id} \ul\boxtimes {\bf id}, {\bf 1} \ul\boxtimes {\bf id})
\oplus \on{LBA}_2({\bf 1} \ul\boxtimes \wedge^2, {\bf 1} \ul\boxtimes {\bf id})
\\ & \simeq \on{LBA}(\wedge^2,{\bf id}) \oplus \on{LBA}({\bf id},T_2) \oplus
\on{LBA}(T_2,{\bf id}) \oplus \on{LBA}({\bf id},\wedge^2).
\end{align*}
Similarly,
\begin{align*}
& D_{add}(\on{LBA})({\bf id},\wedge^2)
\\ & =
\on{LBA}_2({\bf id} \ul\boxtimes {\bf 1},\wedge^2\ul\boxtimes {\bf 1})
\oplus \on{LBA}_2({\bf id} \ul\boxtimes {\bf 1},{\bf id} \ul\boxtimes {\bf id})
\oplus \on{LBA}_2({\bf 1} \ul\boxtimes {\bf id},{\bf id} \ul\boxtimes {\bf id})
\oplus \on{LBA}_2({\bf 1} \ul\boxtimes {\bf id},{\bf 1} \ul\boxtimes \wedge^2)
\\ & \simeq \on{LBA}({\bf id},\wedge^2) \oplus \on{LBA}(\wedge^2,{\bf id})
\end{align*}
as the two intermediate spaces are zero.

We have $M_{add}^+(\on{LBA})({\bf id},{\bf id})
= \on{LBA}_2({\bf id}\ul\boxtimes {\bf 1},
{\bf id}\ul\boxtimes {\bf 1}) \oplus \on{LBA}_2({\bf id} \ul\boxtimes {\bf 1},
{\bf 1} \ul\boxtimes {\bf id})\simeq \on{LBA}({\bf id},{\bf id})$ as the second
space is zero
and $M_{add}^-(\on{LBA})({\bf id},{\bf id})
= \on{LBA}_2({\bf 1}\ul\boxtimes{\bf id},
{\bf id}\ul\boxtimes{\bf 1}) \oplus\on{LBA}_2({\bf 1}\ul\boxtimes {\bf id},
{\bf 1}\ul\boxtimes {\bf id}) \simeq \on{LBA}({\bf id},{\bf id})$ as the first
space is zero.

\begin{lemma} \label{lemma:double}
There are unique prop morphisms $\on{double} : \on{LBA} \to
D_{add}(\on{LBA})$, left prop $D_{add}(\on{LBA})$-module
morphisms $\alpha^\pm_{can} : D_{add}(\on{LBA}) \to M_{add}^\pm(\on{LBA})$
and right prop $\on{LBA}$-module morphisms $\beta^\pm_{can} : \on{LBA}
\to M_{add}^\pm(\on{LBA})$, such that
$$
\on{double}(\mu) \simeq \mu \oplus \delta \oplus (-\mu) \oplus \delta,
\quad \on{double}(\delta) \simeq \delta \oplus (-\mu),
$$
$$
\alpha^+_{can}(\on{id}_{{\bf id}}^{D_{add}(\on{LBA})}) =
\beta^+_{can}(\on{id}_{{\bf id}}^{\on{LBA}}) \simeq
\on{id}_{{\bf id}}^{\on{LBA}},
\quad \alpha^-_{can}(\on{id}_{{\bf id}}^{D_{add}(\on{LBA})}) =
\beta^-_{can}(\on{id}_{{\bf id}}^{\on{LBA}})
\simeq \on{id}_{{\bf id}}^{\on{LBA}}.
$$
The diagrams
$$
\xymatrix 
{& M^\pm_{add}(\on{LBA}) & \\
\on{LBA} \ar[ru]^{\beta_{can}^\pm} \ar[rr]_{\on{double}} & &
D_{add}(\on{LBA})\ar[ul]_{\alpha_{can}^\pm}
}
$$
commute.
\end{lemma}

We define $\on{can}_F^+\in \on{LBA}_2(F\ul\boxtimes {\bf 1},\ul\Delta(F))$
and $\on{can}_F^-\in \on{LBA}_2({\bf 1}\ul\boxtimes F,\ul\Delta(F))$
as follows: $\on{can}_F^\pm$ are the images of the elements of
$\on{Sch}_{1+1}(F\ul\boxtimes {\bf 1},\ul\Delta(F))$ and
$\on{Sch}_{1+1}({\bf 1}\ul\boxtimes F,\ul\Delta(F))$, which are the
universal versions of $F(V)\to F(V\oplus V^*)$ induced by $v\mapsto
v\oplus 0$, resp.  $F(V^*)\to F(V\oplus V^*)$ induced by $v^*\mapsto
0\oplus v^*$.

Then $\alpha_{can}^{\pm},\beta_{can}^{\pm}$ are given by
$\alpha_{can}^\pm(X) = X \circ \on{can}_F^\pm$,
$\beta_{can}^+(x) = \on{can}_G^+ \circ (x\ul\boxtimes
\on{id}_{{\bf 1}}^{\on{LBA}^*})$,
$\beta_{can}^-(x) = \on{can}_G^- \circ (\on{id}_{{\bf 1}}^{\on{LBA}}
\boxtimes *cop(x))$.

The proof is a propic version of (a) the double Lie bialgebra
construction $\A\mapsto {\mathfrak D}(\A)$ and (b) the Lie
bialgebra morphisms $\A\to {\mathfrak D}(\A)$, $\A^{*,cop} \to
{\mathfrak D}(\A)$, where $\A$ is a finite dimensional Lie
bialgebra.

These morphisms have analogues when $\on{LBA}$ is replaced by
$\ul{\on{LBA}}$ or ${\bf LBA}$.

\subsection{Doubles of Hopf algebras}

In this section, we will set $P := \underline{\on{Bialg}}$ or $S({\bf LBA})$.

Define a prop $D_{mult}(P)$
by $D_{mult}(P)(F,G) := P_2(\ul\Delta_{\otimes}(F),\ul\Delta_{\otimes}(G))$,
where $\ul\Delta_\otimes : \on{Sch} \to \on{Sch}_{1+1}$ is given by
$\ul\Delta_\otimes(F)(V,W) := F(V\otimes W)$. We also define prop
$(D_{mult}(P),P)$-bimodules $M^\pm_{mult}(P)$
by
$$M_{mult}^+(P)(F,G) := P_2(F \ul\boxtimes {\bf 1},\ul\Delta_{\otimes}(G))
\quad\on{and}\quad
M_{mult}^-(P)(F,G) := P_2({\bf 1}\ul\boxtimes F,\ul\Delta_\otimes(G)).
$$

The bimodule structures are defined as follows: the left prop module
structure over $D_{mult}(P)$ is obvious; the right prop module structures
of $M_{mult}^\pm(P)$ over $P$ are defined as above, replacing $\ul\Delta$ by
$\ul\Delta_{\otimes}$ and in the case of $M_{mult}^-(P)$, $*cop$ by
$*cop_{\on{Bialg}}$. As above, $M_{mult}^{-}(P)$ is also a
$(D_{mult}(P),P^{*})$-prop bimodule.

We have $D_{mult}(\ul{\on{Bialg}})(T_2,{\bf id})
=\ul{\on{Bialg}}_2(T_2\ul\boxtimes T_2,
{\bf id}\ul\boxtimes{\bf id}) \subset \ul{\on{Bialg}}(T_2\otimes {\bf id},{\bf id}
\otimes T_2)$ and similarly, $D_{mult}(\ul{\on{Bialg}})({\bf id},T_2) \subset
\ul{\on{Bialg}}({\bf id}\otimes T_2,T_2 \otimes {\bf id})$.

\begin{lemma} There are unique prop morphisms $\on{Double} :
\on{Bialg} \to D_{mult}(\underline{\on{Bialg}})$, left
prop $D_{mult}(\ul{\on{Bialg}})$-module morphisms $\alpha_{\on{Bialg}}^\pm :
D_{mult}(\ul{\on{Bialg}}) \to M_{add}^\pm(\ul{\on{Bialg}})$ and
right prop $\on{Bialg}$-module morphisms $\beta_{\on{Bialg}}^\pm :
\on{Bialg} \to M_{mult}^\pm(\ul{\on{Bialg}})$, such that
$$
\on{Double}(m) \simeq (m\boxtimes m^{(2)}\boxtimes
\on{id}_{\bf id}^{\ul{\on{Bialg}}})
\circ (152346) \circ a^{[2]} \circ (\on{id}_{\bf id}^{\ul{\on{Bialg}}}\boxtimes
\Delta^{(2)}\boxtimes\Delta),
$$
$$\on{Double}(\Delta) \simeq \Delta \boxtimes
((21) \circ m), \quad \on{Double}(\varepsilon)
= \varepsilon \ul\boxtimes \eta^{*}, \quad \on{Double}(\eta)  =
\eta\ul\boxtimes \varepsilon^{*},
$$
$$
\alpha_{\on{Bialg}}^+(\on{id}_{{\bf id}}^{D_{mult}(\ul{\on{Bialg}})})  =
\beta_{\on{Bialg}}^+(\on{id}_{{\bf id}}^{{\on{Bialg}}}) =
\on{id}_{\bf id}^{\ul{\on{Bialg}}} \ul\boxtimes \eta^{*},
$$
$$
\alpha_{\on{Bialg}}^-(\on{id}_{{\bf id}}^{D_{mult}(\ul{\on{Bialg}})})
= \beta_{\on{Bialg}}^-(\on{id}_{{\bf id}}^{{\on{Bialg}}})
= \varepsilon \ul\boxtimes \on{id}_{\bf id}^{\ul{\on{Bialg}}^{*}}.
$$
The diagrams
$$
\xymatrix 
{& M^\pm_{mult}(\ul{\on{Bialg}}) & \\
\on{Bialg} \ar[ru]^{\beta_{\on{Bialg}}^\pm} \ar[rr]_{\on{Double}} & &
D_{mult}(\ul{\on{Bialg}})\ar[ul]_{\alpha_{\on{Bialg}}^\pm}
}
$$
commute.
\end{lemma}

Here $a^{[2]} = \on{id}_{\bf id}^{\ul{\on{Bialg}}} \boxtimes
a\boxtimes (\on{id}_{\bf id}^{\ul{\on{Bialg}}})^{\boxtimes 4}$.
$\on{Double}(m)$ is the propic version of the map $A^{\otimes 3}\to
A^{\otimes 3}$, $x\otimes y \otimes z\mapsto xy^{(2)} \otimes
y^{(3)}z^{(1)}a(y^{(1)})\otimes z^{(2)}$; it is obtained by
dualizing the formula of the multiplication formula of the Drinfeld double.

The prop bimodule properties are expressed as follows:
$$
\alpha_{\on{Bialg}}^+(F,H)(Y\circ X) = Y \circ \alpha_{\on{Bialg}}^+(F,G)(X)
$$
for $X\in D_{mult}(\ul{\on{Bialg}})(F,G)$, $Y\in D_{mult}(\ul{\on{Bialg}})(G,H)$,
$$
\beta_{\on{Bialg}}^+(F,H)(y\circ x) = \beta_{\ul{\on{Bialg}}}^+(G,H)(y)
\circ (x\ul\boxtimes \on{id}_{\bf 1}^{\ul{\on{Bialg}}^*}),
$$
$$
\beta_{\on{Bialg}}^-(F,H)(y\circ x) = \beta_{\on{Bialg}}^-(G,H)(y)
\circ (\on{id}_{\bf 1}^{\ul{\on{Bialg}}} \ul\boxtimes *cop_{\on{Bialg}}(x)).
$$
for $x\in D_{mult}(\ul{\on{Bialg}})(F,G)$, $y\in D_{mult}(\ul{\on{Bialg}})(G,H)$.

The proof is the universal version of the proof of the fact that the Drinfeld
double $D(A)$ of a Hopf algebra $A$ is a bialgebra, equipped with bialgebra
morphisms $A\to D(A)$ and $A^{*,cop}\to D(A)$.

Since $\on{Double}(\on{id}_{\bf id}^{\on{Bialg}} -
\eta\circ\varepsilon) = (\on{id}_{\bf id}^{\ul{\on{Bialg}}} -
\eta\circ\varepsilon) \ul\boxtimes \on{id}_{\bf
id}^{\ul{\on{Bialg}}^*} + (\eta\circ\varepsilon) \ul\boxtimes
(\on{id}_{\bf id}^{\ul{\on{Bialg}}} - \eta\circ\varepsilon)^*$,
$\on{Double}$ extends to a prop morphism $\on{Double} :
\ul{\on{Bialg}} \to D_{mult}(\underline{\on{Bialg}})$ with the
same properties.

\subsection{Compatibility of quantization functors with doubling operations}

Let $Q : \ul{\on{Bialg}} \to S({\bf LBA})$ be a quantization functor. As we have seen,
this is a prop isomorphism.

\begin{prop}
$Q$ gives rise to an isomorphism of biprops
$Q_2 : \ul{\on{Bialg}}_2 \to S^{\ul\boxtimes 2}({\bf LBA}_2)$.
\end{prop}

{\em Proof.} Let $F,...,G'\in\on{Sch}$. The morphism $Q$ gives rise to a
continuous linear isomorphism
\begin{align*}
& Q(F\otimes G^{\prime *},F^{\prime *}\otimes G) :
\\ & \ul{\on{Bialg}}(F\otimes G^{\prime *},F^{\prime *}\otimes G) \to
S({\bf LBA})(F\otimes G^{\prime *},F'\otimes G)
= {\bf LBA}(S(F)\otimes S(G^{\prime *}),S(F')\otimes S(G))
\end{align*}
(denoted shortly by $Q$).
Recall that $\ul{\on{Bialg}}_2(F\ul\boxtimes G,F'\ul\boxtimes G')
\subset \ul{\on{Bialg}}(F\otimes G^{\prime *},F'\otimes G)$
and
$$
S^{\ul\boxtimes 2}({\bf LBA}_2)(F\ul\boxtimes G,F'\ul\boxtimes G')
= {\bf LBA}_2(S(F)\ul\boxtimes S(G),S(F')\ul\boxtimes S(G'))
\subset {\bf LBA}(S(F)\otimes S(G^{\prime *}),S(F')\otimes S(G)).
$$
We will prove that $Q$
maps $\ul{\on{Bialg}}_2(F\ul\boxtimes G,F'\ul\boxtimes G')$ bijectively to
${\bf LBA}_2(S(F)\ul\boxtimes S(G),S(F')\ul\boxtimes S(G'))$.

We first prove that $Q$ maps the first space into the second.
Consider the diagram
$$
\xymatrix{
*\txt{ $\hat\oplus_{Z_{XY}\in\on{Irr(Sch)}}$ \\
$\on{Coalg}(F,Z_{FF'}\otimes Z_{FG})$ \\
$\otimes \on{Coalg}(G^{\prime*},Z_{G'G})$
\\ $\otimes \on{Alg}(Z_{FF'},F')$ \\
$\otimes \on{Alg}(Z_{FG}\otimes Z_{G'G},G^*)$}
\ar[r]^>>>>>>>>>>>>>>>>>>>>{(a)}  \ar[d]_>>>>>>{Q'}
&
\ul{\on{Bialg}}(F\otimes G^{\prime*},F'\otimes G^*) \ar[d]^Q
\\
*\txt{
$\hat\oplus_{Z_{XY}\in\on{Irr(Sch)}}$ \\
${\bf LBA}(S(F),S(Z_{FF'})\otimes S(Z_{FG}))$ \\
$\otimes {\bf LBA}(S(G^{\prime*}),S(Z_{G'G}))$ \\
$\otimes {\bf LBA}(S(Z_{FF'}),S(F'))$ \\
$\otimes {\bf LBA}(S(Z_{FG}) \otimes S(Z_{G'G}),S(G^*)) $}
\ar[r]^>>>>>>>>>>>>>{(b)} \ar[dd]^>>>>>>>>>>>>>>{(c)}
&
{\bf LBA}(S(F)\otimes S(G^{\prime*}),S(F') \otimes S(G^*))
\\
& *\txt{$\hat\oplus_{W_{XY}} {\bf LCA}(S(F),W_{FF'} \otimes W_{FG})
\otimes {\bf LCA}(S(G^{\prime*}),W_{G'G})$ \\ $ \otimes{\bf LA}(W_{FF'},S(F'))
\otimes {\bf LA}(W_{FG}\otimes W_{G'G},S(G^*))$} \ar[u]^>>>>>>>>>>{(f)}
\\
*\txt{
$\hat\oplus_{Z_{XY},Z_{X|XY},Z_{XY|Y}\in\on{Irr(Sch)}}$ \\
${\bf LCA}(S(F),Z_{F|FF'}\otimes Z_{F|FG})$ \\
$\otimes {\bf LA}(Z_{F|FF'},S(Z_{FF'}))$ \\
$\otimes {\bf LA}(Z_{F|FG},S(Z_{FG}))$ \\
$\otimes {\bf LCA}(S(G^{\prime*}),Z_{G'|G'G})$ \\
$\otimes {\bf LA}(Z_{G'|G'G},S(Z_{G'G}))$ \\
$\otimes {\bf LCA}(S(Z_{FF'}),Z_{FF'|F'})$ \\
$\otimes {\bf LA}(Z_{FF'|F'},S(F'))$ \\
$\otimes {\bf LCA}(S(Z_{FG}),Z_{FG|G})$ \\
$\otimes {\bf LCA}(S(Z_{G'G}),Z_{G'G|G})$ \\
$\otimes {\bf LA}(Z_{FG|G} \otimes Z_{G'G|G},S(G^*))$
}\ar[r]^>>>>>>>>>>>{(d)}
& *\txt{
$\hat\oplus_{W_{XY},Z_{X|XY},Z_{XY|Y}\in\on{Irr(Sch)}}$ \\
${\bf LCA}(S(F),Z_{F|FF'}\otimes Z_{F|FG})$ \\
$\otimes {\bf LCA}(S(G^{\prime*}),Z_{G'|G'G})$ \\
$\otimes {\bf LA}(Z_{FF'|F'},S(F'))$ \\
$\otimes {\bf LA}(Z_{FG|G} \otimes Z_{G'G|G},S(G^*))$ \\
$\otimes {\bf LCA}(Z_{F|FF'},W_{FF'})
\otimes {\bf LA}(W_{FF'},Z_{FF'|F'})$ \\
$\otimes {\bf LCA}(Z_{F|FG},W_{FG})
\otimes {\bf LA}(W_{FG},Z_{FG|G})$ \\
$\otimes {\bf LCA}(Z_{G'|G'G},W_{G'G})
\otimes {\bf LA}(W_{G'G},Z_{G'G|G})$
}\ar[u]^>>>>>>>>{(e)}
}
$$
where the maps labeled $Q$, $Q'$ are those induced by $Q$,
the map $(a)$ is $c_{F|F'G} \otimes c_{G'|G} \otimes a_{F|F'}
\otimes a_{FG'|G} \mapsto
(a_{F|F'} \boxtimes a_{FG'|G}) \circ
\big( (\on{id} - \eta\circ\varepsilon)_{Z_{FF'}} \boxtimes
(\on{id} - \eta\circ\varepsilon)_{Z_{FG}} \boxtimes
(\on{id} - \eta\circ\varepsilon)_{Z_{G'G}} \big)
\circ (c_{F|F'G} \boxtimes c_{G'|G'G})$
(recall that $inj_{i}\in {\bf Sch}(S^{i},S)$ and $pr_{i}
\in {\bf Sch}(S,S^{i})$ are the canonical injection and projection); the map $(b)$ is
$\lambda_{F|FG} \otimes \lambda_{G'|G} \otimes \lambda_{F|F'}
\otimes \lambda_{FG'|G} \mapsto
(\lambda_{F|F'} \boxtimes \lambda_{FG'|G}) \circ
[(\on{id}_{S} - inj_0 \circ pr_0)(Z_{FF'})^{\on{LBA}}
\boxtimes
(\on{id}_{S} - inj_0 \circ pr_0)(Z_{FG})^{\on{LBA}}
\boxtimes
(\on{id}_{S} - inj_0 \circ pr_0)(Z_{G'G})^{\on{LBA}}]
\circ (\lambda_{F|FG} \boxtimes \lambda_{G'|G'G})$; the map $(c)$
is induced by the isomorphisms (\ref{isoms:LBA}); the map $(d)$
is induced by the composition
$$
\hat\oplus_W {\bf LA}(Z,S(W)) \otimes {\bf LCA}(S(W),Z')
\to {\bf LBA}(Z,Z') \to \hat\oplus_{W'} \on{LCA}(Z,W') \otimes
\on{LA}(W',Z'),
$$
where the first map is $\sum_W \lambda_W \otimes \kappa_W
\mapsto \sum_W \kappa_W \circ (\on{id}_{S} - inj_0 \circ
pr_0)(W)^{\on{LBA}} \circ \lambda_W$ and the second map
is the inverse of $\sum_{W'} \kappa_{W'} \circ \lambda_{W'}
\mapsto \sum_{W'} \kappa_{W'} \circ \lambda_{W'}$; the map $(e)$
is induced by the compositions
$$
{\bf LCA}(S(F),Z_{F|FF'} \otimes Z_{F|FG}) \otimes
{\bf LCA}(Z_{F|FF'},W_{FF'})\otimes
{\bf LCA}(Z_{F|FG},W_{FG}) \to {\bf LCA}(S(F),W_{FF'} \otimes W_{FG}),
$$
$\kappa \otimes \kappa' \otimes \kappa'' \mapsto
(\kappa'\boxtimes \kappa'') \circ \kappa$, and the analogous
compositions with ${\bf LA}$ instead of ${\bf LCA}$. Finally, $(f)$
is the standard composition map.

The upper square of this diagram commutes because $Q$ is
a prop morphism. One also checks that the lower square of
this diagram commutes, so $Q \circ (a) = (f) \circ (e) \circ (d) \circ (c)
\circ Q'$. The image of $(a)$ is $\ul{\on{Bialg}}_2(F\ul\boxtimes G,
F'\ul\boxtimes G')$, so its image by $Q$ is the image of $Q \circ (a)$,
which is therefore contained in the image of $(f)$, which is
${\bf LBA}_2(S(F)\ul\boxtimes S(G),S(F')\ul\boxtimes S(G'))$.

Let us now check that $Q^{-1}$ maps
${\bf LBA}_2(S(F)\ul\boxtimes S(G),S(F')\ul\boxtimes S(G'))$
to $\ul{\on{Bialg}}_2(F\ul\boxtimes G,F'\ul\boxtimes G')$.
Consider the diagram
$$
\xymatrix{
*\txt{ $\hat\oplus_{Z_{XY}\in\on{Irr(Sch)}}$ \\
${\bf LCA}(S(F),Z_{FF'}\otimes Z_{FG})$ \\
$\otimes {\bf LCA}(S(G^{\prime*}),Z_{G'G})$
\\ $\otimes {\bf LA}(Z_{FF'},S(F'))$ \\
$\otimes {\bf LA}(Z_{FG}\otimes Z_{G'G},S(G^*))$}
\ar[r]^>>>>>>>>>{(a)}  \ar[d]_>>>>>>{\tilde Q^{-1}}
&
{\bf LBA}(S(F)\otimes S(G^{\prime*}),S(F')\otimes S(G^*)) \ar[d]^{Q^{-1}}
\\
*\txt{
$\hat\oplus_{Z_{XY}\in\on{Irr(Sch)}}$ \\
$\ul{\on{Bialg}}(F,Z_{FF'}\otimes Z_{FG})$ \\
$\otimes \ul{\on{Bialg}}(G^{\prime*},Z_{G'G})$ \\
$\otimes \ul{\on{Bialg}}(Z_{FF'},F')$ \\
$\otimes \ul{\on{Bialg}}(Z_{FG} \otimes Z_{G'G},G^*) $}
\ar[r]^>>>>>>>>>>>>>>{(b)} \ar[dd]^>>>>>>>>>>>>>{(c)}
&
\ul{\on{Bialg}}(F\otimes G^{\prime*},F' \otimes G^*)
\\
& *\txt{$\hat\oplus_{W_{XY}} \on{Coalg}(F,W_{FF'} \otimes W_{FG})
\otimes \on{Coalg}(G^{\prime*},W_{G'G})$ \\ $ \otimes \on{Alg}(W_{FF'},F')
\otimes \on{Alg}(W_{FG}\otimes W_{G'G},G^*)$} \ar[u]^>>>>>>>>{(f)}
\\
*\txt{
$\hat\oplus_{Z_{XY},Z_{X|XY},Z_{XY|Y}\in\on{Irr(Sch)}}$ \\
$\on{Coalg}(F,Z_{F|FF'}\otimes Z_{F|FG})$ \\
$\otimes \on{Alg}(Z_{F|FF'},Z_{FF'})$ \\
$\otimes \on{Alg}(Z_{F|FG},Z_{FG})$ \\
$\otimes \on{Coalg}(G^{\prime*},Z_{G'|G'G})$ \\
$\otimes \on{Alg}(Z_{G'|G'G},Z_{G'G})$ \\
$\otimes \on{Coalg}(Z_{FF'},Z_{FF'|F'})$ \\
$\otimes \on{Alg}(Z_{FF'|F'},F')$ \\
$\otimes \on{Coalg}(Z_{FG},Z_{FG|G})$ \\
$\otimes \on{Coalg}(Z_{G'G},Z_{G'G|G})$ \\
$\otimes \on{Alg}(Z_{FG|G} \otimes Z_{G'G|G},G^*)$
}\ar[r]^>>>>>>>>>>{(d)}
& *\txt{
$\hat\oplus_{W_{XY},Z_{X|XY},Z_{XY|Y}\in\on{Irr(Sch)}}$ \\
$\on{Coalg}(F,Z_{F|FF'}\otimes Z_{F|FG})$ \\
$\otimes \on{Coalg}(S(G^{\prime*}),Z_{G'|G'G})$ \\
$\otimes \on{Alg}(Z_{FF'|F'},F')$ \\
$\otimes \on{Alg}(Z_{FG|G} \otimes Z_{G'G|G},G^*)$ \\
$\otimes \on{Coalg}(Z_{F|FF'},W_{FF'})$ \\
$\otimes \on{Alg}(W_{FF'},Z_{FF'|F'})$ \\
$\otimes \on{Coalg}(Z_{F|FG},W_{FG})$ \\
$\otimes \on{Alg}(W_{FG},Z_{FG|G})$ \\
$\otimes \on{Coalg}(Z_{G'|G'G},W_{G'G})$ \\
$\otimes \on{Alg}(W_{G'G},Z_{G'G|G})$
}\ar[u]^>>>>>{(e)}
}
$$
Here $\tilde Q^{-1}$ is the composed map
\begin{align*}
& \scriptstyle{ \hat\oplus_{Z_{XY}\in\on{Irr(Sch)}}
{\bf LCA}(S(F),Z_{FF'}\otimes Z_{FG})
\otimes {\bf LCA}(S(G^{\prime*}),Z_{G'G})
\otimes {\bf LA}(Z_{FF'},S(F'))
\otimes {\bf LA}(Z_{FG}\otimes Z_{G'G},S(G^*))}
\\ & \scriptstyle{\to
\hat\oplus_{Z_{XY}\in\on{Irr(Sch)}}
{\bf LCA}(S(F),S(Z_{FF'})\otimes S(Z_{FG}))
\otimes {\bf LCA}(S(G^{\prime*}),S(Z_{G'G}))
\otimes {\bf LA}(S(Z_{FF'}),S(F'))
\otimes {\bf LA}(S(Z_{FG})\otimes S(Z_{G'G}),S(G^*))}
\\ & \scriptstyle{\to
\hat\oplus_{Z_{XY}\in\on{Irr(Sch)}}
\ul{\on{Bialg}}(F,Z_{FF'}\otimes Z_{FG})
\otimes \ul{\on{Bialg}}(G^{\prime*},Z_{G'G})
\otimes \ul{\on{Bialg}}(Z_{FF'},F')
\otimes \ul{\on{Bialg}}(Z_{FG}\otimes S(Z_{G'G}),S(G^*)),}
\end{align*}
in which the first map is the tensor product of maps
${\bf LCA}(S(A),B) \to {\bf LCA}(S(A),S(B))$, $c\mapsto
inj_1(A)^{\on{LCA}} \circ c$ and ${\bf LA}(A,S(B)) \to
{\bf LA}(S(A),S(B))$,
$\ell \mapsto \ell \circ pr_1(A)^{\on{LA}}$, and the second map
is the tensor product of the maps induced by $Q^{-1}$.
Here the map $(a)$ is defined by
$\kappa_{F|F'G} \otimes \kappa_{G'|G} \otimes \alpha_{F|F'}
\otimes \alpha_{FG'|G} \mapsto
(\alpha_{F|F'} \boxtimes \alpha_{FG'|G}) \circ
(\kappa_{F|F'G} \boxtimes \kappa_{G'|G'G})$;
the map $(b)$ is $b_{F|F'G} \otimes b_{G'|G'G} \otimes b_{FF'|F'}
\otimes b_{FG'|G}
\mapsto (b_{FF'|F'} \boxtimes b_{FG'|G})
\circ [(\on{id} - \eta\circ \varepsilon)_{Z_{FF'}}
\boxtimes (\on{id} - \eta\circ \varepsilon)_{Z_{FG}}
\boxtimes (\on{id} - \eta\circ \varepsilon)_{Z_{G'G}}]
\circ (b_{F|F'G}\boxtimes b_{G'|G'G})$; the map $(c)$
is induced by the isomorphisms (\ref{isom:bialg}) ; the map $(d)$ is a tensor
product of composed maps
$$
\on{Alg}(A,Z) \otimes \on{Coalg}(Z,B) \to \ul{\on{Bialg}}(A,B)
\to \hat\oplus_{W\in\on{Irr(Sch)}} \on{Coalg}(A,W) \otimes
\on{Alg}(W,B),
$$
where the first map is $a\otimes c \mapsto c\circ
(\on{id}-\eta\circ\varepsilon)_Z
\circ a$ and the second map is inverse to $\sum_W c_W\otimes a_W \mapsto
\sum_W a_W \circ (\on{id} - \eta\circ \varepsilon) \circ c_W$; the map $(e)$
is a tensor product of the maps $\on{Coalg}(A,Z) \otimes \on{Coalg}(Z,W)
\to \on{Coalg}(A,W)$, $c\otimes c'\mapsto c'\circ (\on{id} - \eta\circ
\varepsilon)_Z \circ c$, their analogues when $(Z,W)$ is replaced by
two pairs $(Z_i,W_i)$ ($i=1,2$) or when $\on{Coalg}$ is replaced by
$\on{Alg}$; the map $(f)$ is the map $\sum_{W_{XY}}
c_{F|F'G} \otimes c_{G'|G} \otimes a_{F|F'} \otimes a_{FG'|G} \mapsto
(a_{F|F'} \boxtimes a_{FG'|G}) \circ [
(\on{id} - \eta\circ \varepsilon)_{W_{FF'}}\boxtimes
(\on{id} - \eta\circ \varepsilon)_{W_{FG}}\boxtimes
(\on{id} - \eta\circ \varepsilon)_{W_{G'G}}] \circ
(c_{F|F'G}\boxtimes c_{G'|G})$.

The map $(e)$ is well-defined, since nonzero elements of the
the space labeled $(W_{XY}) = (W_{FF'},W_{FG},W_{G'G})$ can be in the image
only of the spaces labeled $(W_{XY},Z_{X|XY},Z_{XY}|Y)$,
where\footnote{Recall that $|Z|$ is the degree of a homogeneous Schur functor.}
$|Z_{X|XY}|,
|Z_{XY|Y}| \leq |W_{XY}|$, for each pair $(X,Y) \in \{(F,F'),(F,G),(G',G)\}$.

The upper square of this diagram commutes since it can be split into
commuting triangle, sum over the $Z_{XY}$ of
$$
\xymatrix{
*\txt{ ${\bf LCA}(S(F),Z_{FF'}\otimes Z_{FG})$ \\
$\otimes {\bf LCA}(S(G^{\prime*}),Z_{G'G})$
\\ $\otimes {\bf LA}(Z_{FF'},S(F'))$ \\
$\otimes {\bf LA}(Z_{FG}\otimes Z_{G'G},S(G^*))$}
\ar[r]^>>>>>>>>>{(a)}  \ar[d]_>>>>>>{(g)}
& {\bf LBA}(S(F)\otimes S(G^{\prime*}),S(F')\otimes S(G^*))
\\
*\txt{ ${\bf LCA}(S(F),S(Z_{FF'})\otimes S(Z_{FG}))$ \\
$\otimes {\bf LCA}(S(G^{\prime*}),S(Z_{G'G}))$
\\ $\otimes {\bf LA}(S(Z_{FF'}),S(F'))$ \\
$\otimes {\bf LA}(S(Z_{FG})\otimes S(Z_{G'G}),S(G^*))$}
\ar[ur]_{(h)}
}
$$
where $(g)$ is the tensor product of a map ${\bf LCA}(S(A),Z)
\to {\bf LCA}(S(A),S(Z))$, $\kappa \mapsto inj_1(Z)^{\on{LCA}}
\circ \kappa$ with its analogues with $Z$ replaced by $Z,Z'$,
and with ${\bf LCA}$ replaced by ${\bf LA}$ (and $inj_1$ by $pr_1$);
and $(h)$ is the map $\kappa_{F|F'G}\otimes \kappa_{G'|G}
\otimes \lambda_{F|F'} \otimes \lambda_{FG'|G}
\mapsto (\lambda_{F|F'} \otimes \lambda_{FG'|G}) \circ
[(\on{id}_{S} - inj_0 \circ pr_0)(Z_{FF'})^{\on{LBA}}\boxtimes
(\on{id}_{S} - inj_0 \circ pr_0)(Z_{FG})^{\on{LBA}}\boxtimes
(\on{id}_{S} - inj_0 \circ pr_0)(Z_{G'G})^{\on{LBA}}]
\circ (\kappa_{F|F'G} \boxtimes \kappa_{G'|G})$,
which commutes since $(\on{id}_{S} - inj_0 \circ pr_0) \circ inj_1
= inj_1 \circ (\on{id}_{S} - inj_0 \circ pr_0) =0$; and of the commuting square
$$
\xymatrix{
*\txt{ $\hat\oplus_{Z_{XY}\in \on{Irr(Sch)}}$ \\
${\bf LCA}(S(F),S(Z_{FF'})\otimes S(Z_{FG}))$ \\
$\otimes {\bf LCA}(S(G^{\prime*}),S(Z_{G'G}))$
\\ $\otimes {\bf LA}(S(Z_{FF'}),S(F'))$ \\
$\otimes {\bf LA}(S(Z_{FG})\otimes S(Z_{G'G}),S(G^*))$}
\ar[r]^>>>>>>{(i)}\ar[d]_{Q^{-1}}&
{\bf LBA}(S(F)\otimes S(G^{\prime*}),S(F')\otimes S(G^*))
\ar[d]^{Q^{-1}}
\\
*\txt{$\hat\oplus_{Z_{XY}\in \on{Irr(Sch)}}$ \\
$\on{Coalg}(F,Z_{FF'}\otimes Z_{FG})$ \\
$\otimes \on{Coalg}(G^{\prime*},Z_{G'G})$
\\ $\otimes \on{Alg}(Z_{FF'},F')$ \\
$\otimes \on{Alg}(Z_{FG}\otimes Z_{G'G},G^*)$}\ar[r]_{(j)}
& \ul{\on{Bialg}}(F\otimes G^{\prime*},F'\otimes G^*)
}
$$
where the map $(i)$ is given by $\kappa^{S}_{F|F'G} \otimes \kappa^{S}_{G'|G}
\otimes \lambda^{S}_{F|F'} \otimes \lambda^{S}_{FG'|G'} \to
(\lambda^{S}_{F|F'}\boxtimes \lambda^{S}_{FG'|G}) \circ
[(\on{id}_{S}-inj_{0}\circ pr_{0})(Z_{FF'})^{\on{LBA}}\boxtimes
(\on{id}_{S}-inj_{0}\circ pr_{0})(Z_{FG})^{\on{LBA}}\boxtimes
(\on{id}_{S}-inj_{0}\circ pr_{0})(Z_{G'G})^{\on{LBA}}]\circ
(\kappa^{S}_{F|F'G}\boxtimes \kappa^{S}_{G'|G})$ and the map $(j)$
is given by the same expression, where
$(\on{id}- inj_{0}\circ pr_{0})^{\on{LBA}}_{S(Z)}$
is replaced by $(\on{id}-\eta\circ\varepsilon)_{Z}^{\on{Bialg}}$;
this square commutes, since $Q((\on{id}-\eta\circ\varepsilon)_{Z}) =
(\on{id}_{S} - inj_{0}\circ pr_{0})(Z)^{\on{LBA}}$.  One also checks that
the bottom square of above the diagram  commutes. Using the same argument as
before, one obtains that $Q^{-1}$ maps
${\bf LBA}_2(S(F)\ul\boxtimes S(G),S(F')\ul\boxtimes S(G'))$
to $\ul{\on{Bialg}}_2(F\ul\boxtimes G,F'\ul\boxtimes G')$.

In the same way as one proved above that $Q : \ul{\on{Bialg}}(F\otimes G^{\prime*},
F'\otimes G^{*}) \to S({\bf LBA})(F\otimes G^{\prime *},F'\otimes G^{*})$
restricts to a map $\ul{\on{Bialg}}_{2}(F\ul{\boxtimes} G,F'\ul{\boxtimes}G')
\to S({\bf LBA})_{2}(F\ul{\boxtimes} G,F'\ul{\boxtimes}G')$, one proves
more generally that for any finite sets $I,J$, any $\Sigma\subset I\times J$,
and any $A\in \on{Sch}_{I}$, $B\in \on{Sch}_{J}$,
$Q : \ul{\on{Bialg}}(\otimes(A),\otimes(B)) \to S({\bf LBA})(\otimes(A),
\otimes(B))$ restricts to a linear map
(in fact an isomorphism) $\ul{\on{Bialg}}^{\Sigma}(A,B)\to {\bf LBA}^{\Sigma}
(S(A),S(B))$ (where $S : \on{Sch}_{I} \to {\bf Sch}_{I}$ takes the functor
$A : \on{Vect}^{I}\to \on{Vect}$ to $\on{Vect}^{I} \stackrel{S^{I}}{\to}
{\bf Vect}^{I} \stackrel{\tilde A}{\to} {\bf Vect}$, where
$\tilde A$ is the natural extension of $A$). If moreover $\bullet\notin I\cup J$
and $\bar\Sigma\subset \bar I\times\bar J$ is such that $(\bullet,\bullet)\notin
\bar\Sigma$ (where $\bar I = I \sqcup \{\bullet\}$, $\bar J = J \sqcup
\{\bullet\}$), then the trace maps $\on{tr} : \ul{\on{Bialg}}^{\bar\Sigma}
(A\boxtimes F,B \boxtimes F) \to \ul{\on{Bialg}}^{\Sigma}(A,B)$
and $\on{tr} : {\bf LBA}^{\bar\Sigma}(S(A)\boxtimes S(F),S(B)\boxtimes S(F))
\to {\bf LBA}^{\Sigma}(S(A),S(B))$ are such that the diagram
$$\xymatrix{
 \ul{\on{Bialg}}^{\bar\Sigma}(A\boxtimes F,B \boxtimes F)
 \ar[r]^>>>>>>{Q}\ar[d]_{\on{tr}}&
{\bf LBA}^{\bar\Sigma}(S(A)\boxtimes S(F),S(B)\boxtimes S(F))
\ar[d]_{\on{tr}}
\\
 \ul{\on{Bialg}}^{\Sigma}(A,B)\ar[r]^{Q}  &
{\bf LBA}^{\Sigma}(S(A),S(B))
}$$
commutes.  Since the trace is the basic ingredient of the prop structures of
$\ul{\on{Bialg}}_{2}$ and $S({\bf LBA})_{2}$, it follows that $Q_{2}$
is a biprop morphism, hence isomorphism. \hfill \qed \medskip

\begin{lemma}
$Q$ gives rise to a prop isomorphism
$D_{mult}(Q) :  D_{mult}(\ul{\on{Bialg}}) \to D_{mult}(S({\bf LBA}))$.
\end{lemma}

{\em Proof.}
For $F,G\in\on{Sch}$, we have isomorphisms
$D_{mult}(\ul{\on{Bialg}})(F,G) = \ul{\on{Bialg}}_2(\ul\Delta_{\otimes}(F),
\ul\Delta_{\otimes}(G)) \stackrel{Q_2}{\to}
S^{\ul\boxtimes 2}({\bf LBA}_2)(\ul\Delta_{\otimes}(F),
\ul\Delta_{\otimes}(G)) \simeq S({\bf LBA})_2(\ul\Delta_{\otimes}(F),
\ul\Delta_{\otimes}(G)) = D_{mult}(S({\bf LBA}))(F,G)$. One
checks that they are compatible with the prop operations.
\hfill \qed \medskip

\begin{lemma}
We have a canonical prop isomorphism
$D_{mult}(S({\bf LBA})) \simeq S(D_{add}({\bf LBA}))$.
\end{lemma}

{\em Proof.} This is given by
\begin{align*}
& D_{mult}(S({\bf LBA}))(F,G) = S({\bf
LBA})_2(\ul\Delta_{\otimes}(F), \ul\Delta_{\otimes}(G)) = {\bf
LBA}_2(S^{\ul\boxtimes 2} \circ \ul\Delta_{\otimes}(F),
S^{\ul\boxtimes 2} \circ \ul\Delta_{\otimes}(G))
\\ & =
{\bf LBA}_2(\ul\Delta(S(F)),\ul\Delta(S(G))) =
S(D_{add}({\bf LBA}))(F,G) = D_{add}({\bf LBA})(S(F),S(G));
\end{align*}
the third equality uses the canonical isomorphism
$S^{\ul\boxtimes 2} \circ \ul\Delta_{\otimes}(F) \simeq \ul\Delta(S(F))$,
which follows from the isomorphism $F(S(V) \otimes S(V)^*) \simeq
F(S(V\oplus V^*))$ for any $V\in \on{Vect}$.   \hfill \qed \medskip

We say that $Q$ is compatible with the doubling operations
iff there exists an inner
automorphism $\on{Inn}(\Lambda)$ of $S(D_{add}({\bf LBA}))$, where
$\Lambda\in S(D_{add}({\bf LBA}))({\bf id},{\bf id})^\times$, such
that the diagram
$$
\xymatrix@!0 @R=8ex @C=25ex
{
\ul{\on{Bialg}} \ar[rrr]^{\on{Double}}\ar[d]_{Q} & & &
D_{mult}(\ul{\on{Bialg}})\ar[d]^{D_{mult}(Q)} \\
S({\bf LBA}) \ar[r]_{S(\on{double})} & S(D_{add}({\bf LBA}))
\ar@{..>}[r]_{\on{Inn}(\lambda)} & S(D_{add}({\bf LBA}))
\ar[r]_{\simeq} & D_{mult}(S({\bf LBA}))
}
$$
commutes.

Recall that $Q$ gives rise to a prop morphism $Q^* : \ul{\on{Bialg}}^*
\to S({\bf LBA})^*$, and $Q$ is called compatible with duality iff
the diagram of props
$$
\xymatrix{
\ul{\on{Bialg}} \ar[r]^{Q}\ar[d]_{*cop_{\on{Bialg}}}
& S({\bf LBA})\ar[d]^{S(*cop_{\on{LBA}})} \\
\ul{\on{Bialg}}^* \ar[r]^{Q^*} &
S({\bf LBA})^* & \simeq  S({\bf LBA}^*)
}
$$
commutes.

\begin{thm} \label{thm:doubles}
If a quantization functor $Q$ is compatible with duality, then
it is also compatible with the doubling operations. In particular, the
EK quantization functors are compatible with the doubling operations.
\end{thm}

The rest of this section is devoted to the proof of this theorem.

\subsection{Construction of commuting triangles of prop modules based on
$S({\bf LBA})$}

We have shown that the prop isomorphism $Q : \ul{\on{Bialg}} \to S({\bf LBA})$
induces a prop isomorphism $D_{mult}(Q) : D_{mult}(\ul{\on{Bialg}})
\to D_{mult}(S({\bf LBA}))$.

\begin{lemma}
$Q$ induces isomorphisms
$$
M^\pm_{mult}(Q)(F,G) : M^\pm_{mult}(\ul{\on{Bialg}})(F,G) \to
M^\pm_{mult}(S({\bf LBA}))(F,G),
$$
(shortly denoted $M^\pm_{mult}(Q)$) with the following
compatibilities with $Q$, $Q^*$ and $D_{mult}(Q)$:
$$
M_{mult}^\pm(Q)(X\circ m) = D_{mult}(Q)(X) \circ M_{mult}^\pm(Q)(m)
$$
for $m\in M_{mult}^\pm(\ul{\on{Bialg}})(F,G)$, $X\in D_{mult}(\ul{\on{Bialg}})(G,H)$,
$$
M_{mult}^+(Q)(m \circ (x\ul\boxtimes\on{id}_{\bf 1}^{\on{Bialg}^*}))
= M_{mult}^+(Q)(m) \circ (Q(x)\ul\boxtimes\on{id}_{\bf 1}^{S(\on{LBA})^*}),
$$
$$
M_{mult}^-(Q)(m \circ (\on{id}_{\bf 1}^{\on{Bialg}}\ul\boxtimes x^*))
= M_{mult}^-(Q)(m) \circ (\on{id}_{\bf 1}^{S(\on{LBA})}\ul\boxtimes
Q^*(x^*)),
$$
for $x^*\in \ul{\on{Bialg}}^*(F,G)$, $m\in M_{mult}^\pm(\ul{\on{Bialg}})(G,H)$.

So $M^\pm_{mult}(Q) : M^\pm_{mult}(\ul{\on{Bialg}}) \to
M^\pm_{mult}(S({\bf LBA}))$ are prop bimodule isomorphisms
compatible with the prop morphisms and with $Q$, $Q^*$ and $D_{mult}(Q)$.
\end{lemma}

{\em Proof.}
We have isomorphisms
\begin{align*}
& M^+_{mult}(\ul{\on{Bialg}})(F,G)
= \ul{\on{Bialg}}_2(F\ul\boxtimes {\bf 1},\ul\Delta_{\otimes}(G))
\\ & \stackrel{Q_2}{\to} S^{\ul\boxtimes 2}({\bf LBA}_2)(F\ul\boxtimes {\bf 1},
\ul\Delta_{\otimes}(G))= S({\bf LBA})_2(F\ul\boxtimes {\bf 1},
\ul\Delta_{\otimes}(G)) = M^+_{mult}(S({\bf LBA}))(F,G),
\end{align*}
and
$M^-_{mult}(\ul{\on{Bialg}})(F,G) \simeq M^-_{mult}(S({\bf LBA}))(F,G)$
by replacing $F\ul\boxtimes {\bf 1}$ by ${\bf 1}\ul\boxtimes F$.

The properties of these isomorphisms follow from the fact that $Q_2$
is a biprop morphism such that $Q_2(x\ul\boxtimes
\on{id}_{\bf 1}^{\ul{\on{Bialg}}^*})
= Q(x) \ul\boxtimes \on{id}_{\bf 1}^{S({\bf LBA})^*}$, and
$Q_2(\on{id}_{\bf 1}^{\ul{\on{Bialg}}} \ul\boxtimes x^*)
= \on{id}_{\bf 1}^{S({\bf LBA})}\ul\boxtimes Q^*(x^*)$.
\hfill \qed\medskip

We now transport $\on{Double}$, $\alpha_{\on{Bialg}}^\pm$ and $\beta_{\on{Bialg}}^\pm$
using the isomorphisms $Q$, $D_{mult}(Q)$ and $M_{mult}^\pm(Q)$ as follows.

For $F,G\in\on{Sch}$, define $\tilde\alpha^\pm(F,G) : D_{mult}(S({\bf LBA}))(F,G)
\to M_{mult}^\pm(S({\bf LBA}))$ by
$$
\tilde\alpha^\pm(F,G) := M_{mult}^\pm(Q)(F,G) \circ
\alpha^\pm_{\on{Bialg}}(F,G) \circ D_{mult}(Q)(F,G)^{-1},
$$
define
$\tilde\beta^\pm(F,G) : S({\bf LBA})(F,G)\to
M_{mult}^\pm(S({\bf LBA}))$ by
$$
\tilde\beta^\pm(F,G) := M_{mult}^\pm(Q)(F,G) \circ
\beta^\pm_{\on{Bialg}}(F,G) \circ Q(F,G)^{-1},
$$
and define $\widetilde{\on{Double}}(F,G) : S({\bf LBA})(F,G)\to
D_{mult}(S({\bf LBA}))$ by
$$
\widetilde{\on{Double}}(F,G) := D_{mult}(Q)(F,G) \circ
\on{Double}(F,G) \circ Q(F,G)^{-1}.
$$

\begin{prop} \label{prop:module}
$\widetilde{\on{Double}}$ is a prop morphism. If $Q$ is compatible
with duality, then we have
$$
\tilde\alpha^\pm(F,H)(Y \circ X) = Y \circ \tilde \alpha^\pm(F,G)(X),
$$
for $X\in D_{mult}(S({\bf LBA}))(F,G)$, $Y\in D_{mult}(S({\bf LBA}))(G,H)$,
$$
\tilde\beta^+(F,H)(y \circ x) = \tilde\beta^+(G,H)(y) \circ
(x\ul\boxtimes\on{id}_{\bf 1}^{S({\bf LBA})^*}),
$$
$$
\tilde\beta^-(F,H)(y \circ x) = \tilde\beta^-(G,H)(y) \circ
(\on{id}_{\bf 1}^{S({\bf LBA})}\ul\boxtimes *cop_{S({\bf LBA})}(x)).
$$
for $x\in S({\bf LBA})(F,G)$, $y\in S({\bf LBA})(G,H)$, and
$\tilde\alpha^\pm(F,G)\circ \widetilde{\on{Double}}(F,G)
= \tilde \beta^\pm(F,G)$ for any $F,G$.
\end{prop}

{\em Proof.} The first statement follows from the fact that
$\widetilde{\on{Double}}$ is a composition of prop morphisms.
The first (resp., second) identity follows from the prop isomorphism
property of $D_{mult}(Q)$ (resp., $Q$) and the prop module properties
of $\alpha_{\on{Bialg}}^\pm$ (resp., $\beta_{\on{Bialg}}^+$)
and $M_{mult}^\pm(Q)$ (resp., $M_{mult}^+(Q)$).
Using the prop isomorphism property of $Q$ and the
prop module properties of $\beta_{\on{Bialg}}^-$ and $M^{-}_{mult}(Q)$,
we prove that $\tilde\beta^-(F,H)(y \circ x) = \tilde\beta^-(G,H)(y) \circ
(\on{id}_{\bf 1}^{S({\bf LBA})}\ul\boxtimes (Q^* \circ
*cop_{\on{Bialg}} \circ Q^{-1})(x))$, which implies the third identity
since $Q$ is compatible with duality.

The last identity follows from
$\alpha^\pm_{\on{Bialg}}(F,G)\circ \on{Double}(F,G)
= \beta^\pm_{\on{Bialg}}(F,G)$ for any $F,G$.
\hfill \qed \medskip

These identities mean that we have a prop morphism
$\widetilde{\on{Double}} : S({\bf LBA}) \to D_{mult}(S({\bf LBA}))$,
left prop $D_{mult}(S({\bf LBA}))$-module morphisms
$\tilde\alpha^\pm : D_{mult}(S({\bf LBA})) \to M_{mult}^\pm(S({\bf LBA}))$
and right prop $S({\bf LBA})$-module morphisms
$\tilde\beta^\pm : S({\bf LBA}) \to M_{mult}^\pm(S({\bf LBA}))$, such that
$\tilde\alpha^\pm\circ \widetilde{\on{Double}}
= \tilde \beta^\pm$.

Now using the canonical isomorphisms
$X_{mult}(S({\bf LBA})) \simeq S(X_{add}({\bf LBA}))$ for
$X \in \{D,M^\pm\} $, we view these morphisms as:
a prop morphism
$\widetilde{\on{Double}} : S({\bf LBA}) \to S(D_{add}({\bf LBA}))$,
left prop $S(D_{add}({\bf LBA}))$-module morphisms
$\tilde\alpha^\pm : S(D_{add}({\bf LBA})) \to S(M_{add}^\pm({\bf LBA}))$
and right prop $S({\bf LBA})$-module morphisms
$\tilde\beta^\pm : S({\bf LBA}) \to S(M_{add}^\pm({\bf LBA}))$, such that
the diagrams
$$
\xymatrix 
{& S(M^\pm_{add}({\bf LBA})) & \\
S({\bf LBA}) \ar[ru]^{\tilde\beta^\pm} \ar[rr]_{\widetilde{\on{Double}}} & &
S(D_{add}({\bf LBA}))\ar[ul]_{\tilde\alpha^\pm}
}
$$
commute. The prop module morphism properties are
expressed by the relations of Proposition \ref{prop:module},
where now $X\in S(D_{add}({\bf LBA}))(F,G)$,
$Y\in S(D_{add}({\bf LBA}))(G,H)$.

\subsection{Construction of a prop morphism $\varphi : {\bf LBA} \to
D_{add}({\bf LBA})$}

The diagram of prop modules
$$ \xymatrix{ & S({\bf Sch})
\ar[dl]_{S(i_{{\bf LBA}})}
\ar[dr]^{S(i_{D_{add}({\bf LBA})})}
\\ S({\bf LBA})\ar[rr]_{\widetilde{\on{Double}}  } & &
S(D_{add}({\bf LBA})) }
$$
does not necessarily commute.

In this section, we will construct an inner automorphism $\on{Inn}(\Xi)$
of $S(D_{add}({\bf LBA}))$ (where $\Xi\in 
S({\bf LBA})({\bf id},{\bf id})^\times$), and a prop morphism 
$\varphi : {\bf LBA} \to D_{add}({\bf LBA})$, such that 
$\on{Inn}(\Xi)^{-1} \circ \widetilde{\on{Double}} \circ S(i_{{\bf LBA}})
=S(i_{D_{add}({\bf LBA})})$ and 
$\on{Inn}(\Xi)^{-1} \circ \widetilde{\on{Double}} = S(\varphi)$.
The existence of $\Xi$ and $\varphi$ follows from Proposition 
\ref{gen:prop} below. 

We first define generators $m^{S},\Delta^{S},p_{n}^{S}$ ($n\geq 0$)
of the prop $S({\bf Sch})$. Recall that $m^S\in
S({\bf Sch})(T_2,{\bf id})$ and $\Delta^S\in
S({\bf Sch})({\bf id},T_2)$ are the universal versions of 
the product $S(V)^{\otimes 2}\to S(V)$ and the coproduct 
$S(V) \to S(V)^{\otimes 2}$ (where $V$ is a vector space). We define 
$p^S_n\in S({\bf Sch})({\bf id},{\bf id})={\bf Sch}(S,S)$ as 
$p_{n}^{S}:= inj_{n}\circ pr_{n}$, i.e., $p_{n}^{S}$ is the universal version of
the projector $p_n^V \in \on{End}(S(V))$ onto $S^n(V) \subset
\oplus_{n\geq 0} S^n(V) = S(V)$. One can prove that $S({\bf Sch})$ is generated by
the elements $m^{S},\Delta^{S},p_{n}^{S}$.

\begin{prop} \label{gen:prop}
Let ${\bf P},{\bf Q}$ be ${\bf Sch}$-props graded by $\NN$,
complete and separated for these gradings. Assume that $\Phi :
S({\bf P}) \to S({\bf Q})$ is a prop morphism (not necessarily
compatible with the gradings), such that for any of the generators
$x\in \{p_n^S,m^S,\Delta^S\}$ of $S({\bf Sch})$, $(S(i_{\bf Q}) -
\Phi\circ S(i_{\bf P}))(x)$ has positive degree ($i_{\bf P} : {\bf
Sch} \to {\bf P}$, $i_{\bf Q} : {\bf Sch} \to {\bf Q}$ are the
canonical morphisms).

Then there exists a prop morphism $\varphi : {\bf P} \to {\bf Q}$ and an
invertible element $\Xi\in S({\bf Q})({\bf id},{\bf id})^\times$, such that
$\Phi = \on{Inn}(\Xi) \circ S(\varphi)$.
\end{prop}

{\em Proof.} We first prove that if the smallest degree of all
$(S(i_{\bf Q}) - \Phi\circ S(i_{\bf P}))(x)$ is $N>0$ (for $x\in S({\bf Sch})(F,G)$ 
and $F,G\in \on{Ob(Sch)}$), then one can construct
$\Xi_N\in S({\bf Q})({\bf id},{\bf id})$ of degree $N$, such that
the degree $N$ part of $(S(i_{\bf Q}) - \on{Inn}({\rm id}_{{\bf id}} + \Xi_N)
\circ \Phi\circ S(i_{\bf P}))(x)$ vanishes for any $x$.
Defining $\Xi$ as the product of
all the ${\rm id}_{{\bf id}} + \Xi_N$, we will then have:
$S(i_{\bf Q}) = \on{Inn}(\Xi)^{-1} \circ \Phi \circ S(i_{\bf P})$.

For $x^S\in S({\bf Sch})(F,G)$ and
$F,G\in \on{Ob(Sch)}$, define $\dot x^Q\in S({\bf Q})(F,G)$ as the degree $N$
part of $(S(i_{\bf Q}) - \Phi \circ S(i_{\bf P}))(x^S)$ and set
$x^Q := S(i_{\bf Q})(x^S)$.
The condition on $\Xi_N$ is that for any such $x^S$,
$$
\dot x^Q = (d/dt)_{|t=0}(\on{Inn}(\on{id}_{{\bf id}} + t \Xi_N)(x^Q));
$$
it suffices that this holds when $x^S$ belongs to a set of generators of
$S({\bf Sch})$.

We set $\Xi_N = \sum_{n,m\geq 0} (\Xi_N)_{n,m}$ according to the decomposition
$S({\bf Q})({\bf id},{\bf id}) = \hat\oplus_{n,m\geq 0} {\bf Q}(S^n,S^m)$;
we will first construct the off-diagonal part of $\Xi_N$, then its
diagonal part.

The idempotent relations between $p_n^S$ imply
$p_n^Q \circ \dot p_m^Q + \dot p_n^Q \circ p_m^Q = \delta_{n,m}
\dot p_n^Q$ and $\sum_n \dot p_n^Q = 0$.

For $n\neq m$, we then set $(\Xi_N)_{n,m} := p_m^Q \circ \dot p_n^Q\circ p_n^Q$
and $(\Xi_N)_{off}:= \sum_{n,m|n\neq m} (\Xi_N)_{n,m}$.
We will prove that for any $i\geq 0$, $\dot p_i^Q = [(\Xi_N)_{off},p_i^Q]$
(where $[a,b] = a\circ b - b\circ a$).

We have $(\Xi_N)_{off} = \sum_{n} (\dot p_n^Q \circ p_n^Q -
p_n^Q \circ \dot p_n^Q \circ p_n^Q)$, therefore
$[(\Xi_N)_{off}, p_i^Q] =
\sum_n [\dot p_n^Q \circ p_n^Q,p_i^Q] - \sum_n [p_n^Q \circ \dot p_n^Q \circ
p_n^Q,p_i^Q]
= \dot p_i^Q \circ p_i^Q - \sum_n p_i^Q \circ \dot p_n^Q \circ p_n^Q$ as the
second bracket vanishes. Then $- \sum_n p_i^Q \circ \dot p_n^Q \circ p_n^Q
= \sum_n \dot p_i^Q \circ p_n^Q - \dot p_i^Q \circ p_i^Q$, so that
$[(\Xi_N)_{off}, p_i^Q] = \dot p_i^Q$. It follows that the degree
$N$ part of $(S(i_{\bf Q}) - \on{Inn}({\rm id}_{{\bf id}} + (\Xi_N)_{off})
\circ \Phi\circ S(i_{\bf P}))(p_n^S)$ vanishes.

For $x^S\in S({\bf Sch})(F,G)$, define $\breve{x}^Q$ as the degree $N$ part of
$(S(i_{\bf Q}) - \on{Inn}({\rm id}_{{\bf id}} + (\Xi_N)_{off})
\circ \Phi\circ S(i_{\bf P}))(x^S)$. We have $\breve p_n^Q=0$.
We now construct $(\Xi_N)_{n,n} \in {\bf Q}(S^n,S^n)$, such that if
$(\Xi_N)_{diag} := \sum_{n\geq 0} (\Xi_N)_{n,n}$, then
$\breve x^Q = d/dt_{|t=0} \on{Inn}(\on{id}_{{\bf id}} + t (\Xi_N)_{diag})
(x^Q)$.

Let for $n\geq 0$, $m^S_n :=
(m^S)^{(n)} \circ (p_1^S)^{\boxtimes n}
\in S({\bf Sch})(T_n,{\bf id})$ and
$\Delta^S_n :=
(p_1^S)^{\boxtimes n} \circ (\Delta^S)^{(n)}
\in S({\bf Sch})({\bf id},T_n)$, where
$(m^S)^{(n)} = m^S \circ ... \circ (m^S \boxtimes
\on{id}_{{\bf id}}^{\boxtimes n-2})$
and $(\Delta^S)^{(n)} = (\Delta^S \boxtimes
\on{id}_{{\bf id}}^{\boxtimes n-2}) \circ ... \circ \Delta^S$.
These elements generate $S({\bf Sch})$, as $m^S$ and $\Delta^S$ can be
constructed from them by
$m^S = \sum_{p,q\geq 0} {1\over {p!q!}} m_{p+q}^S \circ (\Delta_p^S \boxtimes
\Delta_q^S)$ and $\Delta^S = \sum_{p,q\geq 0} {1\over {p!q!}} (m_p^S \boxtimes
m_q^S) \circ \Delta_{p+q}^S$.

They satisfy the relations
\begin{equation} \label{eq:1}
m_1^S = \Delta_1^S = p_1^S,
\end{equation}
\begin{equation} \label{eq:2}
p^S_{n'} \circ m^S_n = \delta_{n,n'}
m^S_n, \quad
m^S_n \circ (\boxtimes_{i=1}^n p^S_{\nu_i}) =
\delta_{(\nu_1,...,\nu_n),(1,...,1)} m^S_n,
\end{equation}
\begin{equation} \label{eq:3}
\Delta^S_n \circ p^S_{n'} = \delta_{n,n'}
\Delta^S_n, \quad
(\boxtimes_{i=1}^n p^S_{\nu_i}) \circ \Delta^S_n =
\delta_{(\nu_1,...,\nu_n),(1,...,1)} \Delta^S_n,
\end{equation}
\begin{equation} \label{rel:4}
\Delta^S_n \circ m^S_n
= n! \on{Sym}_n \circ (p_1^S)^{\boxtimes n}, \quad
m^S_n \circ \Delta^S_n = n! p_n^S,
\end{equation}
where $\on{Sym}_n\in S({\bf Sch})(T_n,T_n)$ is the total symmetrization.
(\ref{eq:2}) implies
\begin{equation} \label{eq:sym}
m^S_n \circ \on{Sym}_n = m^S_n, \quad
\on{Sym}_n \circ \Delta^S_n =\Delta^S_n.
\end{equation}

Set
$$
(\Xi_N)_{nn} := {1\over {n!}}\breve m^Q_n \circ \Delta^Q_n.
$$
(\ref{eq:1}) implies that $(\Xi_N)_{11}=0$.

We now prove that
$$
\breve m^Q_n = (\Xi_N)_{n,n} \circ
m^Q_n - m^Q_n \circ
(\sum_{i=1}^n (\Xi_N)_{11}^{(i)}).
$$
(\ref{eq:1}) implies that $(\Xi_N)_{11}=0$, so the r.h.s. is
$(\Xi_N)_{n,n} \circ m^Q_n = {1\over {n!}}\breve m^Q_n \circ
\Delta^Q_n \circ m^Q_n = \breve m^Q_n
\circ \on{Sym}_n \circ (p_1^Q)^{\boxtimes n} =
\breve m^Q_n \circ (p_1^Q)^{\boxtimes n}
\circ \on{Sym}_n =
(m^Q_n \circ (p_1^Q)^{\boxtimes n} \circ \on{Sym}_n)\breve{} =
\breve m^Q_n$, where
the second equality follows from the first part of (\ref{rel:4}),
the third equality follows from
the equality $\on{Sym}_n \circ
(p_1^S)^{\boxtimes n} = (p_1^S)^{\boxtimes n} \circ
\on{Sym}_n$ in $S({\bf Sch})$, the fourth equality follows from
$\breve p_Q^1 = 0$, and the last equality follows from
the second part of (\ref{eq:2}) and the first part of (\ref{eq:sym}).

The second part of (\ref{rel:4}) implies that
$$
(\Xi_N)_{nn} = - m^S_n \circ
\breve \Delta^S_n.
$$
We then prove that $\breve \Delta^Q_n = \sum_{i=1}^n (\Xi_N)_{11}^{(i)}
\circ \Delta^Q_n - \Delta^Q_n \circ (\Xi_N)_{nn}$ as above.

We also have $\breve p_n^Q = [(\Xi_N)_{11},p_n^Q]$ as $\breve p_n^Q = 0$
and $(\Xi_N)_{11} = 0$. It follows that if $(\Xi_N)_{diag} := \sum_{n\geq 1}
(\Xi_N)_{nn}$, then $\breve x^Q = (d/dt)_{|t=0} \on{Inn}(\on{id}_{{\bf id}}
 + t (\Xi_{N})_{diag})(x^Q)$ for $x^S \in \{m_n^S,\Delta_n^S,p_n^S\}$,
hence for any $x^S\in S({\bf Sch})(F,G)$ and any $F,G\in \on{Ob(Sch)}$.

If we now set $\Xi_N := (\Xi_N)_{diag} + (\Xi_N)_{off}$, we then get
$\dot x^Q = (d/dt)_{|t=0} \on{Inn}(\on{id}_{{\bf id}}
 + t \Xi_{N})(x^Q)$  for any $x^S\in S(\on{Sch})(F,G)$, as wanted.

The morphism $\Phi':= \on{Inn}(\Xi)^{-1} \circ \Phi : S({\bf P}) \to S({\bf Q})$
now satisfies $\Phi' \circ S(i_{\bf P}) = S(i_{\bf Q})$. To prove the proposition,
it remains to prove:

\begin{lemma}
Let $\Phi' : S({\bf P}) \to S({\bf Q})$ be a prop morphism such that
$\Phi' \circ S(i_{\bf P}) = S(i_{\bf Q})$, then there exists a prop
morphism $\varphi : {\bf P} \to {\bf Q}$ such that $\Phi' = S(\varphi)$.
\end{lemma}

{\em Proof of Lemma.} We have morphisms of Schur functors
$inj_1 : {\bf id} \to S$ and $pr_1 : S\to {\bf id}$, corresponding to the
direct sum $S = {\bf id} \oplus (\oplus_{i\neq 1} S^i)$.
These induce morphisms of props ${\bf P} = {\bf id}({\bf P}) \to S({\bf P})$
and $S({\bf Q}) \to {\bf id}({\bf Q}) = {\bf Q}$. We then define $\varphi$
as the composed morphism ${\bf P} \to S({\bf P}) \stackrel{\Phi'}{\to} S({\bf Q})
\to {\bf Q}$. If $F,G$ are any Schur functors, we want to prove that
the diagram
$$
\xymatrix@!0 @R=7ex @C=42ex{
S({\bf P})(F,G) \ar[r]^{\Phi'(F,G)} \ar[d]_{\scriptstyle{\sim}} & S({\bf Q})(F,G)
\ar[d]^{\scriptstyle{\sim}}\\
\oplus_{i,j} {\bf P}(S^i(F),S^j(G)) \ar[r]^{\oplus_{i,j} \varphi(S^i(F),S^j(G))}
& \oplus_{i,j} {\bf Q}(S^i(F),S^j(G))
}
$$
commutes, i.e., that for any pair $(i,j)$ of integers, the diagram
$$
\xymatrix@!0 @R=7ex @C=35ex{
{\bf P}(S^i(F),S^j(G)) \ar[r]^{ \varphi(S^i(F),S^j(G))} \ar[d]
& {\bf Q}(S^i(F),S^j(G)) \ar[d] \\
S({\bf P})(F,G) \ar[r]^{\Phi'(F,G)} & S({\bf Q})(F,G)
}
$$
commutes. If $F,G$ are any Schur functors, then the map
$\varphi(F,G) : {\bf P}(F,G) \to {\bf Q}(F,G)$ is the composition
${\bf P}(F,G) \to {\bf P}(S(F),S(G)) \simeq S({\bf P})(F,G)
\stackrel{\Phi'(F,G)}{\to}
S({\bf Q})(F,G) \simeq {\bf Q}(S(F),S(G))\to {\bf Q}(F,G)$, where the initial map
is $x \mapsto  i_{\bf P}(inj_1(G))\circ x \circ i_{\bf P}(pr_1(F))$ and the final
map is $y \mapsto i_{\bf Q}(pr_1(G))\circ y \circ i_{\bf Q}(inj_1(F))$.
So we want to prove that the diagram
$$
\xymatrix@!0 @R=11ex @C=19ex { &
S({\bf P})(S^i(F),S^j(G))
\ar[rr]^{\Phi'(S^i(F),S^j(G))} && S({\bf Q})(S^i(F),S^j(G))
\ar[dr]^{\scriptstyle{(b)}} \\
{\bf Q}(S^i(F),S^j(G)) \ar[ur]^{\scriptstyle{(a)}} \ar[dr]_{\scriptstyle{(c)}}
&&&&  {\bf Q}(S^i(F),S^j(G)) \\
   & S({\bf P})(F,G)\ar[rr]_{\Phi'(F,G)} && S({\bf Q})(F,G)
   \ar[ur]_{\scriptstyle{(d)}}
}
$$
commutes, where the map (a) is $x\mapsto i_{\bf P}(inj_1(S^j(G)))\circ x \circ
i_{\bf P}(pr_1(S^i(F)))$, the map (b) is $y\mapsto i_{\bf Q}(pr_1(S^j(G)))\circ y \circ
i_{\bf Q}(inj_1(S^i(F)))$, the map (c) is $x\mapsto inj_j(G) \circ x\circ pr_i(F)$,
and the map (d) is $x\mapsto pr_j(G) \circ x\circ inj_i(F)$.
Here $inj_i : S^i \hookrightarrow S$ and $pr_i : S \twoheadrightarrow S^i$
are the canonical injection and projection attached to $S
= \oplus_{j\geq 0} S^i$.

To prove the commutativity of this diagram, we construct a surjective
map $S({\bf P})(F^{\oplus i},G^{\oplus j}) \twoheadrightarrow {\bf P}(S^i,S^j)$
and an injective map ${\bf Q}(S^i,S^j) \hookrightarrow
S({\bf Q})(F^{\oplus i},G^{\oplus j})$, and prove that the external
diagram in
$$
\xymatrix@!0 @R=11ex @C=19ex { &
 S({\bf P})(S^i(F),S^j(G))  \ar[rr]^{\Phi'(S^i(F),S^j(G))} &&
 S({\bf Q})(S^i(F),S^j(G))\ar[d] \ar[dr]^{\scriptstyle{u'}}\\
S({\bf P})(F^{\oplus i},G^{\oplus j})\twoheadrightarrow
\ar[ur]^{\scriptstyle{u}} 
 \ar[dr]_{\scriptstyle{v}} & {\bf P}(S^i(F),S^j(G))\ar[u] \ar[d] & &
{\bf Q}(S^i(F),S^j(G))\hookrightarrow
& S({\bf Q})(F^{\oplus i},G^{\oplus j})\\
 &  S({\bf P})(F,G)  \ar[rr]_{\Phi'(F,G)} && S({\bf Q})(F,G)
 \ar[u] \ar[ur]_{\scriptstyle{v}}&
}
$$
commutes, where in this diagram the diagonal arrows are defined by the
condition that the triangles commute and the inner rectangle is the above
diagram. The commutativity of the external diagram then implies
that of the inner rectangle.

We first construct the maps $S({\bf P})(F^{\oplus i},G^{\oplus j})
\twoheadrightarrow {\bf P}(S^i(F),S^j(G))$ and ${\bf Q}(S^i(F),S^j(G)) \hookrightarrow
S({\bf Q})(F^{\oplus i},G^{\oplus j})$. There is a unique morphism of
Schur functors $a_i : S^i \hookrightarrow S \circ {\bf id}^{\oplus i}$,
given by the
composition $S^i \hookrightarrow {\bf id}^{\otimes i}
\stackrel{inj_1^{\boxtimes i}}{\hookrightarrow} S^{\otimes i} =
S \circ {\bf id}^{\oplus i}$ and
$b_i : S \circ {\bf id}^{\oplus i} \twoheadrightarrow S^i$,
given by the
composition $S \circ {\bf id}^{\oplus i} = S^{\otimes i}
\stackrel{pr_1^{\boxtimes i}}{\twoheadrightarrow}
{\bf id}^{\otimes i} \twoheadrightarrow S^i$.
Then $S({\bf P})(F^{\oplus i},G^{\oplus j})
\twoheadrightarrow {\bf P}(S^i(F),S^j(G))$ is
$x \mapsto b_j(G)\circ x \circ a_i(F)$ and
${\bf Q}(S^i(F),S^j(G)) \hookrightarrow
S({\bf Q})(F^{\oplus i},G^{\oplus j})$ is $y\mapsto a_j(G)\circ y \circ b_i(F)$.

We now compute the diagonal maps. We define
$\alpha_i \in S({\bf Sch})({\bf id}^{\oplus i},S^i) \simeq
{\bf Sch}(S \circ {\bf id}^{\oplus i},S \circ S^i)$ as the
composition $S \circ {\bf id}^{\oplus i} \simeq S^{\otimes i}
\stackrel{pr_1^{\boxtimes i}}{\twoheadrightarrow}
{\bf id}^{\otimes i} \twoheadrightarrow S^i
\stackrel{inj_1(S^i)}{\hookrightarrow} S\circ S^i$.
Similarly, we define:

$\bullet$ $\beta_i \in S({\bf Sch})({\bf id}^{\oplus i},{\bf id})$
as the composition $S \circ {\bf id}^{\oplus i} \simeq S^{\otimes i}
\stackrel{pr_1^{\boxtimes i}}{\twoheadrightarrow}
{\bf id}^{\otimes i} \twoheadrightarrow S^i
\stackrel{inj_i}{\hookrightarrow} S$;

$\bullet$ $\alpha'_i \in S({\bf Sch})(S^i,{\bf id}^{\oplus i})$
as the composition $S \circ S^i \stackrel{pr_1(S^i)}{\twoheadrightarrow} S^i
\hookrightarrow {\bf id}^{\otimes i}
\stackrel{inj_1^{\boxtimes i}}{\hookrightarrow} S^{\boxtimes i}\simeq S \circ
{\bf id}^{\oplus i}$;

$\bullet$ $\beta'_i \in S({\bf Sch})({\bf id},{\bf id}^{\oplus i})$
as the composition $S \circ {\bf id} \simeq S
\stackrel{pr_i}{\twoheadrightarrow} S^i
\hookrightarrow {\bf id}^{\otimes i}
\stackrel{inj_1^{\boxtimes i}}{\hookrightarrow} S^{\boxtimes i}\simeq S \circ
{\bf id}^{\oplus i}$.

The diagonal map $u$ is then $x\mapsto S(i_P)(\alpha_j(G))\circ x
\circ S(i_P)(\alpha'_i(F))$; the map $v$ is $x\mapsto
S(i_P)(\beta_j(G)) \circ x \circ S(i_P)(\beta'_i(F))$;
the map $u'$ is $y\mapsto S(i_Q)(\alpha'_j(G) \circ
y \circ S(i_Q)(\alpha_i(F))$; and the map $v'$ is
$y'\mapsto S(i_Q)(\beta'_j(G)) \circ y' \circ S(i_Q)(\beta_i(F))$.

We now prove the commutativity of the external diagram, i.e.,
that
$$
u' \circ \Phi'(S^i(F),S^j(G)) \circ u = v' \circ \Phi'(F,G) \circ v.
$$
Let $x\in S({\bf P})(F^{\oplus i},G^{\oplus j})$.  Then
\begin{align*}
& (u' \circ \Phi'(S^i(F),S^j(G)) \circ u)(x)
\\ & =
(u' \circ \Phi'(S^i(F),S^j(G))) \big( S(i_P)(\alpha_j(G))\circ x
\circ S(i_P)(\alpha'_i(F))\big)
\\ & = u' \big( S(i_{\bf Q})(\alpha_j(G)) \circ
\Phi'(F^{\oplus i},G^{\oplus j})(x)\circ S(i_{\bf Q})(\alpha'_i(F))\big)
\\ & = S(i_{\bf Q})(\alpha'_j(G) \circ S(i_{\bf Q})(\alpha_j(G)) \circ
\Phi'(F^{\oplus i},G^{\oplus j})(x)\circ S(i_{\bf Q})(\alpha'_i(F))
\circ S(i_{\bf Q})(\alpha_i(F))
\\ & = S(i_{\bf Q})(\alpha'_j \circ \alpha_j(G))
\circ \Phi'(F^{\oplus i},G^{\oplus j})(x)\circ
S(i_{\bf Q})(\alpha'_i \circ \alpha_i(F)),
\end{align*}
where the second equality follows from
$\Phi'\circ S(i_{\bf Q}) = S(i_{\bf Q})$. On the other hand,
\begin{align*}
& (v' \circ \Phi'(F,G) \circ v)(x)
\\ & =
(v' \circ \Phi'(F,G)) \big(
S(i_{\bf P})(\beta_j(G)) \circ x \circ S(i_{\bf P})(\beta'_i(F)) \big)
\\ & = v' \big(
S(i_{\bf Q})(\beta_j(G)) \circ \Phi'(F^{\oplus i},G^{\oplus j})(x)
\circ S(i_{\bf Q})(\beta'_i(F)) \big)
\\ & =
S(i_{\bf Q})(\beta'_j(G)) \circ  S(i_{\bf Q})(\beta_j(G)) \circ
\Phi'(F^{\oplus i},G^{\oplus j})(x)
\circ S(i_{\bf Q})(\beta'_i(F))\circ S(i_{\bf Q})(\beta_i(F))
\\ & = S(i_{\bf Q})(\beta'_j \circ\beta_j(G)) \circ
\Phi'(F^{\oplus i},G^{\oplus j})(x) \circ
S(i_{\bf Q})(\beta'_i \circ \beta_i(F)).
\end{align*}
To prove that these terms are equal, we will prove that
$$
\alpha'_i \circ \alpha_i = \beta'_i \circ \beta_i
$$
(an equality of endomorphisms of the Schur functor
$S \circ {\bf id}^{\oplus i}$). We have
$\alpha_i = inj_1(S^i)\circ b_i$, $\beta_i = inj_i\circ b_i$,
$\alpha'_i = a_i \circ pr_1(S^i)$, $\beta'_i = a_i \circ pr_i$,
so it suffices to show the equality $pr_1(S^i) \circ inj_1(S^i) =
pr_i \circ inj_i$. But $pr_i \circ inj_i = \on{id}_{S^i}$, and
$pr_1 \circ inj_1 = \on{id}_{{\bf id}}$, which implies that
$pr_1(S^i) \circ inj_1(S^i) = \on{id}_{S^i}$. This ends the proof
of the lemma, and hence of the proposition.
\hfill \qed \medskip

\subsection{Construction of commuting triangles of prop modules based on
${\bf LBA}$}

Let us then define $\hat\alpha^\pm : S(D_{add}({\bf LBA})) \to
S(M_{add}^\pm({\bf LBA}))$ by $\hat\alpha^\pm(X) := \Xi^{-1}_G \circ
\tilde\alpha^\pm (\on{Inn}(\Xi)(X))$ for $X\in S(D_{add}({\bf LBA}))(F,G)$,
and $\hat\beta^\pm : S({\bf LBA}) \to S(M_{add}^\pm({\bf LBA}))$ by
$\hat\beta^\pm(x) := \hat\alpha^\pm (S(\varphi)(x))$ for
$x\in S({\bf LBA})(F,G)$.

\begin{lemma}
$\hat\alpha^\pm$ is a prop left $S(D_{add}({\bf LBA}))$-module morphism
and $\hat\beta^\pm$ is a prop right $S({\bf LBA})$-module morphism.
\end{lemma}

{\em Proof.} We have for $X\in S(D_{add}({\bf LBA}))(F,G)$
and $Y\in S(D_{add}({\bf LBA}))(G,H)$, $\tilde\alpha^\pm(Y\circ X)
= Y \circ \tilde\alpha^\pm(X)$. Then $\hat\alpha^\pm(Y \circ X) =
\Xi_H^{-1} \circ \tilde\alpha^\pm (\Xi_H \circ Y \circ X \circ
\Xi_F^{-1}) =\Xi_H^{-1} \circ \tilde\alpha^\pm (Y' \circ X') =
\Xi_H^{-1} \circ Y' \circ \tilde\alpha^\pm (X') =
Y \circ \Xi_G^{-1} \circ \tilde \alpha^\pm(\Xi_G
\circ X \circ \Xi_F^{-1})= Y \circ \hat\alpha^\pm(X)$,
where $Y' = \Xi_H \circ Y \circ \Xi_G^{-1}$, $X' = \Xi_G
\circ X \circ \Xi_F^{-1}$. On the other hand, we have
for $x\in S({\bf LBA})(F,G)$, $\hat\beta^\pm(x) = \Xi_G^{-1} \circ
\tilde\beta^\pm(x)$, so the prop module properties of
$\tilde\beta^\pm$ imply that for $y\in S({\bf LBA})(F,G)$, we have
$\hat\beta^+(y\circ x) = \hat\beta^+(y) \circ (x \ul\boxtimes
\on{id}_{\bf 1}^{{\bf LBA}^*})$ and $\hat\beta^-(y\circ x)
= \hat\beta^+(y) \circ (\on{id}_{\bf 1}^{{\bf LBA}}
\ul\boxtimes *cop(x))$. \hfill \qed \medskip

So we have commuting triangles of prop morphisms
$$
\xymatrix 
{& S(M^\pm_{add}({\bf LBA})) & \\
S({\bf LBA}) \ar[ru]^{\hat\beta^\pm} \ar[rr]_{S(\varphi)} & &
S(D_{add}({\bf LBA}))\ar[ul]_{\hat\alpha^\pm}
}
$$

For $F\in \on{Sch}$, set $\hat\alpha_F^\pm :=
\hat\alpha^\pm(\on{id}_F^{S(D_{add}({\bf LBA}))})$,
$\hat\beta_F^\pm := \hat\beta^\pm(\on{id}_F^{S({\bf LBA})})$;
then
$$
\hat\alpha_F^\pm,\hat\beta_F^\pm\in S(M_{add}^\pm({\bf LBA}))(F,F).
$$

The lemma implies that if $x\in S({\bf LBA})(F,G)$, then
$\hat\beta^+(x) = \hat\beta^+_G \circ (x\ul\boxtimes
\on{id}_{\bf 1}^{{\bf LBA}^*})$,
$\hat\beta^-(x) = \hat\beta^-_G \circ
(\on{id}_{\bf 1}^{{\bf LBA}}\ul\boxtimes *cop(x))$,
and if $X\in S(D_{add}({\bf LBA}))(F,G)$, then
$\hat\alpha^\pm(X) = X \circ \hat\alpha^\pm_F$.

We have therefore
\begin{equation} \label{rel:alpha:beta}
\hat\beta_G^+ \circ (x\ul\boxtimes \on{id}_{\bf 1}^{{\bf LBA}^*})
= S(\varphi)(x) \circ \hat\alpha_F^+,
\quad \hat\beta_G^- \circ (\on{id}_{\bf 1}^{{\bf LBA}}\ul\boxtimes *cop(x))
= S(\varphi)(x) \circ \hat\alpha_F^-,
\end{equation}
for $x\in S({\bf LBA})(F,G)$.

For $F = G$ and $x = \on{id}_F^{S({\bf LBA})}$,
we get $\hat\beta_F^\pm = \hat\alpha_F^\pm$.

For $X_i\in S(D_{add}({\bf LBA}))(F_i,G_i)$
($i=1,2$), $\hat\alpha^\pm(X_1\boxtimes X_2) = \hat\alpha^\pm(X_1)\boxtimes
\hat\alpha^\pm(X_2)$ so $(X_1\boxtimes X_2) \circ \hat\alpha_{F\otimes G}^\pm
= (X_1 \circ \hat\alpha_F^\pm) \boxtimes (X_2 \circ \hat\alpha_G^\pm)$, so
$\hat\alpha^\pm_{F\otimes G} = \hat\alpha^\pm_F \boxtimes \hat\alpha^\pm_G$.
Moreover, for $f\in \on{Sch}(F,G)$, (\ref{rel:alpha:beta}) together with
$S(\varphi) \circ S(i_{\bf LBA}) = S(i_{D_{add}({\bf LBA})})$
implies $\hat\alpha_G^\pm \circ i_{S({\bf LBA})}(f) = i_{S(D_{add}({\bf LBA}))}(f)
\circ \hat\alpha_F^\pm$. This implies that the $\hat\alpha_F^\pm$ are obtained
from $\hat\alpha_{{\bf id}}^\pm$ by inflation (see Proposition \ref{inner}).

We now study $\hat\alpha_{{\bf id}}^\pm\in S(M_{add}^\pm({\bf LBA}))
({\bf id},{\bf id}) = M_{add}^\pm({\bf LBA})(S,S)$.

Let $inj_1 \in {\bf Sch}({\bf id},S)$ and $pr_1 \in {\bf Sch}(S,{\bf id})$
be the canonical morphisms. Set $\alpha_{\bf id}^\pm:=
pr_1^{S(D_{add}({\bf LBA}))} \circ \hat\alpha^\pm_{{\bf id}}\circ
(inj_1^{S({\bf LBA})}\ul\boxtimes \on{id}_{\bf 1}^{{\bf LBA}^*}) \in
M_{add}^\pm({\bf LBA})({\bf id},{\bf id})$. Then:

\begin{lemma} We have
\begin{equation} \label{form:hatalpha}
\hat\alpha_{{\bf id}}^\pm = (\alpha^\pm_{\bf id})_S,
\end{equation}
i.e.,
\begin{equation} \label{id:alpha+}
\hat\alpha_{{\bf id}}^+ = \sum_{k\geq 0}
(inj_k \circ  p_k)^{D_{add}({\bf LBA})} \circ
(\alpha_{\bf id}^+)^{\boxtimes k} \circ
((i_k \circ pr_k)^{{\bf LBA}}\ul\boxtimes \on{id}^{{\bf LBA}^*}_{\bf 1}),
\end{equation}
\begin{equation} \label{id:alpha-}
\hat\alpha_{{\bf id}}^- = \sum_{k\geq 0}
(inj_k \circ  p_k)^{D_{add}({\bf LBA})} \circ
(\alpha_{\bf id}^-)^{\boxtimes k} \circ
(\on{id}_{\bf 1}^{{\bf LBA}} \ul\boxtimes(i_k \circ pr_k)^{{\bf LBA}^*}),
\end{equation}
where $S \stackrel{pr_k}{\twoheadrightarrow} S^k
\stackrel{i_k}{\hookrightarrow} {\bf id}^{\otimes k}$,
${\bf id}^{\otimes k} \stackrel{p_k}{\twoheadrightarrow} S^k
\stackrel{inj_k}{\hookrightarrow} S$ are the natural morphisms
in ${\bf Sch}$ (with $p_k \circ i_k = pr_k \circ inj_k = \on{id}_{S^k}$).
\end{lemma}

{\em Proof.} We prove (\ref{id:alpha+}).
For $x = S(i_{{\bf LBA}})(p_n^S)\in S({\bf LBA})({\bf id},{\bf id})$,
(\ref{rel:alpha:beta}) together with
$S(\varphi) \circ S(i_{{\bf LBA}}) = S(i_{D_{add}({\bf LBA})})$ gives
$\hat\alpha_{{\bf id}}^+ \circ (S(i_{{\bf LBA}})(p_n^S)\ul\boxtimes
\on{id}_{\bf 1}^{{\bf LBA}^*}) =
S(i_{D_{add}({\bf LBA})})(p_n^S) \circ \hat\alpha_{{\bf id}}^+$.
Recall that $\hat\alpha_{{\bf id}}^+$ decomposes uniquely
as $\sum_{n,m\geq 0} inj_m^{D_{add}({\bf LBA})} \circ \alpha_{n,m}
\circ (pr_n^{{\bf LBA}} \ul\boxtimes \on{id}_{\bf 1}^{{\bf LBA}^*})$,
where $\alpha_{n,m}\in
M^+_{add}({\bf LBA})(S^n,S^m)$; then the above identity implies
$\sum_{q\geq 0} inj_q^{D_{add}({\bf LBA})} \circ \alpha_{n,q}
\circ (pr_n^{{\bf LBA}} \ul\boxtimes \on{id}_{\bf 1}^{{\bf LBA}^*})
= \sum_{p\geq 0} inj_n^{D_{add}({\bf LBA})}
\circ \alpha_{p,n}\circ (pr_p^{{\bf LBA}}\ul\boxtimes
\on{id}_{\bf 1}^{{\bf LBA}^*})$, so $\alpha_{n,m}=0$
for $n\neq m$. This implies that $\hat\alpha_{{\bf id}}^+$ has
the form $\sum_{n\geq 0} inj_n^{D_{add}({\bf LBA})} \circ \alpha_n
\circ (pr_n^{{\bf LBA}}\ul\boxtimes \on{id}_{\bf 1}^{{\bf LBA}^*})$,
where $\alpha_n\in M^+_{add}({\bf LBA})(S^n,S^n)$.

For $x = S(i_{{\bf LBA}})(m_n^S)$, (\ref{rel:alpha:beta}) gives
$\hat\alpha_{{\bf id}}^+ \circ S(i_{{\bf LBA}})(m_n^S) =
S(i_{D_{add}({\bf LBA})})(m_n^S) \circ
(\hat\alpha_{{\bf id}}^+)^{\boxtimes n}$.
Composing this equality with $inj_n^{{\bf LBA}}
\ul\boxtimes\on{id}_{\bf 1}^{{\bf LBA}^*}$
from the right and with $pr_n^{D_{add}({\bf LBA})}$
from the left, we get
$$
\alpha_n \circ (p_n^{{\bf LBA}}\ul\boxtimes \on{id}_{\bf 1}^{{\bf LBA}^*})
= p_n^{D_{add}({\bf LBA})} \circ (\alpha_{\bf id}^+)^{\boxtimes n}
$$
(an identity in $M_{add}^+({\bf LBA})({\bf id}^{\otimes n},S^n)$).
Composing from the right by $(i_n \circ pr_n)^{{\bf LBA}} \ul\boxtimes
\on{id}_{\bf 1}^{{\bf LBA}^*}$ and from the left by $inj_n^{D_{add}({\bf LBA})}$,
we find that $inj_n^{D_{add}({\bf LBA})} \circ \alpha_n \circ
(pr_n^{{\bf LBA}}\ul\boxtimes \on{id}_{\bf 1}^{{\bf LBA}^*})$ is the
$n$th term of the r.h.s. of (\ref{id:alpha+}), which implies this identity.
The identity on
$\hat\alpha_{\bf id}^-$ is proved in the same way, which implies
(\ref{form:hatalpha}). \hfill \qed \medskip

By computing classical limits, one shows that $\alpha_{\bf id}^\pm 
= \on{can}^\pm_{\bf id}$ + terms of higher degree, and similarly 
$\beta_{\bf id}^\pm = \on{can}^\pm_{\bf id}$ + terms of higher degree
($\on{can}_F^\pm$ are defined after Lemma \ref{lemma:double}).

Define $\alpha_F^\pm \in M_{add}^\pm({\bf LBA})(F,F)$ by applying inflation
(Proposition \ref{inner}) to $\alpha_{\bf id}^\pm\in
M^\pm_{add}({\bf LBA})({\bf id},{\bf id})$. Since we have
$\alpha_{F\otimes G}^\pm = \alpha_F^\pm \boxtimes \alpha_G^\pm$, the maps
$D_{add}({\bf LBA})(F,G) \ni X\mapsto X \circ \alpha_F^\pm\in
M_{add}^\pm({\bf LBA})(F,G)$
and ${\bf LBA}(F,G) \ni x \mapsto \alpha_G^+ \circ
(x\ul\boxtimes \on{id}_{\bf 1}^{{\bf LBA}^*}) \in D_{add}({\bf LBA})(F,G)$,
resp., ${\bf LBA}(F,G) \ni x \mapsto \alpha_G^+ \circ
(\on{id}_{\bf 1}^{{\bf LBA}} \ul\boxtimes cop\circ *(x))\in
D_{add}({\bf LBA})(F,G)$ define a left prop $D_{add}({\bf LBA})$-module
morphism $\alpha^\pm : D_{add}({\bf LBA}) \to M_{add}^\pm({\bf LBA})$
and a right prop ${\bf LBA}$-module morphism
$\beta^\pm : {\bf LBA} \to M_{add}^\pm({\bf LBA})$.

\begin{lemma}
We have $\hat\alpha^\pm = S(\alpha^\pm)$,
$\hat\beta^\pm = S(\beta^\pm)$.
\end{lemma}

{\em Proof.} This follows from the fact that if $F\in \on{Ob(Sch)}$,
the identification $S(M^\pm_{add}({\bf LBA}))(F,F) \simeq
M^\pm_{add}({\bf LBA})(F\circ S,F\circ S)$ takes $\hat\alpha^\pm_F$
to $\alpha^\pm_{F\circ S}$. Let $\gamma^\pm_F\in
M^\pm_{add}({\bf LBA})(F\circ S,F\circ S)$ be the image of $\hat\alpha^\pm_F$.
Then (\ref{form:hatalpha}) implies that $\gamma^\pm_{\bf id} = \alpha_S^\pm$.
On the other hand, the relations satisfied by the $\hat\alpha_F^\pm$
imply that $\gamma^\pm_{F\boxtimes S} = \gamma^\pm_F \boxtimes \gamma^\pm_G$,
and if $f\in \on{Sch}(F,G)$, then $(f\circ S)^{D_{add}({\bf LBA})}
\circ \gamma_F^+ = \gamma_G^+ \circ
( (f\circ S)^{{\bf LBA}}\ul\boxtimes \on{id}_{\bf 1}^{{\bf LBA}^*})$,
$(f\circ S)^{D_{add}({\bf LBA})} \circ \gamma_F^+ = \gamma_G^+ \circ
( \on{id}_{\bf 1}^{{\bf LBA}}\ul\boxtimes(f\circ S)^{{\bf LBA}^*})$.
All these conditions are satisfied by the family $(\alpha^\pm_{F\circ S})_{F\in
\on{Ob(Sch)}}$, and since they determine the family $(\gamma_F^\pm)_{F\in
\on{Ob(Sch)}}$ uniquely, we obtain  $\gamma^\pm_{F} = \alpha_{F\circ S}^\pm$,
as wanted. \hfill \qed \medskip

\begin{lemma}
The diagrams of prop module morphisms
$$
\xymatrix 
{& M^\pm_{add}({\bf LBA}) & \\
{\bf LBA} \ar[ru]^{\beta^\pm} \ar[rr]_{\varphi} & &
D_{add}({\bf LBA})\ar[ul]_{\alpha^\pm}
}
$$
commute.
\end{lemma}

{\em Proof.} Let $F,G\in \on{Ob(Sch)}$. We know that
the diagram
$$
\xymatrix 
{& S(M^\pm_{add}({\bf LBA}))(F,G) & \\
S({\bf LBA})(F,G) \ar[ru]^{\hat\beta^\pm(F,G)} \ar[rr]_{S(\varphi)(F,G)} & &
S(D_{add}({\bf LBA}))(F,G)\ar[ul]_{\hat\alpha^\pm(F,G)}
}
$$
commutes, i.e., the inner triangle of the following diagram
$$
\xymatrix  @!0 @R=8ex @C=19ex{
 & & M_{add}^\pm({\bf LBA})(F,G)\ar[d] \\
  & & M^\pm_{add}({\bf LBA})(F\circ S,G\circ S) \\
  & {\bf LBA}(F\circ S,G\circ S)\ar[ur]^{\beta^\pm(F\circ S,G\circ S)}
  \ar[rr]_{\varphi(F\circ S,G\circ S)} & &
  D_{add}({\bf LBA})(F\circ S,G\circ S)\ar[ul]_{\; \alpha^\pm(F\circ S,G\circ S)}
  \\
  {\bf LBA}(F,G)\ar[rrrr]_{\varphi(F,G)}\ar@/^{2pc}/[uuurr]^{\beta^\pm(F,G)}\ar[ur]
  & & & & D_{add}({\bf LBA})(F,G)\ar[ul]\ar@/_{2pc}/[uuull]_{\alpha^\pm(F,G)} }
$$
commutes. Here the maps between the triangles are
$$
{\bf LBA}(F,G) \to {\bf LBA}(F\circ S,F\circ G), \quad
x\mapsto G(inj_1)^{{\bf LBA}} \circ x \circ (F(pr_1))^{{\bf LBA}},
$$
$$D_{add}({\bf LBA})(F,G) \to D_{add}({\bf LBA})(F\circ S,F\circ G),
\quad
X\mapsto G(inj_1)^{D_{add}({\bf LBA})} \circ X \circ
(F(pr_1))^{D_{add}({\bf LBA})},
$$
$$
M_{add}^+({\bf LBA})(F,G) \to M_{add}^+({\bf LBA})(F\circ S,F\circ G), \;
m\mapsto G(inj_1)^{D_{add}({\bf LBA})} \circ m \circ
((F(pr_1))^{{\bf LBA}} \ul\boxtimes \on{id}_{\bf 1}^{{\bf LBA}^*}),
$$
and
$$
M_{add}^-({\bf LBA})(F,G) \to M_{add}^-({\bf LBA})(F\circ S,F\circ G), \;
m\mapsto G(inj_1)^{D_{add}({\bf LBA})} \circ m \circ
(\on{id}_{\bf 1}^{{\bf LBA}} \ul\boxtimes (F(pr_1))^{{\bf LBA}^*}).
$$
Since these maps are injective, the commutativity of the outer triangle
follows from that of the three side quadrilaterals. We now check
the commutativity of these diagrams. The bottom diagram obviously commutes.
The commutativity of the upper right diagram follows from the
identities
\begin{align} \label{ids:right}
& (F(pr_1))^{D_{add}({\bf LBA})} \circ \alpha_{F\circ S}^+ =
\alpha_F^+ \circ ( (F(pr_1))^{{\bf LBA}}\ul\boxtimes
\on{id}_{\bf 1}^{{\bf LBA}^*}),
\\ & \nonumber
(F(pr_1))^{D_{add}({\bf LBA})} \circ \alpha_{F\circ S}^- =
\alpha_F^- \circ ( \on{id}_{\bf 1}^{{\bf LBA}} \ul\boxtimes
(F(pr_1))^{{\bf LBA}^*}),
\end{align}
and the commutativity of the upper left diagram follows
from the identities
\begin{align} \label{ids:left}
& (G(inj_1))^{D_{add}({\bf LBA})} \circ \alpha_{G}^+ =
\alpha_{G\circ S}^+ \circ ( (G(inj_1))^{{\bf LBA}}\ul\boxtimes
\on{id}_{\bf 1}^{{\bf LBA}^*}),
\\ & \nonumber
(G(inj_1))^{D_{add}({\bf LBA})} \circ \alpha_{G}^- =
\alpha_{G\circ S}^- \circ ( \on{id}_{\bf 1}^{{\bf LBA}} \ul\boxtimes
(G(inj_1))^{{\bf LBA}^*}).
\end{align}
We now prove (\ref{ids:right}), (\ref{ids:left}). When $F = {\bf id}$
or $G = {\bf id}$, these identities follow from the equations
(\ref{id:alpha+}), (\ref{id:alpha-}) relating $\alpha_S^\pm$ and
$\alpha_{\bf id}^\pm$. Let us prove the first identity of (\ref{ids:right}).
We denote its l.h.s. by $\xi_F$ and its r.h.s. by $\eta_F$. Both
are elements of $M^+_{add}({\bf LBA})(F\circ S,F)$; $(\xi_F)_{F\in \on{Sch}}$
satisfies the relations $\xi_{F\otimes G} = \xi_F \boxtimes \xi_G$
and for $f\in \on{Sch}(F,G)$, $f^{D_{add}({\bf LBA})}\circ \xi_F =
\xi_G \circ (f^{S({\bf LBA})}\ul\boxtimes\on{id}_{\bf 1}^{{\bf LBA}^*})$,
and $(\eta_F)_{F\in \on{Sch}}$ satisfies the same relations; we also have
$\xi_{\bf id} = \eta_{\bf id}$. Applying the uniqueness proposition \ref{inner}
to the biprop $(F,G) \mapsto M^+_{add}(\bf LBA)(F\circ S,G)$ (over the props
$D^{add}({\bf LBA})$ and $S({\bf LBA})$),  we get $\xi_F = \eta_F$ for any $F$.
The other identities are proved in the same way.
\hfill \qed \medskip

\subsection{The relation between $\varphi$ and $\on{double}$}

In particular, we have for $x\in {\bf LBA}(F,G)$,
$$
\varphi(x) \circ \alpha^+_F = \alpha^+_G \circ (x\ul\boxtimes
\on{id}^{{\bf LBA}^*}_{\bf 1}),
\quad
\varphi(x) \circ \alpha^-_F = \alpha^-_G \circ (\on{id}_{\bf 1}^{{\bf LBA}}
\ul\boxtimes *cop(x)).
$$
We construct $\gamma \in D_{add}({\bf LBA})({\bf id},{\bf id})^\times$
such that $\gamma \circ \alpha_{\bf id}^\pm = \on{can}^\pm_{\bf id}$, where
$\on{can}_{\bf id}^\pm \in M_{add}^\pm({\bf LBA})({\bf id},{\bf id})$
have been defined in Section \ref{sec:double:LBA}.

We have canonical identifications $M_{add}^+({\bf LBA})({\bf id},{\bf id}) \simeq
{\bf LBA}({\bf id},{\bf id})$ and $M_{add}^-({\bf LBA})({\bf id},{\bf id}) \simeq
{\bf LBA}^*({\bf id},{\bf id})$; they send $\on{can}_{\bf id}^\pm$ to
$\on{id}^{{\bf LBA}}_{\bf id}$, resp., $\on{id}^{{\bf LBA}^*}_{\bf id}$.
Let us again denote by $\alpha_{\bf id}^\pm$ the images of
$\alpha_{\bf id}^\pm$ in ${\bf LBA}({\bf id},{\bf id})$, resp., in
${\bf LBA}^*({\bf id},{\bf id})$. These elements expand as
$\on{id}^{{\bf LBA}}_{\bf id}$ + higher terms, resp.,
$\on{id}^{{\bf LBA}^*}_{\bf id}$ + higher terms, hence they are invertible.
Let us set
$$
\gamma := (\alpha^+_{\bf id})^{-1} \ul\boxtimes \on{id}^{{\bf LBA}^*}_{\bf 1}
+ \on{id}^{{\bf LBA}}_{\bf 1} \ul\boxtimes (\alpha^-_{\bf id})^{-1},
$$
then $\gamma\in ({\bf LBA}\ul\boxtimes {\bf LBA}^*)(\ul\Delta({\bf id}),
\ul\Delta({\bf id})) \subset D_{add}({\bf LBA})({\bf id},{\bf id})$
is invertible, and satisfies
$\gamma \circ \alpha_{\bf id}^\pm = \on{can}^\pm_{\bf id}$.

We then set $\hat D := \on{Inn}(\gamma) \circ \varphi$. Then
$\hat D : {\bf LBA} \to D_{add}({\bf LBA})$ is a prop morphism.
For any $x\in {\bf LBA}(F,G)$, $\hat D(x)\in
{\bf LBA}_2(\ul\Delta(F),\ul\Delta(G))$ and
$$
\hat D(x) \circ \on{can}^+_F = \on{can}^+_G \circ (x\ul\boxtimes
\on{id}_{{\bf 1}}^{{\bf LBA}^*}), \quad
\hat D(x) \circ \on{can}^-_F = \on{can}^-_G \circ (\on{id}_{{\bf
1}}^{{\bf LBA}} \ul\boxtimes *cop(x)),
$$
(the first identity is in ${\bf LBA}_2(F\ul\boxtimes {\bf 1},\ul\Delta(G))$
and the second identity is in ${\bf LBA}_2({\bf 1}\ul\boxtimes F,\ul\Delta(G))$).

Indeed, $\hat D(x) \circ \on{can}^+_F = \gamma_G \circ \varphi(x) \circ
\gamma_F^{-1} \circ \on{can}^+_F = \gamma_G \circ \varphi(x) \circ \alpha_F^+
= \gamma_G \circ \alpha_G^+ \circ (x\ul\boxtimes \on{id}_{\bf 1}^{{\bf LBA}^*}) =
\on{can}_G^+ \circ (x\ul\boxtimes \on{id}_{\bf 1}^{{\bf LBA}^*})$; the second
identity is proved in the same way.

This means that we have a commutative diagram of prop modules
$$
\xymatrix 
{& M^\pm_{add}({\bf LBA}) & \\
{\bf LBA} \ar[ru]^{\beta_{can}^\pm} \ar[rr]_{\hat D} & &
D_{add}({\bf LBA})\ar[ul]_{\alpha_{can}^\pm} }
$$

We will show:

\begin{prop} For $\lambda$ a scalar, set $Z_\lambda :=
\on{exp}(\lambda(\mu\circ\delta))\in {\bf LBA}({\bf id},{\bf id})$.
There exists a scalar $\lambda$ such that
$\hat D = \on{Inn}( \on{id}_{\bf id}^{{\bf LBA}} \ul\boxtimes
\on{id}_{\bf 1}^{{\bf LBA}^*}
+ \on{id}_{\bf 1}^{{\bf LBA}} \ul\boxtimes *(Z_{\lambda})) \circ \on{double}$.
\end{prop}

{\em Proof of Proposition.}
$\hat D$ is uniquely determined by $\hat D(\mu)
\in D_{add}({\bf LBA})(\wedge^2,{\bf id}) = {\bf LBA}_2(
\wedge^2 \ul\boxtimes {\bf 1} \oplus {\bf id} \ul\boxtimes {\bf id} \oplus
{\bf 1} \ul\boxtimes \wedge^2, {\bf id} \ul\boxtimes {\bf 1}
\oplus {\bf 1}\ul\boxtimes {\bf id})$ and
$\hat D(\delta) \in D_{add}({\bf id},\wedge^2)
= {\bf LBA}_2({\bf id} \ul\boxtimes {\bf 1} \oplus {\bf 1}\ul\boxtimes {\bf id},
\wedge^2 \ul\boxtimes {\bf 1} \oplus {\bf id} \ul\boxtimes {\bf id} \oplus
{\bf 1}\ul\boxtimes \wedge^2) = {\bf LBA}_2({\bf id}\ul\boxtimes {\bf 1},
\wedge^2\ul\boxtimes {\bf 1}) \oplus {\bf LBA}_2({\bf 1}\ul\boxtimes {\bf id},
{\bf 1}\ul\boxtimes \wedge^2)$.

In the finite dimensional case, the semiclassical limit of the
Drinfeld double of the quantization of a Lie bialgebra $\A$
is ${\mathfrak D}(\A)$. The propic counterpart of this fact is that
we have the expansions $\hat D(\mu) = \on{double}(\mu)$ + terms of
higher degree, $\hat D(\delta) = \on{double}(\delta)$ + terms of
higher degree.

We decompose $\hat D(\mu)$ as $\mu_{\wedge^2\ul\boxtimes {\bf 1},{\bf
id}\ul\boxtimes {\bf 1}}
+ \mu_{{\bf id}\ul\boxtimes {\bf id},{\bf id}\ul\boxtimes {\bf
1}} + \mu_{{\bf id}\ul\boxtimes {\bf id},{\bf 1}\ul\boxtimes {\bf
id}} + \mu_{{\bf 1}\ul\boxtimes {\bf id}, {\bf 1}
\ul\boxtimes \wedge^2}$, where $\mu_{\wedge^2\ul\boxtimes {\bf 1},{\bf
id}\ul\boxtimes {\bf 1}} \in {\bf LBA}_2({\wedge^2\ul\boxtimes {\bf 1},{\bf
id}\ul\boxtimes {\bf 1}})$, etc.
Similarly, we decompose $\hat D(\delta)$ as
$\delta_{{\bf id}\ul\boxtimes{\bf 1},\wedge^2\ul\boxtimes {\bf 1}}
+ \delta_{{\bf 1}\ul\boxtimes {\bf id},{\bf 1}\ul\boxtimes \wedge^2}$.

The relation $\hat D(\mu) \circ \on{can}^+_{\wedge^2}
= \on{can}^+_{{\bf id}} \circ \mu$ gives
$\mu_{\wedge^2\ul\boxtimes {\bf 1}, {\bf id}\ul\boxtimes {\bf 1}} =
\mu$. Similarly,  $\hat D(\mu) \circ \on{can}^-_{\wedge^2}
= \on{can}^-_{{\bf id}} \circ \mu$  gives
$\mu_{{\bf 1}\ul\boxtimes \wedge^2, {\bf 1}\ul\boxtimes {\bf id}} =
\delta$.

The relation $\on{can}^+_{\wedge^2} \circ \delta = \hat D(\delta) \circ
\on{can}^+_{{\bf id}}$ gives
$\delta_{{\bf id}\ul\boxtimes {\bf 1}, \wedge^2\ul\boxtimes {\bf 1}} =
\delta$; similarly,  $\on{can}^-_{\wedge^2} \circ \delta = \hat D(\delta) \circ
\on{can}^-_{{\bf id}}$ gives $\delta_{{\bf 1}\ul\boxtimes {\bf id},
{\bf 1}\ul\boxtimes \wedge^2} = -\mu$.

Let us set $\rho:=\mu_{{\bf id}\ul\boxtimes{\bf id},
{\bf 1}\ul\boxtimes{\bf id}}$, $\check\rho := \mu_{{\bf id}\ul\boxtimes
{\bf id},{\bf id}\ul\boxtimes {\bf 1}}$;
we identify $-\rho$ with an element $\sigma\in {\bf LBA}(T_2,{\bf id})$ and
$\check\rho$ with an element $\check\sigma\in{\bf LBA}({\bf id},T_2)$.
We will write some equations satisfied by $\sigma$ and $\check\sigma$.

If $(\A,\mu_\A,\delta_\A)$ is a finite dimensional Lie bialgebra, then its
double is ${\mathfrak D}(\A)= \A \oplus \A^*$ with the bracket
$$
[x,y]_{{\mathfrak D}(\A)} = [x,y]_\A, \; [\xi,\eta]_{{\mathfrak D}(\A)} 
= [\xi,\eta]_{\A^*}, \; [x,\xi]_{{\mathfrak D}(\A)} 
= \on{ad}^*_x(\xi) - \underline{\on{ad}}^*_\xi(x)
$$
where $[x,y]_\A = \mu_\A(x,y)$, $[-,-]_\A^* = {}^t\delta_\A$,
$\on{ad}^*$ is the coadjoint action of $(\A,[-,-]_\A)$ on $\A^*$,
$\underline{\on{ad}}^*$ is the coadjoint action of $(\A^*,[-,-]_{\A^*})$ on
$\A$; and with the cobracket
$$
\delta_{{\mathfrak D}(\A)}(x) = \delta_\A(x), \quad 
\delta_{{\mathfrak D}(\A)}(\xi) = -{}^t\mu_\A(\xi);
$$
here $x,y\in \A$ and $\xi,\eta\in\A^*$.

On the other hand, $\hat D$ gives rise to a composed prop morphism
${\bf LBA} \to D_{add}({\bf LBA}) \to \on{Prop}(\A\oplus \A^*)$, hence to a
Lie bialgebra structure on $\A\oplus \A^*$, denoted
$(\mu'_{{\mathfrak D}(\A)},\delta'_{{\mathfrak D}(\A)})$. 
The brackets and cobrackets for this structure are the same as above, except for
$$
\mu'_{{\mathfrak D}(\A)}(x,\xi) = \rho^\A_x(\xi) + \check\rho^\A_\xi(x),
$$
where $\rho^\A : \A\otimes \A^* \to \A^*$, $x\otimes\xi\mapsto \rho^\A_x(\xi)$
and $\check\rho^\A : \A\otimes \A^* \to \A$,
$x\otimes\xi\mapsto \check\rho^\A_\xi(x)$
are the realizations of $\rho$ and $\check\rho$. The realization of
$\sigma$ is then $\sigma^\A : \A^{\otimes 2} \to \A$,
$x\otimes y \mapsto \sigma^\A_x(y)$, such that for any $x\in\A$,
$\sigma^\A_x = - {}^t(\rho^\A_x)$.

The fact that $\mu'_{{\mathfrak D}(\A)}$ satisfies the Jacobi identity implies that
$\sigma^\A_{[x,y]} = [\sigma^\A_x,\sigma^\A_y]$ (identity in $\on{End}(\A)$)
and the cocycle identity for the bracket $\mu'_{{\mathfrak D}(\A)}(x,\xi)$ implies
$\sigma^\A_x([y,z]) = [\sigma^\A_x(y),z] + [y,\sigma^\A_x(z)]$
for any $x,y,z\in \A$.

At the propic level, one shows that the Jacobi identity satisfied by
$\hat D(\mu)$ and the cocycle identity satisfied by $(\hat D(\mu),
\hat D(\delta))$ imply that $\sigma$ satisfies the universal versions
of these identities
$$
\sigma \circ (\mu\boxtimes \on{id}_{\bf id}^{{\bf LBA}})
= \sigma \circ (\on{id}_{\bf id}^{{\bf LBA}} \boxtimes \sigma)
- \sigma \circ (\on{id}_{\bf id}^{{\bf LBA}} \boxtimes \sigma) \circ (213),
$$
$$
\sigma \circ (\on{id}_{\bf id}^{{\bf LBA}}\boxtimes \mu)
= \mu \circ (\sigma \boxtimes \on{id}_{\bf id}^{{\bf LBA}})
- \mu \circ (\sigma \boxtimes \on{id}_{\bf id}^{{\bf LBA}}) \circ (132).
$$
$\sigma\in {\bf LBA}(T_2,{\bf id})$ has the expansion $\mu$ + terms of higher
degree. Let $\sigma'$ be the nonzero homogeneous term of smallest
possible degree of $\sigma - \mu$ if $\sigma - \mu\neq 0$, or $0$ if
$\sigma=\mu$.  Then $\sigma'$
satisfies the equations
$$
\sigma' \circ (\on{id}_{{\bf id}}^{{\bf LBA}} \boxtimes \mu)
= \mu \circ (\sigma'
\boxtimes \on{id}_{{\bf id}}^{{\bf LBA}}) - \mu \circ (\sigma'
\boxtimes \on{id}_{{\bf id}}^{{\bf LBA}}) \circ (132),
$$
$$
\sigma' \circ (\mu\boxtimes \on{id}_{{\bf id}}^{{\bf LBA}}) = \mu \circ
(\on{id}_{\bf id}^{{\bf LBA}} \boxtimes \sigma') + \sigma' \circ
(\on{id}_{\bf id}^{{\bf LBA}} \boxtimes \mu)- \mu \circ
(\on{id}_{\bf id}^{{\bf LBA}} \boxtimes \sigma') \circ (213)
- \sigma' \circ (\on{id}_{\bf id}^{{\bf LBA}} \boxtimes \mu) \circ (213).
$$

\begin{lemma} \label{lemma:sigma'}
Let $\sigma'\in {\bf LBA}({\bf id}^{\otimes 2},{\bf id})$ satisfy the above
identities; then $\sigma'$ is proportional to
$\mu \circ ((\mu\circ \delta) \boxtimes \on{id}_{\bf id}^{{\bf LBA}})$.
\end{lemma}

{\em Proof of Lemma.} These identities imply
\begin{equation} \label{main:id:theta}
\sigma' \circ (\mu \boxtimes \on{id}_{{\bf id}}^{{\bf LBA}}) =
\mu \circ \Big( \big( \sigma' - \sigma' \circ (21) \big)
\boxtimes \on{id}_{{\bf id}}^{{\bf LBA}} \Big).
\end{equation}
The proof is a propic version of the following argument.
The above identities are the universal versions of
\begin{equation}\label{der}
\sigma^{\prime\A}_{x}([y,z]) = [\sigma^{\prime\A}_{x}(y),z]
+ [y,\sigma^{\prime\A}_{x}(z)]
\end{equation}
(identity of maps $\A^{\otimes 3} \to \A$) and
$$
\sigma^{\prime\A}_{[x,y]} = [\sigma^{\prime\A}_{x},\on{ad}_y]
- [\sigma^{\prime\A}_{y},\on{ad}_x]
$$
(identity in $\on{End}(\A)$, where $\on{ad}_x(y) = [x,y]$) which is written
$$
\sigma^{\prime\A}_{[x,y]}(z) = \sigma^{\prime\A}_{x}([y,z])
- [y,\sigma^{\prime\A}_{x}(z)]
- \big( \sigma^{\prime\A}_{y}([x,z]) - [x,\sigma^{\prime\A}_{y}(z)]\big).
$$
Taking into account (\ref{der}), this gives
$$
\sigma^{\prime\A}_{[x,y]}(z) = [\sigma^{\prime\A}_{x}(y)
- \sigma^{\prime\A}_{y}(x),z].
$$
We have:

\begin{fact} \label{lemma:inject}
If $A$ is any Schur functor, the map $\on{LBA}({\bf id}\otimes A,{\bf id})
\to \on{LBA}(T_2 \boxtimes A,{\bf id})$,
$x\mapsto x\circ (\mu \boxtimes \on{id}_{{\bf id}}^{\on{LBA}})$ is injective.
\end{fact}

{\em Proof of Fact.}
The composed map
$\oplus_{Z,W\in\on{Irr(Sch)}} \on{LCA}({\bf id},W) \otimes \on{LCA}(A,Z)
\otimes \on{LA}(W\otimes Z,{\bf id}) \simeq
\on{LBA}({\bf id} \otimes A,{\bf id}) \to
\on{LBA}(T_2 \otimes A,{\bf id}) \simeq \oplus_{Z_1,Z,W\in\on{Irr(Sch)}}
\on{LCA}({\bf id},Z_1) \otimes \on{LCA}({\bf id},Z) \otimes \on{LCA}(A,W)
\otimes \on{LA}(Z_1 \otimes Z \otimes W,{\bf id}) \twoheadrightarrow
\oplus_{Z,W\in\on{Irr(Sch)}} \on{LCA}({\bf id},Z) \otimes \on{LCA}(A,W)
\otimes \on{LA}({\bf id} \otimes Z \otimes W,{\bf id})$,
where the last projection is onto the sum of
components with $Z_1 = {\bf id}$, is the direct sum over $Z,W$ of
the maps $\on{LA}(Z\otimes W,{\bf id}) \to \on{LA}({\bf id} \otimes Z \otimes
W,{\bf id})$, $\lambda \mapsto \lambda \circ (\mu_Z
\boxtimes \on{id}_W^{\on{LBA}})$, with the identity; here
$\mu_Z\in \on{LA}({\bf id} \otimes Z,Z)$
is induced by $\mu_{|Z|} \in \on{LA}({\bf id} \otimes T_{|Z|},T_{|Z|})$
given by $x\otimes x_1\otimes .. \otimes x_{|Z|} \mapsto
\sum_{i=1}^{|Z|} x_1 \otimes ... \otimes [x,x_i] \otimes ... \otimes x_{|Z|}$.
It follows from the first statement of \cite{EH}, Prop. 3.2, that
this map is injective. It follows that the map
$\on{LBA}({\bf id} \otimes A,{\bf id}) \to
\on{LBA}(T_2 \otimes A,{\bf id})$ is injective.
\hfill \qed \medskip

We now prove the following fact:

\begin{fact} \label{fact}
Assume that $\psi_1,\psi_2\in \on{LBA}(T_2,{\bf id})$ are such that
$\psi_1 \circ (\mu \boxtimes \on{id}_{{\bf id}}^{\on{LBA}}) = \mu \circ
(\psi_2\boxtimes \on{id}_{{\bf id}})$, then there exists $\psi
\in \on{LBA}({\bf id},{\bf id})$, such that $\psi_1 = \mu \circ (\psi\boxtimes
\on{id}_{{\bf id}}^{\on{LBA}})$ and $\psi_2 = \psi \circ \mu$.
\end{fact}

{\em Proof of Fact.} We decompose $\psi_i$ as
$\psi_i = \sum_{Z,W\in \on{Irr(Sch)}} \psi^i_{Z,W}$, where $\psi^i_{Z,W}$
belongs to the image of
$\on{LCA}({\bf id},Z) \otimes
\on{LCA}({\bf id},W) \otimes \on{LA}(Z\otimes W,{\bf id}) \to
\on{LBA}(T_2,{\bf id})$. Then $\psi^1_{Z,W}
\circ(\mu \boxtimes \on{id}_{{\bf id}}^{\on{LBA}})$ belongs to the image of
$\oplus_{Z_1,Z_2\in \on{Irr(Sch)}} \on{LCA}({\bf id},Z_1)\otimes
\on{LCA}({\bf id},Z_2) \otimes \on{LCA}({\bf id},W)
\otimes \on{LA}(Z_1\otimes Z_2 \otimes W,{\bf id}) \to \on{LBA}(T_3,{\bf id})$
while $\mu \circ (\psi_2\boxtimes \on{id}_{{\bf id}}^{\on{LBA}})$ belongs
to the image of
$$
\oplus_{Z_1,Z_2\in \on{Irr(Sch)}} \on{LCA}({\bf id},Z_1)\otimes
\on{LCA}({\bf id},Z_2) \otimes \on{LCA}({\bf id},{\bf id}) \otimes
\on{LA}(Z_1\otimes Z_2 \otimes {\bf id},{\bf id}) \to \on{LBA}(T_3,{\bf id}).
$$
Since this map is injective,
this implies that for $W\neq {\bf id}$, $(\sum_Z \psi^1_{Z,W})
\circ(\mu \boxtimes \on{id}_{{\bf id}}^{\on{LBA}})=0$. Fact \ref{lemma:inject}
then implies that for such $W$, $\sum_Z\psi^1_{Z,W}=0$, so
$\psi^1 = \sum_Z \psi^1_{Z,{\bf id}}$.

Let us denote by $std\in {\bf LA}(S\otimes {\bf id},{\bf id})$
the direct sum of all prop elements corresponding to $x_1...x_k \otimes
x\mapsto {1\over {k!}} \sum_{\sigma\in \SG_k}
[x_{\sigma(1)},...,[x_{\sigma(k)},x]]$. The PBW theorem for
free algebras implies that the map
$$
\on{LA}(Z,S) \to
\on{LA}(Z\otimes {\bf id},{\bf id}),
$$
$\lambda \mapsto std \circ (\lambda \boxtimes \on{id}_{\bf id}^{\on{LBA}})$
is a linear isomorphism. The decomposition of ${\bf LBA}(F,G)$ using
LA and LCA then implies that the map
\begin{equation} \label{map:std}
{\bf LBA}(A,S) \to {\bf LBA}(A\otimes {\bf id},{\bf id}), \quad
x\mapsto std \circ (x\boxtimes \on{id}_{\bf id}^{{\bf LBA}})
\end{equation}
is injective.

Then $\psi^1_{Z,{\bf id}} \in \on{LCA}({\bf id},Z) \otimes
\on{LA}(Z\otimes {\bf id},{\bf id}) \simeq \on{LCA}({\bf id},Z) \otimes
{\bf LA}(Z,S)$; we denote by $\psi^1_{Z,S}$ the image of $\psi^1_{Z,{\bf id}}$
in the latter space (and by $\psi^1_{Z,k}$ its $S^k$ component), which
may be viewed as the element of ${\bf LBA}({\bf id},S)$ such that
$std \circ (\psi^1_{Z,S} \boxtimes \on{id}_{\bf id}^{\on{LBA}})
= \psi^1_{Z,{\bf id}}$.

We now have
$$
std \circ ((\psi^1_{Z,S} \circ \mu )\boxtimes \on{id}_{\bf id}^{\on{LBA}})
= std \circ (\psi_2 \boxtimes \on{id}_{{\bf id}}^{\on{LBA}}).
$$
It follows from the injectivity of (\ref{map:std}) that
$\psi^1_{Z,k} \circ \mu=0$ for $k\neq 1$, and $\psi^1_{Z,1} \circ \mu
= \psi_2$. Now Fact \ref{lemma:inject} implies that $\psi^1_{Z,k} = 0$
for $k\neq 1$, so if we set $\psi := \psi^1_{Z,1}$, we have
$\psi_1 = \mu \circ (\psi \boxtimes
\on{id}_{\bf id}^{\on{LBA}})$ and $\psi_2 = \psi\circ \mu$.
\hfill \qed \medskip

According to Fact \ref{fact}, (\ref{main:id:theta}) implies that there exists
$f\in {\bf LBA}({\bf id},{\bf id})$, such that
\begin{equation} \label{theta:f}
\sigma' - \sigma' \circ (21) = f \circ \mu, \quad
\sigma' = \mu \circ (f \boxtimes \on{id}_{{\bf id}}^{{\bf LBA}}).
\end{equation}
These identities imply that $f$ satisfies the derivation
identity $f \circ \mu = \mu \circ (f\boxtimes \on{id}_{{\bf id}}^{{\bf LBA}} +
\on{id}_{{\bf id}}^{{\bf LBA}} \boxtimes f)$. We proved in \cite{univ:der}
that this implies that  $f$ is proportional to $\mu\circ \delta$.
Together with (\ref{theta:f}), this implies Lemma \ref{lemma:sigma'}.
\hfill \qed\medskip

It follows that $\sigma$ has the expansion $\mu \circ
(1 - \lambda (\mu \circ \delta) \boxtimes \on{id}_{\bf id}^{{\bf LBA}})$
+ terms of higher degree, where $\lambda$ is a scalar.

Since $\sigma \circ (\on{exp}(\lambda(\mu\circ\delta))\boxtimes
\on{id}_{\bf id}^{{\bf LBA}})$ satisfies the same identities as $\sigma$,
and has the expansion $\mu$ + terms of $\delta$-degree $\geq 2$,
it is equal to $\mu$, hence
$$
\sigma = \mu \circ (\on{exp}(\lambda(\mu\circ\delta))\boxtimes
\on{id}_{\bf id}^{{\bf LBA}}).
$$
One proves similarly that for some scalar $\lambda'$,
$$
\check\sigma = (\on{id}_{\bf id}^{{\bf LBA}}
\boxtimes \on{exp}(-\lambda'(\mu\circ\delta)))
\circ \delta.
$$
The Jacobi identity for the bracket $\hat D(\mu)$ then implies
that $\lambda = \lambda'$.

To give an idea of the proof, we write conditions for
the bracket $[-,-]_{\lambda,\lambda'}$ on ${\mathfrak D}(\A)$ (where $\A$
is a finite dimensional Lie bialgebra) to be Jacobi, where
$[-,-]_{\lambda,\lambda'}$ is given by the same formulas as above,
except for
$$
[x,\xi]_{\lambda,\lambda'} = \on{ad}^*_{e^{-\lambda D_\A}(x)}\xi
- \underline{\on{ad}}^*_{e^{\lambda' D_\A}(\xi)}x,
$$
where $D_\A$ is the derivation given by
$D_\A = \mu_\A\circ\delta_\A \oplus [-{}^t(\mu_\A \circ \delta_\A)]
\in \on{End}(\A)\oplus \on{End}(\A^*) \subset \on{End}({\mathfrak D}(\A))$.

The Jacobi identity for ${\mathfrak D}(\A)$ yields
$$
\underline{\on{ad}}^*_\xi([x_1,x_2]) =
[\underline{\on{ad}}^*_\xi x_1,x_2] + [x_1,\underline{\on{ad}}^*_\xi x_2]
+ \underline{\on{ad}}^*_{\on{ad}^*_{x_2}\xi}x_1
- \underline{\on{ad}}^*_{\on{ad}^*_{x_1}\xi}x_2
$$
for $\xi\on \A^*$, $x_1,x_2\in \A$.
On the other hand, the Jacobi identity for $[-,-]_{\lambda,\lambda'}$
implies that
$$
\underline{\on{ad}}^*_{\xi'}([x_1,x_2]) =
[\underline{\on{ad}}^*_{\xi'} x_1,x_2] + [x_1,\underline{\on{ad}}^*_{\xi'} x_2]
+ \underline{\on{ad}}^*_{\on{ad}^*_{e^{(\lambda' - \lambda)D_\A}(x_2)}\xi'}x_1
- \underline{\on{ad}}^*_{\on{ad}^*_{e^{(\lambda' - \lambda)D_\A}(x_1)}\xi'}x_2
$$
where $\xi' = e^{\lambda'D_\A}(\xi)$; replacing $\xi'$ by $\xi$,
the difference of these identities is
$$
\underline{\on{ad}}^*_{\on{ad}^*_{e^{(\lambda' - \lambda)D_\A}(x_2)}\xi}x_1
- \underline{\on{ad}}^*_{\on{ad}^*_{e^{(\lambda' - \lambda)D_\A}(x_1)}\xi}x_2
= \underline{\on{ad}}^*_{\on{ad}^*_{x_2}\xi}x_1
- \underline{\on{ad}}^*_{\on{ad}^*_{x_1}\xi}x_2,
$$
which translates as an identity in $\on{End}(\A^{\otimes 2})$
$$
(\on{id}_\A \otimes \mu_\A) \circ (\delta_\A \otimes \on{id}_\A) \circ
(\on{id}_\A \otimes (e^{(\lambda' - \lambda)D_\A}-1))
=
(\on{id}_\A \otimes \mu_\A) \circ (\delta_\A \otimes \on{id}_\A) \circ
(\on{id}_\A \otimes (e^{(\lambda' - \lambda)D_\A}-1)) \circ (21).
$$
One shows that at the propic level, the Jacobi identity for $\hat D(\mu)$
implies similarly
$
(\on{id}_{\bf id}^{{\bf LBA}} \boxtimes \mu) \circ (\delta \boxtimes
\on{id}_{\bf id}^{{\bf LBA}}) \circ
(\on{id}_{\bf id}^{{\bf LBA}} \boxtimes (e^{(\lambda' - \lambda)(\delta \circ \mu)}-1))
=
(\on{id}_{\bf id}^{{\bf LBA}} \boxtimes \mu) \circ
(\delta \boxtimes \on{id}_{\bf id}^{{\bf LBA}}) \circ
(\on{id}_{\bf id}^{{\bf LBA}} \boxtimes (e^{(\lambda' - \lambda)(\delta \circ \mu)}-1))
\circ (21).
$
The lowest degree term of this identity is the product of
$\lambda' - \lambda$ by
$$
\Big( (\on{id}_{\bf id}^{{\bf LBA}} \boxtimes \mu) \circ
(\delta \boxtimes \on{id}_{\bf id}^{{\bf LBA}}) \circ
(\on{id}_{\bf id}^{{\bf LBA}} \boxtimes (\delta \circ \mu))\Big)
\circ ((12) - (21));
$$
one checks that this is a nonzero element of ${\bf LBA}(T_2,T_2)$, so
$\lambda = \lambda'$.

We now prove that $\hat D = \on{Inn}(\on{id}_{\bf id}^{{\bf LBA}} \ul\boxtimes
\on{id}_{\bf 1}^{{\bf LBA}^*}
+ \on{id}_{\bf 1}^{{\bf LBA}} \ul\boxtimes *(Z_{\lambda}))
\circ \on{double}$. For this,
one has to check that
\begin{align*}
& \hat D(\mu) = \big( \on{id}_{\bf id}^{{\bf LBA}} \ul\boxtimes
\on{id}_{\bf 1}^{{\bf LBA}^*}
+ \on{id}_{\bf 1}^{{\bf LBA}} \ul\boxtimes *(Z_{\lambda}) \big)
\circ \on{double}(\mu) \\ & \circ
\Big(\on{id}_{\wedge^2}^{\bf {LBA}}\ul\boxtimes \on{id}_{\bf 1}^{{\bf LBA}^*}
+ \on{id}_{{\bf id}}^{{\bf LBA}} \ul\boxtimes *(Z_{\lambda})
+ \on{id}_{\bf 1}^{{\bf LBA}}\ul\boxtimes \wedge^2(*(Z_{\lambda}))
\Big)^{-1},
\end{align*}
\begin{align*}
& \hat D(\delta) = \Big(\on{id}_{\wedge^2}^{{\bf LBA}} \ul\boxtimes
\on{id}_{\bf 1}^{{\bf LBA}^*} +
\on{id}_{{\bf id}}^{{\bf LBA}} \ul\boxtimes *(Z_{\lambda})
+ \on{id}_{\bf 1}^{{\bf LBA}}\ul\boxtimes \wedge^2(*(Z_{\lambda}))
\Big) \\ & \circ \on{double}(\delta) \circ
\big( \on{id}_{\bf id}^{{\bf LBA}} \ul\boxtimes \on{id}_{\bf 1}^{{\bf LBA}^*}
+ \on{id}_{\bf 1}^{{\bf LBA}} \ul\boxtimes *(Z_{\lambda}) \big)^{-1}.
\end{align*}
The second identity, as well as the components $\wedge^2 \ul\boxtimes {\bf 1}
\to {\bf id}\ul\boxtimes {\bf 1}$ and ${\bf 1}\ul\boxtimes \wedge^2 \to
{\bf 1}\ul\boxtimes {\bf id}$ of the first identity, follow from the
fact that $Z_{\lambda}$ is central in ${\bf LBA}$, i.e.,
$\on{Inn}(Z_{\lambda})$ is the identity.

We now show that
\begin{equation} \label{id:1}
\mu_{{\bf id}\ul\boxtimes {\bf id},{\bf id}\ul\boxtimes {\bf 1}} =
\on{double}(\mu)_{{\bf id}\ul\boxtimes {\bf id},{\bf id}\ul\boxtimes {\bf 1}}
\circ (\on{id}_{\bf id}^{{\bf LBA}}\ul\boxtimes *(Z_{\lambda}^{-1})),
\end{equation}
\begin{equation} \label{id:2}
\mu_{{\bf id}\ul\boxtimes {\bf id},{\bf 1}\ul\boxtimes {\bf id}} =
(\on{id}_{\bf 1}^{{\bf LBA}} \ul\boxtimes *(Z_{\lambda})) \circ
\on{double}(\mu)_{{\bf id}\ul\boxtimes {\bf id},{\bf 1}\ul\boxtimes {\bf id}}
\circ (\on{id}_{\bf id}^{{\bf LBA}} \ul\boxtimes *(Z_{\lambda}^{-1})).
\end{equation}

Let us denote by $\xi\mapsto \tilde\xi$ the canonical map
${\bf LBA}_2(F\ul\boxtimes G,F'\ul\boxtimes G') \to
{\bf LBA}(F\otimes G^{\prime *},F'\otimes G^*)$. Then for
$x\in {\bf LBA}(G^*,G^*)$ and $y\in
{\bf LBA}(G^{\prime *},G^{\prime *})$, we have
$[(\on{id}_{F'} \ul\boxtimes *(y)) \circ
\xi \circ (\on{id}_F \ul\boxtimes *(x))]\tilde{} =(\on{id}_{F'}\boxtimes x)
\circ \tilde\xi \circ (\on{id}_F\boxtimes y)$.

Then $\tilde\mu_{{\bf id}\ul\boxtimes {\bf id},{\bf id}\ul\boxtimes {\bf 1}} =
\check\sigma = (\on{id}_{\bf id} \boxtimes Z_{-\lambda}) \circ \delta$;
on the other hand, $\{\on{double}(\mu)_{{\bf id}\ul\boxtimes {\bf id},
{\bf id}\ul\boxtimes {\bf 1}}\}\tilde{} = \delta$, so
$\{\on{double}(\mu)_{{\bf id}\ul\boxtimes {\bf id},
{\bf id}\ul\boxtimes {\bf 1}} \circ
(\on{id}_{\bf id}^{{\bf LBA}}\ul\boxtimes *(Z_{-\lambda}))\}
\tilde{} = (\on{id}_{\bf id}^{{\bf LBA}} \boxtimes Z_{-\lambda}) \circ \delta$,
which implies (\ref{id:1}).

Similarly, $\tilde\mu_{{\bf id}\ul\boxtimes {\bf id},{\bf 1}\ul\boxtimes
{\bf id}}
= -\sigma = -\mu \circ (Z_{-\lambda} \boxtimes \on{id}_{\bf id}^{{\bf LBA}})$;
on the other hand, $\{\on{double}(\mu)_{{\bf id}\ul\boxtimes {\bf id},
{\bf 1}\ul\boxtimes {\bf id}}\}\tilde{} = -\mu$, so
$\{(\on{id}_{\bf 1}^{{\bf LBA}}\ul\boxtimes *(Z_{\lambda})) \circ
\on{double}(\mu)_{{\bf id}\ul\boxtimes {\bf id},
{\bf 1}\ul\boxtimes {\bf id}} \circ (\on{id}_{\bf id}^{{\bf LBA}}\ul\boxtimes
*(Z_{-\lambda}))\}\tilde{} = Z_{-\lambda} \circ (-\mu) \circ
(\on{id}_{\bf id}^{{\bf LBA}} \boxtimes Z_\lambda)=
-\mu \circ (Z_{-\lambda} \boxtimes
\on{id}_{\bf id}^{{\bf LBA}})$, where we used the fact
that $Z_{-\lambda}$ is central in ${\bf LBA}$. This proves (\ref{id:2}).

\subsection{Compatibility of $Q$ with doubling operations}

We now summarize our results. We have
$$
D_{mult}(Q) \circ \on{Double} \circ Q^{-1} = \widetilde{\on{Double}} =
\on{Inn}(\Xi) \circ S(\varphi),
$$
$$
\varphi = \on{Inn}(\gamma^{-1}) \circ \hat D, \quad
\hat D = \on{Inn}(\on{id}_{\bf id}^{{\bf LBA}} \ul\boxtimes \on{id}_{\bf
1}^{{\bf LBA}^*} + \on{id}_{\bf 1}^{{\bf LBA}^*} \ul\boxtimes *(Z_\lambda))
\circ \on{double},
$$
which implies that
$$
D_{mult}(Q) \circ \on{Double} \circ Q^{-1} = \on{Inn}(\Lambda)\circ
S(\on{double}),
$$
where
$$
\Lambda = \Xi \circ (\gamma^{-1})_S \circ
(\on{id}_{\bf id}^{{\bf LBA}} \ul\boxtimes \on{id}_{\bf
1}^{{\bf LBA}^*} + \on{id}_{\bf 1}^{{\bf LBA}^*} \ul\boxtimes *(Z_\lambda))_S
\in S(D_{add}({\bf LBA}))({\bf id},{\bf id})^\times,
$$
which means that $Q$ is compatible with the doubling operations, as wanted.

\section{Applications}

\subsection{Finite dimensional Lie bialgebras}

Let us say that the QUE algebra $U$ is of finite type 
iff\footnote{$\on{Prim}(H)$ is the Lie algebra of primitive elements of 
a bialgebra $H$.} $\on{Prim}(U/\hbar U)$ is finite dimensional. 
Then a quantization functor gives rise to a category equivalence 
$\tilde Q : \{$finite dimensional Lie bialgebras$\} \to \{$QUE algebras
of finite type$\}$. 

We also have an action of $D_{4}$ on $\{$finite dimensional Lie bialgebras$\}$
and on $\{$QUE algebras of finite type$\}$. 
We have also self-maps of $\{$finite dimensional Lie bialgebras$\}$
and of $\{$QUE algebras of finite type$\}$, $\A\to {\mathfrak D}(\A)$ 
and $U\mapsto D(U)$. 

Then it follows from Theorem \ref{thm:duality} that if $Q$ is an 
Etingof-Kazhdan functor, then $Q(\A^{*cop})\simeq Q(\A)^{\star cop}$ 
for any finite dimensional 
Lie bialgebra $\A$, and from Theorem \ref{thm:doubles} that if $Q$ is an 
Etingof-Kazhdan functor, then $\tilde Q({\mathfrak D}(\A))\simeq 
D(\tilde Q(\A))$ for any finite dimensional Lie bialgebra $\A$. 
We recover in this way the result of \cite{EK1}.

\subsection{$\pm\NN$-graded Lie bialgebras}

Let us say that a $\pm\NN$-graded Lie bialgebra is of finite type 
iff each fixed degree component is finite dimensional. Let us 
say that a $\pm\NN$-graded QUE algebra $U = \oplus_{n\in\pm\NN}U_{n}$
is of finite type if: (a) $U_{0}$ is of finite type, and (b) each $U_{n}$ is 
finitely generated as a module over $U_{0}$. 

The permutation $*cop$ of $\{$finite dimensional Lie bialgebras$\}$
extends to a bijection $\{\NN$-graded Lie bialgebras of finite 
type$\}\to\{-\NN$-graded Lie bialgebras of finite type$\}$. The 
self-map ${\mathfrak D}$ of $\{$finite dimensional Lie bialgebras$\}$ 
extends to a map $\{\NN$-graded Lie bialgebras of finite 
type$\}\to\{\ZZ$-graded Lie bialgebras of finite type$\}$. 

Then the permutation $\star cop$ of $\{$QUE algebras of finite type$\}$
extends to a bijection $\star cop:\{\NN$-graded QUE algebras of finite 
type$\}\to\{-\NN$-graded QUE algebras of finite type$\}$. Moreover, the 
self-map $D$ of $\{$QUE algebras of finite type$\}$ extends to a map 
$D : \{\NN$-graded QUE algebras of finite type$\}\to\{\ZZ$-graded QUE 
algebras$\}$. 

A quantization functor $Q$ gives rise to functors $\tilde Q : \{\pm\NN$-graded 
Lie bialgebras of finite type$\}\to \{\pm\NN$-graded QUE algebras of 
finite type$\}$ and $\tilde Q : \{\ZZ$-graded 
bialgebras$\}\to\{\ZZ$-graded QUE algebras$\}$. 

Then Theorem \ref{thm:duality} implies that if $\A$ is a $\NN$-graded 
Lie bialgebra of finite type and $Q$ if an Etingof-Kazhdan functor, 
then $\tilde Q(\A^{*cop}) \simeq \tilde Q(\A)^{\star cop}$ and 
Theorem \ref{thm:doubles} implies that with 
the same assumptions, $\tilde Q({\mathfrak D}(\A))\simeq D(\tilde Q(\A))$. 

Theses result and those of the previous Subsection extend immediately 
to the super case. 

\subsection{Affine Lie superalgebras}\label{SS:g}   

In the rest of this section, $\kk = \CC$. Lie (bi)superalgebras, 
QUE superalgebras, etc., mean Lie (bi)superalgebras, etc., in the 
category of vector superspaces, i.e., the category whose objects 
are $\ZZ/2\ZZ$-graded vector spaces; 
morphisms are even linear maps; the symmetry constraint is given by the 
usual sign rule. 

Let $\g$ be an affine Lie superalgebra with symmetrized Cartan matrix 
$(A,\tau)$, where $A=(a_{ij})_{1\leq i,j\leq s}$ is a matrix 
with coefficients in $\CC$ and $\tau$ is a 
subset of $I:= \{1,...,s\}$ determining the parity of the generators 
(see \cite{vL}). Let $d_{1},\dots,d_{s}$ be nonzero rational numbers 
such that $d_{i}a_{ij}=d_{j}a_{ji}$ for $i,j\in I$.    

Let $\h$ be the Cartan sub-superalgebra of $\g$.  There exists linearly
 independent sets $\{\alpha_i\}_{i\in I}\subset \h^*$ and 
$\{h_i\}_{i\in I}\subset \h$ such that $\alpha_j(h_i)=a_{ij}$; up to
isomorphism these sets are determined by $A$.  Let $(-,-)$ be the 
non-degenerate supersymmetric invariant bilinear form on $\g$ defined in 
Proposition 4.2 of \cite{vL}.  The restriction of this form to the Cartan 
sub-superalgebra $\h$ is non-degenerate.  Furthermore, 
$(a,h_{i})=d_{i}^{-1}\alpha_{i}(a)$ for all $i\in I$ and $a\in \h$.

Let $\ghat$ be the Lie superalgebra presented by generators $h_{i}$, $e_{i}$, 
and $f_{i}$ for $i\in I$ (which are all even, except for the $e_i,f_i$
for $i\in \tau$, which are odd) and relations: 
\begin{align}
\label{R:LieSuperalg}
   [h_{i}, h_{j}]=0, \quad [h_{i}, e_{j}] =a_{ij}e_{j}, \quad 
   [ h_{i},f_{j}]=-a_{ij} f_{j}, \quad [e_{i},f_{j}]=\delta_{ij}h_{i}, 
   \quad i,j\in I. 
\end{align}
Then $\ghat$ has a unique maximal ideal $ \mathfrak{r}$ which intersects 
$\h$ trivially. We set $\g:=\ghat/\mathfrak{r}$. When $\g$ is not of type 
$A(n,n)^{(i)}$, Yamane (\cite{Yam}) wrote down explicit Serre-type
relations generating ${\mathfrak r}$. 

The Lie algebra $\g$ is equipped with a Lie bi-superalgebra structure, 
which we now describe. Let 
$\nilp_{+} \: (\text{resp., } \nilp_{-}) $ be the nilpotent Lie 
sub-superalgebra of $\g$ generated by $e_{i} \text{'s }$ 
$(\text{resp., } f_{i}\text{'s})$. Let $\borel_{\pm} := 
\nilp_{\pm}\oplus\h$ be the Borel Lie sub-superalgebras of $\g$.  
Let $\eta_{\pm}:\borel_{\pm} \rightarrow \g \oplus \h$ be defined by
$\eta_{\pm}(x)=x\oplus (\pm \bar x)$, where $\bar x$ is the image of 
$x$ in $\h$.  Using this embedding we can regard $\borel_{+}$ and 
$\borel_{-}$ as Lie sub-superalgebras of $\g \oplus \h$. 
Let $(-,-)_{\g \oplus \h}:=(-,-)\oplus [-(-,-)_{\h}]$, where $(-,-)_{\h}$ 
is the restriction of $(-,-)$ to $\h$.

\begin{prop}\label{P:glManinT}
$(\g \oplus \h, \borel_{+}, \borel_{-})$ is a super Manin triple 
(where $\g$ is equipped with $(-,-)_{\g \oplus \h}$).  
\end{prop}

\begin{proof}
Under the embedding $\eta_{\pm}$, the Lie subsuperalgebras $\borel_{\pm}$ 
are isotropic with respect to $(-,-)_{\g \oplus \h}$.  Since $(-,-) $ 
and $(-,-)_{\h}$ 
both are invariant super-symmetric nondegenerate bilinear forms, then 
so is $(-,-)_{\g \oplus \h}$.  Therefore the proposition follows.
\end{proof}

The proposition implies that $\g \oplus \h, \borel_{+}$ and 
$ \borel_{-}$ are Lie bi-superalgebras.  Moreover, we have that 
$\borel_{+}^{*}\cong\borel_{-}^{cop}$ as Lie bi-superalgebras, 
where $^*$ means the graded dual (for the principal grading on $\g$, given by 
$|a|=0$ for $a\in \h$ and $|e_i|=1$ for $i\in I$) and $^{cop}$ means 
changing the cobracket into its negative. The cobrackets of these Lie 
bi-superalgebras are given by the following formulas (see \cite{G3}):
$$ 
\delta(e_{i})=\frac{d_{i}}{2}(e_{i}\otimes h_{i}-h_{i}\otimes e_{i})
=\frac{d_{i}}{2}e_{i}\wedge h_{i}, \quad 
\delta(f_{i}) =-\frac{d_{i}}{2}f_{i}\wedge h_{i}, \quad \delta(a) =0 
$$ 
for any $i\in I$ and $a \in \h\subset \borel_{\pm}$.

Then the map $\borel_{+} \rightarrow \borel^{*}_+$ given by 
$e_i\mapsto e_i^*$ and $h_i\mapsto -\frac{2}{d_i}\sum_{j=1}^{s}a_{ji}h_j^*$
is an isomorphism of Lie bi-superalgebras (where $h_j^*,e_i^*\in\borel_+^*$
are defined by: $h_j^*$ has degree $0$, $h_j^*(h_i)=\delta_{ij}$ and 
$e_i^*$ has degree $-1$, $e_i^*(e_j)=\delta_{ij}$, and in $\borel_+^*$ the
$\NN$-grading is changed into its opposite).  This isomorphism 
restricts to the identification $\h \rightarrow \h^*$ coming from 
the form $-2(-,-)_\h$.

Define $\tilde\B_\pm\subset\tilde\g$ as the Lie sub-superalgebras 
generated by $h_i,e_i$ (resp., $h_i,f_i$), $i\in I$. One checks that 
the above formulas define Lie bialgebra structures on $\tilde\B_\pm$.
We have Lie bi-superalgebra morphisms $\tilde\B_\pm\to\B_\pm$.  

\subsection{Quantized affine Lie superalgebras}
In this subsection we use the previous results to show that 
Yamane's \cite{Yam} Drinfeld-Jimbo type superalgebras associated 
to certain affine Lie superalgebras are isomorphic to the corresponding 
Etingof-Kazhdan quantizations. The methods of this section are 
inspired by \cite{EK3}.

Let $Q_\Phi : \on{Bialg} \to S({\bf LBA})$ be an Etingof-Kazhdan functor.  
As above, $Q_\Phi$ gives rise to a functor $\tilde Q_\Phi : \{\pm\NN$-graded 
Lie bi-superalgebras of finite type, topologically free over 
$\CC[[\hbar]]\}\to \{\pm\NN$-graded QUE superalgebras of finite type$\}$; 
we then define $\EKh : \{\NN$-graded Lie bi-superalgebras of finite type over 
$\CC\}\to \{$QUE superalgebras over $\CC\}$, 
by $\EKh(\A,\mu_\A,\delta_\A):=\tilde Q_\Phi(\A[[\hbar]],
\mu_\A,\hbar\delta_\A)$. 

\begin{thm}\label{T:GRofUhbhat}
The quantized universal enveloping (QUE) superalgebra $\Uhbhat$ 
is isomorphic to the QUE superalgebra generated over 
$\CC[[\hbar]]$ by $\h$ and the elements $E_{i}$ for $i\in I$ (all 
generators are even except for $E_{i}$, $i\in \tau$ which are odd) 
satisfying the relations
\begin{align}
\label{pres:T}
    [a,a']=0, \quad [a,E_{i}]=\alpha_{i}(a)E_{i}, 
\end{align} 
with coproduct
\begin{align*}
\label{}
    \Delta({a})=a\otimes 1 + 1 \otimes a, \quad 
    \Delta({E_{i}})=E_{i}\otimes q^{d_{i}h_{i}}+ 1\otimes E_{i},
\end{align*}   
for all $a,a'\in \h$ and $i,j\in I$. (Here $q = e^{\hbar/2}$.) 
\end{thm}

The theorem follows from the following two lemmas.  

\begin{lemma}\label{L:UQESbhat}
The QUE superalgebra $\Uhbhat$ is isomorphic to the QUE 
superalgebra generated over $\CC[[\hbar]]$ by $\h$ and the elements 
$E_{i}, \: i\in I$ (all generators are even except for 
$E_{i}$, $i\in \tau$ which are odd) satisfying the relations
\begin{equation}
\label{E:relUhbhatCar}
[a,a']=0, \quad  
[a,E_{i}]=\alpha_{i}(a)E_{i},
\end{equation}  
with coproduct
\begin{equation}
\label{E:coprd}
 \Delta({a})=a\otimes 1 + 1 \otimes a, \quad  
  \Delta({E_{i}})=E_{i}\otimes q^{\gamma_{i}}+ 1\otimes E_{i},
\end{equation}   
for all $a,a'\in \h$ and $i,j\in I$ and suitable elements 
$\gamma_{i}\in \h[[\hbar]]$.
\end{lemma}
\begin{proof}
An analogous lemma was proved in \cite{G3} for finite dimensional 
Lie superalgebras of type A-G.  The proof in \cite{G3} does not 
use the finite dimensional condition.   Thus, the proof of 
Lemma \ref{T:GRofUhbhat} follows word for word as in the 
analogous lemma of \cite{G3}.   
\end{proof}

\begin{lemma}\label{L:gamma}
$ \gamma_{i}=d_{i}h_{i}$. 
\end{lemma}
\begin{proof} 
By the definition of $\borel_+$, we have a (surjective) Lie bi-superalgebra 
morphism $\pi : \borelhat_{+}\rightarrow\borel_{+}$. Since $\EKh$ is a 
functor, this induces a morphism of Hopf superalgebras 
$\EKh(\borelhat_{+}) \rightarrow \EKh(\borel_{+})$.  
This morphism is onto as it is onto modulo $\hbar$.   
Therefore, $\EKh(\borel_{+})$ is generated by $\h$ and $E_{i}$ 
satisfying the relations \eqref{E:relUhbhatCar}-\eqref{E:coprd} 
(and possibly other relations).  So it suffices to show that the images of 
$\gamma_i,h_i$ by $\EKh(\pi)$ (again denoted $\gamma_{i},h_{i}$) satisfy 
$ \gamma_{i}=d_{i}h_{i}$ in $\EKh(\borel_{+})$.

As mentioned above $\borelp$ is self dual in the graded sense,
i.e., $\borelp\cong \borel_+^{*}$ (as graded Lie bi-superalgebras, where 
one of the gradings is changed to its negative). Since $\EKh$ is a 
functor, we get a $\NN$-graded isomorphism  
\begin{equation}
\label{E:quantb}
\EKh(\borelp)\cong \EKh(\borel_+^{*})
\end{equation} 
(again one of the gradings is changed to its negative). 
It follows from the discussion in Subsection \ref{SS:g} that the 
restriction of this isomorphism to $\EKh(\h)$ comes from the 
identification $\h\rightarrow \h^*$ using the form $-2(-,-)_\h$ on $\h$.   

Using Theorem \ref{thm:duality} with $\theta=*cop$, we have 
\begin{align}
\label{E:compdual2}
 \EKh(\borel_+^*) &\cong \EKh(\borel_+^{cop})\qdo.
\end{align}  
Moreover, this theorem says that if $V$ is a finite dimensional 
commutative and cocommutative Lie bialgebra, then the isomorphism 
$\EKh(V^*)\simeq \EKh(V)^\star$ is induced by the bialgebra 
pairing $\EKh(V)\otimes \EKh(V^*) \simeq S(V)[[\hbar]]\otimes S(V^*)[[\hbar]]
\to \CC((\hbar))$ given by $V \otimes V^* \in \xi\otimes x\mapsto 
\hbar^{-1}\xi(x)$. 

Combining Equations \eqref{E:quantb} and \eqref{E:compdual2} with 
$\EKh(\borel_+^{cop})\cong \EKm(\borel_+)$  we have
$$\EKh(\borelp)\cong \EKm(\borelp)\qdo.$$ 
This isomorphism gives rise to a nondegenerate bilinear form 
$B: \EKh(\borelp) \otimes \EKm(\borelp) \rightarrow \CC((\hbar))$ 
which satisfies the following conditions 
\begin{eqnarray*}
  B(xy,z)= B(x \otimes y,\Delta(z)), & B(x, yz)= B(\Delta(x), y \otimes z) 
\end{eqnarray*}
$$ B(q^{a},q^{b})=q^{-(a,b)}, $$
where $a,b \in \h$ and $x,y,z\in\EKh(\borelp)$ (the second line follows from 
the properties of the pairing $\EKh(\h)\otimes \EKh(\h^*)\to\CC((\hbar))$
recalled above). 

Let $a\in\h$ and $i\in I$.  Set $B_{i}=B(E_{i},E_{i})$, which is nonzero.   
Using the above properties of $B$ we have 
\begin{align}
\label{E:GammaB}
    B_{i}q^{(a,\gamma_{i})}=&B(E_{i},q^{a}E_{i}q^{-a}).  
\end{align}  
From relation $[a,E_{i}]=\alpha_{i}(a)E_{i}$ it follows that 
$q^{h_{j}}E_{i}q^{-h_{j}}=q^{\alpha_{i}(h_{j})}E_{i}$ and so 
\begin{align}
\label{E:quantea}
   q^{a}E_{i}q^{-a}=&q^{\alpha_{i}(a)}E_{i}.  
\end{align}
Combining (\ref{E:GammaB}) and (\ref{E:quantea}) we have 
$$B_{i}q^{(a,\gamma_{i})}=B(E_{i},q^{a}E_{i}q^{-a})
=B(E_{i},q^{\alpha_{i}(a)}E_{i})=B_{i}q^{\alpha_{i}(a)}.$$
Hence, $(a, \gamma_{i})=\alpha_{i}(a)$, but $ \alpha_{i}(a)= d_{i}(a,h_{i})$, 
and so $\gamma_{i}=d_{i}h_{i}$, which completes the proof. 
\end{proof}

\begin{lemma}
Let $I_\hbar:= \on{Ker}(\EKh(\pi):\EKh(\tilde\B_+)\to \EKh(\B_+))$. 
Then $I_\hbar\subset \EKh(\tilde\B_+)$ is maximal with the following 
properties: is a graded bialgebra ideal,
$\hbar$-divisible (i.e., the quotient by it is $\hbar$-torsion free),  
whose components of degree $0$ and $1$ are zero. 
\end{lemma}

{\em Proof.} This is equivalent to saying that the 
maximal ideal of $\EKh(\B_+)$  with the same properties is $0$. Indeed,
such an ideal $J_\hbar\subset \EKh(\B_+)$ gives rise to a QUE superalgebra
$\EKh(\B_+)/J_\hbar$ (because this is $\hbar$-torsion free, and by the 
Milnor-Moore theorem any quotient of an enveloping superalgebra is again 
an enveloping superalgebra) and to a QUE superalgebra morphism $\EKh(\B_+)
\to \EKh(\B_+)/J_\hbar$ (because a bi-superalgebra morphism between 
universal enveloping superalgebras is induced by a morphism of Lie 
superalgebras). The classical limit of this morphism is a surjective morphism 
$\B_+\to \x$ of $\NN$-graded Lie superalgebras, which is the identity in 
degrees $0$ and $1$. Let $\iii:= \on{Ker}(\B_+\to \x)$ is a Lie ideal and 
coideal. Then $\iii^\perp\subset \B_+^* \simeq \B_-$ is a Lie sub-superalgebra, 
containing $\h$ and the $e_i^*$. Since $\B_-$ is generated by these elements, 
$\iii=0$, so $\B_+\to\x$ is the identity of $\B_+$, so $J_\hbar=0$. 
\hfill\qed \medskip

\begin{thm}\label{T:FormB}
There exists a unique bilinear form on $\Uhbhat$, taking 
values in $\CC((\hbar))$ and with the following properties
\begin{align*}
\label{}
  C(xy,z)=& C(x \otimes y,\Delta(z)), & C(x, yz)=& C(\Delta(x), y \otimes z), \\
  C(q^{a},q^{b})=&q^{-(a,b)}, a,b \in \h, & C(E_{i},E_{j})
  =&\delta_{ij}B(E_{i},E_{i}),
\end{align*}
where $B$ is given in the proof of Lemma \ref{L:gamma}.   
Moreover $\Uhb \cong \Uhbhat/\on{Ker}(C)$ as QUE superalgebras.   
\end{thm}
\begin{proof}
The existence and uniqueness follows from the fact that the 
superalgebra generated by the $E_{i}$ is free (see Theorem \ref{T:GRofUhbhat}).  

We will show that there is a nondegenerate bilinear form on $\Uhb$ with the 
same properties as $B$. 

We first construct an isomorphism $\EKm(\borel_{+})^{cop} \rightarrow \Uhb$ 
of $\NN$-graded QUE algebras. 
Let $T$ be the superalgebra presented by (\ref{pres:T}). We have
an isomorphism $\EKh(\tilde\B_+)\simeq T$ of $\NN$-graded superalgebras. 
Let $\Delta_\hbar : T\to T^{\otimes 2}$ be the coproduct on $T$ 
obtained from the coproduct of $\EKh(\tilde\B_+)$. Let $c\in \on{Aut}(T)$ 
be the conjugation by $q^{\sum x_{i}^{2}/2}$, where $x_{i}$ is a 
orthonormal basis for $\h$. Then $\Delta_{-\hbar}^{2,1} = (c\otimes c)\circ 
\Delta_\hbar \circ c^{-1}$. Let $I_\hbar := \on{Ker}(T\to \EKh(\B_+))$. 
We have $\Delta_\hbar(I_\hbar)\subset I_\hbar \otimes T + T \otimes I_\hbar$. 
Therefore $\Delta_{-\hbar}(c(I_\hbar)) \subset c(I_\hbar)\otimes T 
+ T \otimes c(I_\hbar)$, 
so by the maximality of $I_{-\hbar}$ w.r.t. $\Delta_{-\hbar}$, $c(I_\hbar)
\subset I_{-\hbar}$. 
Replacing $\hbar$ by $-\hbar$, we get similarly $c^{-1}(I_{-\hbar})\subset
I_\hbar$, so $c(I_\hbar)=I_{-\hbar}$. It follows that $c$ induces an 
isomorphism $T/I_\hbar \to T/I_{-\hbar}$, intertwining the coproduct maps 
$T/I_{\pm \hbar}\to (T/I_{\pm \hbar})^{\otimes 2}$ induced by $\Delta_\hbar$ 
and $\Delta_{-\hbar}^{2,1}$. 
This is an isomorphism $\EKm(\borel_{+})^{cop} \rightarrow \Uhb$. 
  
From the proof of Lemma \ref{L:gamma}, we have 
that $\EKh(\borel_{+}) \cong \EKm(\borel_{+})\qdo$. 
Combining this with the above isomorphism 
$\EKm(\borel_{+})^{cop} \rightarrow \Uhb$, we obtain 
an isomorphism $\Uhb \cong \Uhb\qd$.
This isomorphism gives rise to the desired form on $\Uhb$.  

It follows that $C$ is the pull-back of the bilinear form $B$ 
constructed in Lemma \ref{L:gamma}. Thus the image of the kernel 
of $C$ is contained in the kernel of $B$ under the natural projection.
But the kernel of $B$ is zero since the form is nondegenerate.   
Thus we have $\Uhb \cong \Uhbhat/\on{Ker}(C)$. 
\end{proof}

In \cite{Yam}, Yamane introduced what we will call here the 
Drinfeld-Jimbo type quantization $\DJh(\g)$ of $\g$ as follows: 
$\DJh(\B_\pm):= \EKh(\tilde\B_\pm)/\on{Ker}(C_\pm)$ (where $C_+=C$ and $C_-$
is its analogue for $\EKh(\tilde\B_-)$); we have a non-degenerate
pairing $\DJh(\B_+)\otimes\DJh(\B_-)\to\CC((\hbar))$, which allows 
to identify $\DJh(\B_-) = \DJh(\B_+)^{\star cop}$; then 
$\DJh(\g):= D(\DJh(\B_+))/(\h\simeq \h^*)$. The superalgebra $\DJh(\g)$ is a 
braided (i.e., quasitriangular) quantized universal enveloping superalgebra 
where the braiding is given by the universal $R$-matrix. 
An immediate corollary of the above theorem is: 

\begin{cor}\label{C:QUEBorelIso}
$\EKh(\B_\pm)\simeq \DJh(\B_\pm)$ as QUE superalgebras. 
\end{cor}

We also have: 

\begin{thm}\label{T:UhIsoDJ}
There exists an 
isomorphism of braided quantized universal enveloping (QUE) superalgebras:
$$\alpha: \DJh(\g) \rightarrow \EKh(\g)$$
such that $\alpha_{|\h}={\rm id}_{\h}$.
\end{thm}
\begin{proof}
We first construct an isomorphism between these QUE superalgebras which is the 
identity on $\h$ (disregarding the braiding). 
We have $\DJh(\g) \simeq D(\DJh(\B_{+}))/(\h\simeq\h^{*}) = 
\DJh(\B_{+})\otimes\DJh(\B_{+})^{\star cop}/(\h\simeq\h^{*})$. 
Corollary \ref{C:QUEBorelIso} yields an isomorphism $\DJh(\B_{+})\simeq
\EKh(\B_{+})$ inducing the identity on $\h$, hence an isomorphism 
$D(\DJh(\B_{+}))/(\h\simeq\h^{*}) \simeq D(\EKh(\B_{+}))/(\h\simeq\h^{*})$
inducing the identity on $\h$. Since quantization commutes with the double, 
we have an isomorphism $D(\EKh(\B_{+})) \simeq \EKh({\mathfrak D}(\B_{+}))$. 
The fact that this isomorphism induces the identity on $\h\oplus\h^{*}$
is a consequence of the following property of the isomorphism 
$D(\EKh(\A))\simeq \EKh({\mathfrak D}(\A))$: 

\begin{lemma}
For $\NN$-graded Lie bi(super)algebra morphisms $\B\stackrel{i}{\to}\A$ and $\A
\stackrel{p}{\to}\B$ 
such that $\B\to\A\to\B$ is the identity, we get a morphism 
${\mathfrak D}(\B)\to{\mathfrak D}(\A)$; in the same way, if $A,B$ are 
$\NN$-graded QUE algebra of finite type, QUE algebra morphisms $B\to A$ 
and $A\to B$ such that $B\to A \to B$ is the identity give rise to a 
morphism $D(B)\to D(A)$; taking $(B\to A):= \EKh(\B\to\A)$, 
$(A\to B):= \EKh(\A\to\B)$, we get a morphism $D(\EKh(\B))\to D(\EKh(\A))$; 
the diagram 
$\begin{matrix} \EKh({\mathfrak D}(\B))& \to & \EKh({\mathfrak D}(\A)) \\
\downarrow & & \downarrow\\
D(\EKh(\B))  & \to & D(\EKh(\A))\end{matrix}$  
commutes.
\end{lemma}

{\em Proof of Lemma.} We have $\EKh(\A) = (S(\A)[[\hbar]],m(\mu_\A,\hbar\delta_\A),
\Delta(\mu_\A,\hbar\delta_\A))$. The isomorphism 
$\EKh(\A) \otimes \EKh(\A)^{\star cop} \simeq 
D(\EKh(\A))\simeq \EKh({\mathfrak D}(\A))$ is induced by a map 
$\lambda^{-1}(\mu_\A,\hbar\delta_\A) : S(\A)\otimes S(\A^*)[[\hbar]]
\to S(\A\oplus\A^*)[[\hbar]]$, where $\lambda\in{\bf LBA}_2(S\ul\boxtimes S,
S\ul\boxtimes S)^\times$. We thus have to show that 
$\lambda^{-1}(\mu_\A,\hbar\delta_\A)\circ S(i\oplus p^t) = 
S(i\oplus p^t)\circ \lambda^{-1}(\mu_\B,\hbar\delta_\B)$
(equality of maps $S(\B\oplus\B^*)[[\hbar]]\to S(\A\oplus\A^*)[[\hbar]]$). 
Let $\tilde\lambda\in {\bf LBA}(S^{\otimes 2},S^{\otimes 2})$ be the 
element corresponding to $\lambda^{-1}$, then we should show that 
$(\on{id}_{S(\A)}\otimes S(p)) \circ \tilde\lambda(\mu_\A,\hbar\delta_\A) 
\circ (S(i)\otimes \on{id}_{S(\A)}) = 
(S(i)\otimes \on{id}_{S(\B)})\circ \tilde\lambda(\mu_\B,\hbar\delta_\B) 
\circ (\on{id}_{S(\B)}\otimes S(p))$ (equality of maps
$S(\B\oplus\A)[[\hbar]] \to S(\A\oplus\B)[[\hbar]]$). 

$\tilde\lambda$ is a sum $\sum_{Z_{1,2,3}\in\on{Irr(Sch)}} \sum_\alpha
(\lambda_1^\alpha \boxtimes \lambda_2^\alpha) \circ (\kappa_1^\alpha
\boxtimes \kappa_2^\alpha)$, where $\sum_\alpha \kappa_1^\alpha \otimes
\kappa_2^\alpha \otimes \lambda_1^\alpha\otimes \lambda_2^\alpha
\in \on{LCA}(S,Z_1\otimes Z_2) \otimes 
\on{LCA}(S,Z_3) \otimes \on{LA}(Z_1,S)\otimes \on{LA}(Z_2\otimes Z_3,S)$. 

Then 
\begin{align*}
& (\on{id}_{S(\A)}\otimes S(p)) \circ \tilde\lambda(\mu_\A,\hbar\delta_\A) 
\circ (S(i)\otimes \on{id}_{S(\A)})
\\ & 
= (\on{id}_{S(\A)}\otimes S(p)) \circ \Big(\sum_{Z_{1,2,3},\alpha}
(\lambda_1^\alpha(\mu_\A,\hbar\delta_\A) \otimes \lambda_2^\alpha(\mu_\A,\hbar\delta_\A))
\circ (\kappa_1^\alpha(\mu_\A,\hbar\delta_\A) \otimes 
\kappa_2^\alpha(\mu_\A,\hbar\delta_\A))
\Big) \circ (S(i)\otimes \on{id}_{S(\A)})
\\ & 
= 
\sum_{Z_{1,2,3},\alpha}
(\lambda_1^\alpha(\mu_\A,\hbar\delta_\A) \otimes
\lambda_2^\alpha(\mu_\B,\hbar\delta_\B))
\circ (\on{id}_{Z_1(\A)} \otimes Z_2(p)\otimes Z_3(p))
\circ (Z_1(i)\otimes Z_2(i)\otimes \on{id}_{Z_3}(i))
\\ & \circ (\kappa_1^\alpha(\mu_\B,\hbar\delta_\B) 
\otimes \kappa_2^\alpha(\mu_\A,\hbar\delta_\A))
\\ & 
= 
\sum_{Z_{1,2,3},\alpha}
(\lambda_1^\alpha(\mu_\A,\hbar\delta_\A) \otimes
\lambda_2^\alpha(\mu_\B,\hbar\delta_\B))
\circ (Z_1(i) \otimes \on{id}_{Z_2(\B)}\otimes Z_3(p))  
\circ (\kappa_1^\alpha(\mu_\B,\hbar\delta_\B) 
\otimes \kappa_2^\alpha(\mu_\A,\hbar\delta_\A))
\\ & 
= 
(S(i)\otimes \on{id}_{S(\B)}) \circ 
\Big( \sum_{Z_{1,2,3},\alpha}
(\lambda_1^\alpha(\mu_\B,\hbar\delta_\B) \otimes
\lambda_2^\alpha(\mu_\B,\hbar\delta_\B))
\circ (\kappa_1^\alpha(\mu_\B,\hbar\delta_\B) 
\otimes \kappa_2^\alpha(\mu_\B,\hbar\delta_\B)) \Big) 
\circ (\on{id}_{S(\A)}\otimes S(p))
\\ & 
= 
(S(i)\otimes \on{id}_{S(\B)}) \circ 
\tilde\lambda(\mu_\B,\hbar\delta_\B)
\circ (\on{id}_{S(\A)}\otimes S(p)),   
\end{align*}
as wanted; the second and fourth equalities follow from the fact that 
$p$ and $i$ are 
Lie bialgebra morphisms, the third one follows from $p\circ i=\on{id}_\B$. 
\hfill \qed \medskip 

We thus get an isomorphism $D(\EKh(\B_{+}))/(\h\simeq\h^{*})
\simeq \EKh({\mathfrak D}(\B_{+}))/(\h\simeq\h^{*}) = \EKh(\g)$
inducing the identity on $\h$. We thus have an QUE isomorphism 
$\DJh(\g)\simeq \EKh(\g)$, inducing the isomorphism on $\h$. 

This isomorphism is actually an isomorphism of braided QUE superalgebras.  
This follows from the fact that the $R$-matrix is uniquely determined by its 
degree 1 part and its relations with the comultiplication of the Hopf 
superalgebra. For a proof of this fact see \cite{KhT}.  
\end{proof}

When $\g$ is not of type $A(n,n)^{(i)}$, Yamane wrote down 
generators for $\on{Ker}(C_\pm)$ in terms of 
the $E_i,F_i$. This gives explicit presentations of $\DJh(\B_\pm)$
and $\DJh(\g)$. He also gave an explicit formula for the $R$-matrix. 
The above corollaries then imply that these QUE (braided) superalgebras
identify with $\EKh(\B_\pm)$, $\EKh(\g)$.  

\subsection{Category $\cato$}
In this subsection we show that the category $\cato$ for $\g$ can be 
deformed to an analogous category over $\DJh(\g)$.  We then give a 
Drinfeld-Kohno type theorem for affine Lie superalgebras.  This subsection 
follows \cite{EK3} where similar results are proved for Kac-Moody Lie algebras.
 
Let $\cato[[\hbar]]$ be the deformed category $\cato$ for $\g$, i.e.,
representations of $\g$ on topologically free $\CC[[\hbar]]$-modules 
which are $\h$-diagonalizable, have finite dimensional weight spaces 
and whose weights belong to a union of finitely many cones 
$\lambda -\sum\NN\alpha_i$, $\lambda \in \h^*[[\hbar]]$.  

In a similar way one defines $\cato_\hbar$ to be the category of topologically 
free $\DJh(\g)$-modules whose weights satisfy the above conditions.  
The category $\cato_\hbar$ is a braided tensor category where the braiding 
comes from the universal $R$-matrix of $\DJh(\g)$.

Let $\Omega\in \g^{\hat\otimes 2}$ be the inverse element corresponding 
to the bilinear form $(-,-)$ ($\g^{\hat\otimes 2}=\hat\oplus_{a,b}
\g_a\otimes\g_b$, where $\g$ decomposes as $\oplus_a \g_a$ for the principal 
grading). The action of $\Omega$ is well-defined on 
the tensor product of two modules from the category $\cato[[\hbar]]$.  
Following Drinfeld, we give $\cato[[\hbar]]$ a braided tensor structure 
with associator $\Phi(\hbar\Omega_{12},\hbar\Omega_{23})$ 
and braiding $q^{\Omega}$.  

\begin{thm}\label{T:cato}
There exists a functor $F:\cato[[\hbar]]\rightarrow \cato_\hbar$ which is an 
equivalence of braided tensor categories and commuting with the forgetful 
functors $\cato[[\hbar]],\cato_\hbar \to \{\h$-graded 
$\CC[[\hbar]]$-modules$\}$. 
\end{thm}   
\begin{proof}
Theorem \ref{T:UhIsoDJ} allows us to replace $\DJh(\g)$ with 
$D(\Uhb)/(\h =\h^*)$.  Therefore, it is enough to construct a functor 
between the corresponding categories for the superalgebras 
${\mathfrak D}(\borel_+)$ and $D(\Uhb)$.

In the non-super case such a functor was defined in 
Theorem 4.1 of \cite{EK3}.  Since Theorem 4.1 of \cite{EK3} is formulated 
in the general categorical language of \cite{EK2} it is easy to check 
it generalizes to the super setting and so the theorem follows.

\end{proof}

Let $V$ be a module of $\cato[[\hbar]]$ and let $V_q$ be its image under $F$ 
in $\cato_\hbar$. Let $B_{n}=\langle\sigma_{i}\rangle $ be the braid group.  
Define $\rho_{n}$ to be the representation of $B_{n}$ on $V_q^{\otimes n}$ 
given by $$\sigma_{i}\mapsto \tau_{i, i+1}R_{i i+1}$$ where $\tau_{i,i+1}$ 
is the super permutation of the $i$th and the $(i+1)$th component and 
$R$ is the universal $R$-matrix of $\DJh(\g)$.  

On the other hand, consider the Knizhnik-Zamolodchikov system of differential 
equations with respect to a function $\omega(z_{1},...,z_{n}) $ of complex 
variables $z_{1},...,z_{n}$ with values in $V^{\otimes n}$:
\begin{equation}
\label{E:KZvalueV}
\frac{\partial \omega}{\partial z_{i}} = \frac{\hbar}{2\pi i}
\sum_{i\neq j}\frac{\Omega_{ij}\omega}{z_{i}-z_{j}}.
\end{equation}

This system of equations defines a flat connection on the trivial 
bundle $Y_{n}\times V^{\otimes n}$, where $Y_n=\{(z_1,...,z_n) | i\neq j 
\text{ implies } z_i \neq z_j \}\subset \CC^{n}$. This connection determines 
a monodromy representation from $\pi_{1}(Y_{n})$ to 
$\Aut_{\CC[[\hbar]]}(V^{\otimes n})$ (as it is a direct sum of connections in the
weight spaces, which are finite dimensional).  Moreover, since the system 
of equations \eqref{E:KZvalueV} is invariant under the action of the symmetric 
group we obtain a monodromy representation 
$$\rho_{n}^{KZ}: \pi_{1}(X_{n},p)\rightarrow \Aut_{\CC[[\hbar]]}(V^{\otimes n})$$
where $X_{n}=Y_{n}/S_{n}$ and $p=(1,2,...,n)\in \CC^{n}$.  Finally, we 
identify $\pi_{1}(X_{n},p)$ with the braid group $B_{n}$ to get a monodromy 
representation of $B_{n}$.

\begin{thm} (The Drinfeld-Kohno theorem for affine Lie superalgebras) 
\label{T:D-KT} The representations $\rho_{n}$ and $\rho_{n}^{KZ}$ are 
equivalent.
\end{thm}
\begin{proof}
Let $\rho_{n}^{R_{KZ}}$ be the representation of $B_{n}$ on $V^{\otimes n}$ 
induced by the $R$-matrix $R_{KZ}=e^{\hbar\Omega/2}$ in the category
$\cato[[\hbar]]$ 
with the Knizhnik-Zamolodchikov associator $\Phi_{KZ}$. From Theorem 
\ref{T:cato} we have that $\rho_{n}^{R_{KZ}}$ and $\rho_{n}$ correspond 
to each other under the braided tensor functor $F$.  The theorem follows 
since $\rho^{KZ}_{n}$ coincides with $\rho_{n}^{R_{KZ}}$ when $\Phi=\Phi_{KZ}$. 
\end{proof}

\subsection*{Acknowledgements} 
We would like to thank P. Etingof, D. Kazhdan and H. Yamane for useful 
discussions. The work of N. Geer was supported by NSF grant no. DMS-0706725. 
He also thanks ULP and IRMA (CNRS, Strasbourg) for the invitations during 
which this work was done.

\end{document}